# Random Width and Brightness

*Polyhedral Density Theory, Reconstruction, and Gaussian Identifiability*

Omri Abas
*Independent researcher, Kadima-Zoran, Israel*
omri.abas1@gmail.com
ORCID: 0009-0003-4857-0790
10 August 2026

## Abstract

Let $U$ be uniformly distributed on the unit sphere. We develop a self-contained forward and inverse theory for the random width $w_K(U)$ and brightness $b_K(U)$ of three-dimensional convex bodies. For every full-dimensional polytope, a global spherical co-area formula expresses the width density as a finite sum of angular apertures determined by the normal fan of its difference body; in particular, the density is piecewise real analytic with a finite geometrically determined critical set. This theory yields exact densities for the width of the regular tetrahedron, resolving a question of Finch, and for the regular truncated octahedron, together with the tetrahedral brightness law and the equivalent rhombic-dodecahedral width law. On the inverse side, second- and third-order polarized cosine-transform moments reconstruct finite labelled direction systems whenever the observed triangles span the cycle space of the correlation graph; signed-graph switching describes the unavoidable ambiguity. In contrast, equal three-dimensional intrinsic volumes do not determine either the width law or the brightness law, even for centrally symmetric bodies. Removing the spatial rank constraint gives a dimension-free identifiability theorem for centered multivariate folded-normal vectors: pairwise absolute moments and an anchored family of triple absolute moments, comprising $(m-1)^2$ labelled observations for a complete correlation graph, determine the correlation matrix up to diagonal sign conjugacy without fourth-order moments. A harmonic decomposition further identifies the degree-two variance contribution as a constant multiple of the squared Frobenius norm of the traceless part of the weighted frame operator and explains why this contribution vanishes under irreducible symmetry.

## 1. Introduction

### *1.1. Directional laws and geometric information*

Let $K \subset \mathbb{R}^3$ be a convex body and let $U$ be uniformly distributed on the unit sphere. The random variables $W_K = w_K(U)$ and $B_K = b_K(U)$ record, respectively, the width of $K$ and the area of its orthogonal projection in a random direction. Their expectations are classical intrinsic quantities. Their full distributions are substantially finer: they retain the arrangement of directional values over the sphere, the singularities created when level sets cross changes in the normal structure, and geometric information that is invisible to finitely many rotation-invariant valuations.

The central problem is to determine what survives the compression from directional geometry to a one-dimensional probability law. On the forward side, one asks for an exact density, its support, analytic branches and transition singularities; on the inverse side, one asks which geometric data can be reconstructed from directional transforms, polarized moments, or the scalar law itself. These questions lead to spherical co-area and normal fans, zonoids and cosine transforms, harmonic analysis, signed graphs, and Gaussian absolute moments.

For an arbitrary convex body, width is the support function of its difference body: $w_K(u)=h_{K-K}(u)$. If $\nu$ is a finite measure on $S^2$, write $A_\nu(u)=\int_{S^2}|u\cdot n|d\nu(n)$. Widths of zonoids and brightness functions are both instances of this spherical cosine transform: $w_Z=A_{\mu_Z}$ and $b_K=A_{S_K/2}$, where $\mu_Z$ is the generating measure of $Z$ and $S_K$ is the surface-area measure of $K$; see Gardner [9] and Schneider [13].

This common representation has an exact geometric boundary. Proposition 2.6 proves that $w_K(u)=\sum_{i=1}^{m} c_i|u\cdot n_i|$ for all $u\in S^2$, with $c_i>0$, if and only if $K-K$ is a zonotope. Consequently, the normal-fan method applies to every three-dimensional polytope width, whereas the polarized absolute-product calculus applies on the zonotopal side and, through projection bodies, to polyhedral brightness.

### 1.2. A density theorem for arbitrary three-dimensional polytopes

The first principal result is a global density formula requiring neither central symmetry nor a zonotopal representation. Let $P\subset\mathbb{R}^3$ be a full-dimensional polytope and put $Q=P-P$. On the spherical normal cell $C_q$ corresponding to a vertex $q$ of $Q$, one has $w_P(u)=q\cdot u$.

The level set $w_P(u)=t$ is therefore the intersection of $C_q$ with a small circle. Its contribution to the density is controlled by the angular aperture $\Theta_q(t)$ of the portion lying inside the cell. Applying the spherical co-area formula in this piecewise-linear setting yields Theorem 3.4:

$$f_{W_P}(t)=\frac{1}{4\pi}\sum_{\substack{q\in vert(P-P)\\ |t|<\|q\|}}\frac{\Theta_q(t)}{\|q\|}.$$

Thus the width law of every full-dimensional three-dimensional polytope is absolutely continuous and is reduced to finitely many spherical-circle intersection problems. The normal fan supplies a finite critical set, and between consecutive critical values the aperture functions, and hence the density, are real analytic. The formula therefore converts discrete incidence data from the normal fan of $P-P$ into analytic information about a probability density: the spherical stratification determines its branches, while changes in the active boundary incidences govern its possible losses of smoothness. Exact laws for particular solids are thereby placed within a general passage from finite polyhedral geometry to analytic pushforward measures, rather than obtained by unrelated coordinate integrations. Figure 1 displays the aperture mechanism on which this reduction rests.

The regular tetrahedron is the first non-zonotopal test of the theory. Its difference body is a cuboctahedron, whose triangular facets obstruct any representation of the width as a finite sum of absolute linear forms. Finch computed the mean and mean-square width of the regular tetrahedron and identified its density as an outstanding problem [8]. Theorem 4.3 determines that density exactly for the unit-edge tetrahedron, on support $[1/\sqrt{2}, 1]$ with critical values $1/\sqrt{2}$, $\sqrt{2/3}$, $\sqrt{3}/2$ and $1$.

The corresponding transitions are a linear onset, a corner, a square-root fold and a terminal step. These singularities arise from distinct incidences between a growing level circle and the boundary of a normal cell. The calculation therefore resolves the tetrahedral density problem while also exhibiting the local geometric mechanisms behind the different branches.

The same spherical-cell method applies to brightness, but with cells obtained from the projected-area function rather than from a difference body. The projected-area distribution of the regular tetrahedron was previously derived by Vickers and Brown [16], following earlier work on projected-area laws [14, 15, 17]. Section 5 gives an independent derivation from the common cell-aperture framework. The brightness chamber divides into two noncongruent families, shown in Figure 2, and this division explains an interior jump that has no analogue in the tetrahedral width law.

The calculation also reveals an exact pointwise relation with a second polyhedron. If $R$ is the unit-edge rhombic dodecahedron and $T$ the unit-edge regular tetrahedron, then Proposition 5.8 proves $w_R(u) = 8\, b_T(u)/\sqrt{3}$ for every $u \in S^2$.

Theorem 5.9 consequently transfers the complete tetrahedral brightness law to the width law of $R$. This is not merely an equality of distributions: it identifies the two directional functions up to scale.

The third exact model is the regular truncated octahedron. It lies on the zonotopal side of Proposition 2.6 and is simultaneously a three-dimensional permutohedron, connecting the calculation with graphical zonotopes and their combinatorics [22]. It is also distinguished among parallelohedra by Lángi’s fixed-volume mean-width extremal theorem [23]. Theorem 8.3 gives its complete unit-edge width density on $[\sqrt{6}, \sqrt{10}]$, with critical values $\sqrt{6}$, $2\sqrt{2}$, $3$ and $\sqrt{10}$.

Its four transitions again realize the onset, corner, fold and terminal-step mechanisms predicted by the normal-cell analysis. The three flagship densities are compared in Figure 3. Their different interior singularities make visible the distinction between the non-zonotopal tetrahedral width and the zonotopal laws arising from brightness and the truncated octahedron.

### 1.3. Functional and harmonic calculi

Section 6 develops two complementary calculi for directional functionals. The normal-cell decomposition applies to every three-dimensional polytope: Theorem 6.1 gives a Borel functional calculus for integrals of the form $\int_{S^2} \Phi(w_P(u))\, d\sigma(u)$, so distribution functions, moments and transforms follow from one cellwise identity, even when the density itself is not the simplest object to compute.

On the zonotopal side, moments of the cosine transform $A_\nu(u)=\int_{S^2}|u\cdot n|\,d\nu(n)$ polarize into kernels determined by the Gram geometry of finite direction systems. The second-order kernel is an explicit strictly increasing function of $|n_i\cdot n_j|$; the third-order kernel distinguishes the two switching classes of a signed triple and remains valid on its rank-two boundary. Thus low-order absolute-product data become geometric reconstruction data.

The harmonic form of the theory gives a complementary explanation. Theorem 6.12 decomposes $Var\, A_\nu(U)$ into even spherical-harmonic components, while Appendix C identifies the degree-two term with the failure of the weighted direction system to be a unit-norm tight frame. If its symmetry group acts irreducibly on $\mathbb{R}^3$, Schur's lemma makes the frame operator scalar and annihilates this contribution. The truncated-octahedral and tetrahedral generating systems therefore begin at degree four. This explains a feature of their exact variances that direct integration does not reveal and supplies a criterion beyond the examples treated here.

### 1.4. Finite normal tomography

Let $n_1,\ldots,n_m\in S^2$ be labelled directions and let $\Gamma$ specify the observed pairs. Pair data recover the absolute Gram matrix, but they are unchanged under $n_i\mapsto\varepsilon_i Q n_i$, where $\varepsilon_i\in\{\pm1\}$ and $Q\in O(3)$. Triple Gram-product signs are invariant under these vertex sign changes, and their compatibility is governed by the cycle space of $\Gamma$.

Theorem 7.3 proves that, when the observed triangles span that cycle space, labelled second- and third-order polarized moments determine the directions up to the unavoidable orthogonal and sign gauges. For the complete graph, all pairs and an anchored family of triangles give the explicit sufficient count $\binom{m}{2}+\binom{m-1}{2}=(m-1)^2$. The reconstruction is stable when the observed correlations stay uniformly separated from zero: second-order data determine magnitudes, third-order data determine cycle signs, and signed-graph switching completes the argument [21].

This decomposition into absolute edge weights and switching-invariant cycle signs is also the algebraic mechanism behind the unrestricted Gaussian reconstruction of Theorem 10.1. The rank constraint changes, but the gauge group and cycle-space obstruction do not. For the truncated octahedron, the correlation graph of the six generating directions is $K_{2,2,2}$. Theorem 8.6 identifies the triangle switching classes with the coplanar triples and hence the hexagonal facets, while the orthogonal pairs identify the square facets. The switching data, the rank-two boundary of the three-direction kernel and the facet incidence structure are three manifestations of one combinatorial object.

The finite tomography theorem concerns labelled polarized data; it does not assert that the scalar distribution of one width functional determines its generating configuration. Integration over the sphere contracts labelled tensors into ordinary moments and can erase information retained before contraction.

### 1.5. What the scalar law does not determine

The loss of information is already visible at the level of classical valuations. Mean width, mean brightness and volume determine the intrinsic volumes $V_1, V_2, V_3$ in dimension three, but they do not determine the corresponding probability laws.

Section 9 constructs two centrally symmetric right prisms $K_1$ and $K_2$ with $V_j(K_1)=V_j(K_2)$ for $j=1,2,3$, while both their width laws and brightness laws differ. The planar bases are smooth and strictly convex, are visibly distinct in Figure 4, and have equal area and perimeter. Their second directional moments separate the prisms explicitly: Theorem 9.2 gives $E\,W_{K_1}^2 - E\,W_{K_2}^2 = 16\,t^2/15 > 0$ and $E\,B_{K_1}^2 - E\,B_{K_2}^2 = 16\,H^2 t^2/15 > 0$.

Thus equality of all intrinsic volumes does not determine either law, even in a centrally symmetric family. This complements the positive finite tomography theorem: labelled low-order data can be reconstructive, whereas the classical scalar summaries of a directional law need not be.

Relative to the intrinsic-volume nonclosure established in [1, Theorem 7.5], Theorem 9.2 adds four features: the spatial witnesses are centrally symmetric; their planar bases are smooth and strictly convex; the same pair separates both the width and brightness laws; and both separations are certified by exact positive second-moment gaps. The resulting three-dimensional bodies are right prisms rather than smooth strictly convex bodies, and the corresponding smooth three-dimensional strengthening remains open.

### 1.6. A dimension-free Gaussian reconstruction theorem

The signed-graph reconstruction mechanism persists after the rank-three geometric constraint is removed. Let $X=(X_1,\dots,X_m)$ be a centered Gaussian vector with unit marginal variances and correlation matrix $\Sigma=(\rho_{ij})$. Here $\Sigma$ may have arbitrary rank, so the resulting theorem is not an application of the three-dimensional normal tomography of Section 7. It is a dimension-free inverse theorem for Gaussian absolute moments, governed by the same separation between edge magnitudes and cycle signs.

The pairwise absolute moments $E|X_i X_j|$ determine $|\rho_{ij}|$, but not the correlation signs. The remaining ambiguity is diagonal sign conjugacy, $\Sigma \mapsto D_\varepsilon \Sigma D_\varepsilon$, where $D_\varepsilon = diag(\varepsilon_1,\dots,\varepsilon_m)$.

Third absolute moments provide the cycle-sign information. Theorem 10.1 proves, without a rank restriction, that labelled pair moments and labelled triple moments determine $\Sigma$ up to this unavoidable gauge whenever the observed triangles span the cycle space of the correlation graph. For the complete graph, the same anchored construction gives $(m-1)^2$ sufficient observations. No fourth-order moment is required.

Equivalently, a standard centered multivariate folded-normal vector is identifiable from these labelled moments up to coordinate sign conjugacy. The classical literature provides formulae and computational methods for forward absolute moments of Gaussian vectors [11, 12, 20, 24] and studies the distribution, estimation and transforms of multivariate folded-normal laws [25, 26]. Theorem 10.1 addresses the inverse question: which finite collection of labelled absolute moments reconstructs the underlying correlation structure?

The triangle-span hypothesis has genuine content. Proposition 10.2 constructs positive-definite rank-four correlation matrices on a four-cycle that are not sign-gauge equivalent but have identical labelled absolute moments on every set of at most three coordinates. Hence no dimension-free theorem based only on such local moments can omit the cycle condition. This example does not settle the corresponding rank-three directional problem, and the count $(m-1)^2$ is asserted as an explicit sufficient count rather than as a global optimum among all possible observations.

Section 10 concludes by comparing these finite inverse statements with the full cosine transform. The classical Aleksandrov injectivity theorem, for which a self-contained harmonic proof is supplied in Proposition 10.3, gives $A_\nu = A_\eta$ if and only if $\nu_e = \eta_e$, where $\nu_e$ is the even part of $\nu$; see Gardner [9], Schneider [13] and Ournycheva–Rubin [19]. This places the results in an information hierarchy:

$$\text{labelled cosine transform} \Leftrightarrow \text{even generating measure} \Rightarrow$$

$$\text{labelled polarized moments} \Rightarrow \text{scalar law} \Leftrightarrow \text{all ordinary moments}.$$

The reverse implications generally fail. Theorems 7.3 and 10.1 identify structured situations in which finitely many labelled contractions recover information that scalar averaging would otherwise discard.

### *1.7. Relation to earlier work and organization*

The present work continues the study of width distributions begun with the exact rectangular-box laws in [1] and the reduced-zonotope and spectral framework in [2]. Those papers establish the tractable behaviour of particular zonotopal families and show, in a planar smooth setting, that scalar width laws need not determine the underlying body. The present paper is self-contained and does not use either earlier work as a proof premise. Its general results precede the named examples: Proposition 2.6 identifies the exact zonotopal boundary, Theorems 3.1 and 3.4 establish the global spherical co-area and density formulae in the normalization used here, and the subsequent sections develop the functional, harmonic and inverse theories from first principles.

Section 2 gives the measure dictionary for widths, brightness and cosine transforms, and proves the zonotopal characterization. Section 3 establishes the global spherical density theory. Sections 4 and 5 determine the tetrahedral width and brightness laws and derive the rhombic-dodecahedral law. Section 6 develops the functional, polarized-moment and harmonic calculi. Section 7 proves finite normal tomography. Section 8 determines the width law and facet reconstruction of the regular truncated octahedron. Section 9 proves that intrinsic volumes do not determine width or brightness laws. Section 10 gives the dimension-free Gaussian reconstruction theorem, its four-cycle limitation and the cosine-transform information hierarchy. Appendices A and B contain exact auxiliary calculations; Appendix C gives the tight-frame criterion and the harmonic variance decompositions.

## 2. Width, brightness, and spherical measures

Width and brightness encode different directional data. Width measures the separation of supporting hyperplanes, whereas brightness is the $(d-1)$-dimensional volume of an orthogonal projection. For zonoids and projection bodies, however, both are expressed through the classical spherical cosine transform. We fix the normalizations needed below and keep the generating measure of a zonoid distinct from the surface-area measure of a general convex body.

### 2.1. Support functions, width, and difference bodies

Let $K \subset \mathbb{R}^d$ be a compact convex body, with support function

$$h_K(u) = \max_{x \in K} u \cdot x, u \in S^{d-1}.$$

Its width in the deterministic direction $u$ is $w_K(u) = h_K(u) + h_K(-u)$, and its difference body is $DK = K + (-K)$.

**Proposition 2.1 (difference-body identity).** For every compact convex body $K \subset \mathbb{R}^d$,

$$w_K(u) = h_{DK}(u), u \in S^{d-1}.$$

In particular, the width function is translation invariant and depends on $K$ only through $DK$.

*Proof.* Support functions are Minkowski additive and satisfy $h_{-K}(u) = h_K(-u)$. Hence $h_{DK}(u) = h_K(u) + h_K(-u) = w_K(u)$. Translation leaves $DK$ unchanged. ▫

### 2.2. Centered zonotopes

Let $n_1, \dots, n_m \in S^{d-1}$ and $c_1, \dots, c_m > 0$, and put

$$Z = \sum_{i=1}^{m} \left[ -\frac{c_i}{2} n_i, \frac{c_i}{2} n_i \right].$$

The centered convention is essential: the support function of the $i$th segment is $c_i |u \cdot n_i| / 2$.

**Proposition 2.2 (centered-zonotope width).** For the centered zonotope $Z$,

$$h_Z(u) = \frac{1}{2} \sum_{i=1}^{m} c_i |u \cdot n_i|, w_Z(u) = \sum_{i=1}^{m} c_i |u \cdot n_i|.$$

*Proof.* The first identity follows by Minkowski additivity of support functions. Since $Z = -Z$, its support function is even and $w_Z = 2h_Z$. ▫

For general zonoids the same representation holds with a unique even generating measure; see Bolker [3] and Schneider [13, Section 3.5].

### 2.3. Cosine-transform notation

**Definition 2.3 (cosine-transform notation).** For a finite Borel measure $\nu$ on $S^{d-1}$, we write

$$A_\nu(u) = \int_{S^{d-1}} |u \cdot n| d\nu(n). (1)$$

This is the classical spherical cosine transform of the measure $\nu$; the notation $A_\nu$ is used only to keep the two geometric weightings below visibly distinct. If $\nu^-(E) = \nu(-E)$ and $\nu_e = (\nu + \nu^-)/2$, then the evenness of the kernel gives $A_\nu = A_{\nu_e}$.

For the zonotope in Proposition 2.2, define

$$\mu_Z = \frac{1}{2}\sum_{i=1}^{m} c_i\left(\delta_{n_i} + \delta_{-n_i}\right).$$

Then $w_Z = A_{\mu_Z}$. More generally, a zonoid has a unique even generating measure $\mu_Z$ and width $A_{\mu_Z}$ [3, 13].

### 2.4. Brightness and the measure dictionary

For a convex body $K \subset \mathbb{R}^d$, let

$$b_K(u) = vol_{d-1}\left(K \mid u^{\perp}\right)$$

be its brightness. If $S_K$ denotes the surface-area measure, Cauchy's projection formula gives

$$b_K(u) = \frac{1}{2}\int_{S^{d-1}} |u \cdot n|\, d\, S_K(n),$$

with the standard normalization used in Gardner [9, Section 4.1] and Schneider [13, Sections 5.3 and 10.9].

**Proposition 2.4 (classical width-brightness measure dictionary).** Let $Z$ be a zonoid with its unique even generating measure $\mu_Z$, and let $K$ be a convex body with surface-area measure $S_K$. Then

$$w_Z = A_{\mu_Z}, b_K = A_{S_K/2} = A_{(S_K)_e/2}. \quad (2)$$

*Proof.* The first identity is the classical generating-measure representation of a zonoid [3, 13]. Cauchy's projection formula gives the second, and (1) is unchanged by even symmetrization. ▫

Equation (2) is the measure dictionary used throughout the paper. The measures $\mu_Z$ and $S_K$ encode different geometry: the first records the generating directions and weights of a zonoid, while the second records boundary-area distribution. Their cosine transforms have the same analytic form, but the measures themselves must not be conflated.

### 2.5. Polyhedral brightness as zonotopal width

For a polyhedron $P \subset \mathbb{R}^3$, let $F$ range over its facets, with area $A_F$ and outward unit normal $n_F$. Then

$$S_P = \sum_F A_F \delta_{n_F}, b_P(u) = \frac{1}{2}\sum_F A_F |u \cdot n_F|.$$

The projection body $\Pi P$ is the centrally symmetric body with support function $h_{\Pi P} = b_P$ [9, Section 4.1; 13, Section 10.9].

**Proposition 2.5 (polyhedral brightness as zonotopal width).** Define

$$Z_P = \sum_F \left[-\frac{A_F}{4} n_F, \frac{A_F}{4} n_F\right].$$

Then $Z_P = \Pi P / 2$ and, pointwise on $S^2$,

$$w_{Z_P}(u) = b_P(u).$$

Consequently, for the same uniform random direction $U$,

$$W_{Z_P} = B_P$$

pointwise and hence also in distribution.

*Proof.* Proposition 2.2 gives $h_{Z_P}(u) = \frac{1}{4}\sum_F A_F |u \cdot n_F| = b_P(u)/2$. Equality of support functions yields $Z_P = \Pi P/2$. Since $Z_P$ is centrally symmetric, $w_{Z_P} = 2h_{Z_P} = b_P$. ▫

This pointwise identity is the bridge from polyhedral brightness to non-orthogonal zonotopal width used in Section 5.

### 2.6. The zonotopal width dichotomy

The preceding identities separate two genuinely different regimes for width. The distinction is intrinsic to the difference body and does not depend on a particular choice of generators.

**Proposition 2.6 (zonotopal width dichotomy).** Let $K \subset \mathbb{R}^3$ be a compact convex body. The following are equivalent:

1. there exist $c_1, \dots, c_m > 0$ and $n_1, \dots, n_m \in S^2$ such that
$$w_K(u) = \sum_{i=1}^{m} c_i |u \cdot n_i|, u \in S^2;$$
2. the difference body $DK = K - K$ is a zonotope.

If these conditions hold, the width law admits the generator-based sign-cell calculus of the zonotopal regime. If $DK$ has a two-dimensional face that is not centrally symmetric, then neither condition can hold.

*Proof.* Suppose first that

$$w_K(u) = \sum_{i=1}^{m} c_i |u \cdot n_i|.$$

Let

$$Z = \sum_{i=1}^{m} \left[ -\frac{c_i}{2} n_i, \frac{c_i}{2} n_i \right].$$

By Proposition 2.2, $w_Z = w_K$. Proposition 2.1 gives $h_{DK} = w_K$, while the central symmetry of $Z$ gives $w_Z = 2h_Z = h_{2Z}$. Hence $h_{DK} = h_{2Z}$, so $DK = 2Z$ and is a zonotope.

Conversely, if

$$DK = \sum_{i=1}^{m} \left[ -\frac{d_i}{2} n_i, \frac{d_i}{2} n_i \right],$$

then Proposition 2.1 and the support-function formula for a centered zonotope give

$$w_K(u) = h_{DK}(u) = \frac{1}{2}\sum_{i=1}^{m} d_i |u \cdot n_i|.$$

Finally, every face of a zonotope is itself a zonotope, and every two-dimensional zonotope is centrally symmetric. Thus a non-centrally-symmetric two-face of $DK$ obstructs zonotopality and therefore obstructs every finite absolute-sum representation of $w_K$. ▫

This equivalence appeared, with the normalized central symmetral $1/2(K - K)$, in [1, Proposition 7.7]. Proposition 2.6 records the form used here. It is the structural boundary between the two calculi developed below: the normal-fan cell calculus applies to arbitrary polytope widths, whereas the polarized absolute-product calculus of Section 6 requires the zonotopal side of the dichotomy.

# 3. Global spherical co-area and arbitrary polytope width laws

From this section onward we work in $\mathbb{R}^3$, with directions on $S^2$. Proposition 2.6 identifies the exact boundary between zonotopal and non-zonotopal width laws. In *Width distributions for rectangular boxes* [1, Proposition 7.7], we proved the equivalent statement using the normalized central symmetral. That paper uses the normalization $D_0 K = 1/2(K - K)$, whereas the present paper writes $DK = K - K$; this is only a positive dilation and therefore does not affect the zonotope property. The corresponding density formula for the sign cells of the generator arrangement was established in [1, Proposition 7.9], and the great-circle formulation was developed further in *Width Laws and Spectral Geometry: Universality, Geometric Memory, and Reconstruction* [2, Section 8.1].

We isolate the part of the earlier co-area argument that does not depend on a generator representation, replacing zonotopal sign cells by arbitrary finite spherical polyhedral cells with cellwise linear data. The cancellation between level-circle radius and spherical-gradient magnitude is the same three-dimensional mechanism used for rectangular boxes in [1, Theorem 3.1]; only the indexing geometry changes.

## *3.1. Cellwise spherical co-area*

Let $C_2$ be a finite family of two-dimensional spherical polyhedral cells in $S^2$. We assume that their interiors are pairwise disjoint and that

$$\sigma\left(S^2 \setminus \bigcup_{C \in C_2} \int C\right) = 0.$$

Their boundaries lie in finitely many great circles and hence have spherical area zero.

Let $F : S^2 \to \mathbb{R}$ be measurable. Suppose that for every $C \in C_2$ there is a nonzero vector $a_C \in \mathbb{R}^3$ such that

$$F(u) = a_C \cdot u$$

for $u \in \int C$. Values assigned on cell boundaries do not affect the law. Write $D_C = \| a_C \|$, and for $|t| < D_C$ let $\Theta_C(t)$ be the total angular aperture, about the axis $a_C / D_C$, of the part of the level circle $a_C \cdot u = t$ lying in $\int C$.

**Theorem 3.1 (global spherical co-area).** Let $U$ be uniformly distributed on $S^2$. Under the preceding hypotheses, the law of $F(U)$ is absolutely continuous. For almost every $t$, a density representative is

$$f_F(t) = \frac{1}{4\pi} \sum_{\substack{C \in C_2 \\ |t| < D_C}} \frac{H^1\left(\int C \cap \{u \in S^2 : a_C \cdot u = t\}\right)}{\sqrt{D_C^2 - t^2}}, (3)$$

or, equivalently,

$$f_F(t) = \frac{1}{4\pi} \sum_{\substack{C \in C_2 \\ |t| < D_C}} \frac{\Theta_C(t)}{D_C}. (4)$$

Density representatives are defined only up to Lebesgue-null sets. In particular, a jump of a chosen representative at an endpoint or critical value does not by itself represent an atom.

*Proof.* Fix $C \in C_2$. On its interior,

$$\nabla_{S^2}(a_C \cdot u) = a_C - (a_C \cdot u)u.$$

Along the level $a_C \cdot u = t$, its norm is $\sqrt{D_C^2 - t^2}$. For a Borel set $E \subset (-D_C, D_C)$, the co-area formula [6] applied to $\int C$ gives

$$\sigma\left(\int C \cap F^{-1}(E)\right) = \frac{1}{4\pi} \int_E \frac{H^1\left(\int C \cap \{a_C \cdot u = t\}\right)}{\sqrt{D_C^2 - t^2}} dt.$$

The endpoint levels $t = \pm D_C$ contain at most the axial points $\pm a_C / D_C$ and therefore have zero spherical area. Since the two-dimensional cell interiors are disjoint, their complement has spherical measure zero by hypothesis, and there are only finitely many cells, summing proves absolute continuity and (3). If closed cells were used instead, a boundary arc could be counted twice only at exceptional levels; this changes a density representative on at most a finite set and not the pushforward measure.

For (4), put $\hat{a}_C = a_C / D_C$ and choose an orthonormal basis $e_1, e_2$ of $\hat{a}_C^{\perp}$. The level circle is parametrized by

$$u(\theta) = \frac{t}{D_C} \hat{a}_C + \sqrt{1 - \frac{t^2}{D_C^2}} (\cos\theta\, e_1 + \sin\theta\, e_2).$$

Its arc element is $ds = \sqrt{1 - t^2 / D_C^2}\, d\theta$, while

$$\sqrt{D_C^2 - t^2} = D_C \sqrt{1 - \frac{t^2}{D_C^2}}.$$

The level-circle radius therefore cancels the spherical-gradient factor, leaving $\Theta_C(t) / D_C$. Substitution into (3) gives (4). ▫

The theorem is precisely the form needed below: [1, Proposition 7.9] is recovered when the cells are the sign cells of a zonotope generator arrangement, whereas no generator representation is required here. Its proof uses only the spherical co-area formula; references [1] and [2] locate earlier special cases and do not supply mathematical premises for Theorem 3.1.

### *3.2. Exact aperture calculus*

A spherical polyhedral cell is cut out by finitely many hemispherical inequalities. On a fixed level circle, each inequality becomes a circular interval condition.

**Proposition 3.2 (interval representation of a cell aperture).** Let

$$\overline{C}=\{u\in S^2 : b_j\cdot u\ge 0,\, j=1,\dots,N\},\, C=\int_{S^2}\overline{C},$$

and suppose $F(u)=a\cdot u$ on $C$, with $a\neq 0$. Put $D=\|a\|$, $\hat{a}=a/D$, and choose an orthonormal basis $e_1, e_2$ of $\hat{a}^{\perp}$. On $a\cdot u=t$, $|t|<D$, write

$$u(\theta)=\frac{t}{D}\hat{a}+\sqrt{1-\frac{t^2}{D^2}}\left(\cos\theta\, e_1+\sin\theta\, e_2\right).$$

For each $j$, let

$$\beta_j=\sqrt{(b_j\cdot e_1)^2+(b_j\cdot e_2)^2}.$$

If $\beta_j>0$, choose $\psi_j$ by

$$b_j\cdot e_1=\beta_j\cos\psi_j,\, b_j\cdot e_2=\beta_j\sin\psi_j.$$

Then the $j$th cell inequality is equivalent to

$$\cos(\theta-\psi_j)\ge c_j(t),\, c_j(t)=-\frac{t\, b_j\cdot\hat{a}}{D\beta_j\sqrt{1-t^2/D^2}}.$$

Thus the aperture is the angular measure of a finite intersection of explicit circular intervals. If $\beta_j=0$, the corresponding inequality is constant on the level circle and either retains the whole circle or removes it.

*Proof.* Substituting the level parametrization into $b_j\cdot u\ge 0$ gives

$$\frac{t}{D}b_j\cdot\hat{a}+\sqrt{1-\frac{t^2}{D^2}}\left[(b_j\cdot e_1)\cos\theta+(b_j\cdot e_2)\sin\theta\right]\ge 0.$$

For $\beta_j>0$, the bracket is $\beta_j\cos(\theta-\psi_j)$. If $-1<c_j(t)<1$, the interval endpoints are

$$\theta=\psi_j\pm\arccos c_j(t)\,(mod\, 2\pi);$$

if $c_j(t)\le -1$ the entire circle is admitted, while if $c_j(t)>1$ none of it is. The case $\beta_j=0$ is independent of $\theta$. Intersecting the finitely many admissible angular sets gives the aperture. ▫

The interval mechanism of Proposition 3.2 is illustrated in Figure 1. Panel (a) shows the level-circle geometry of the tetrahedral normal cell used in Section 4. Panels (b) and (c) show the two nonconstant forbidden-arc regimes that arise in the truncated-octahedral normal cell of Section 8. In both cases the aperture is obtained as the angular measure of the complement of finitely many explicitly determined circular intervals.

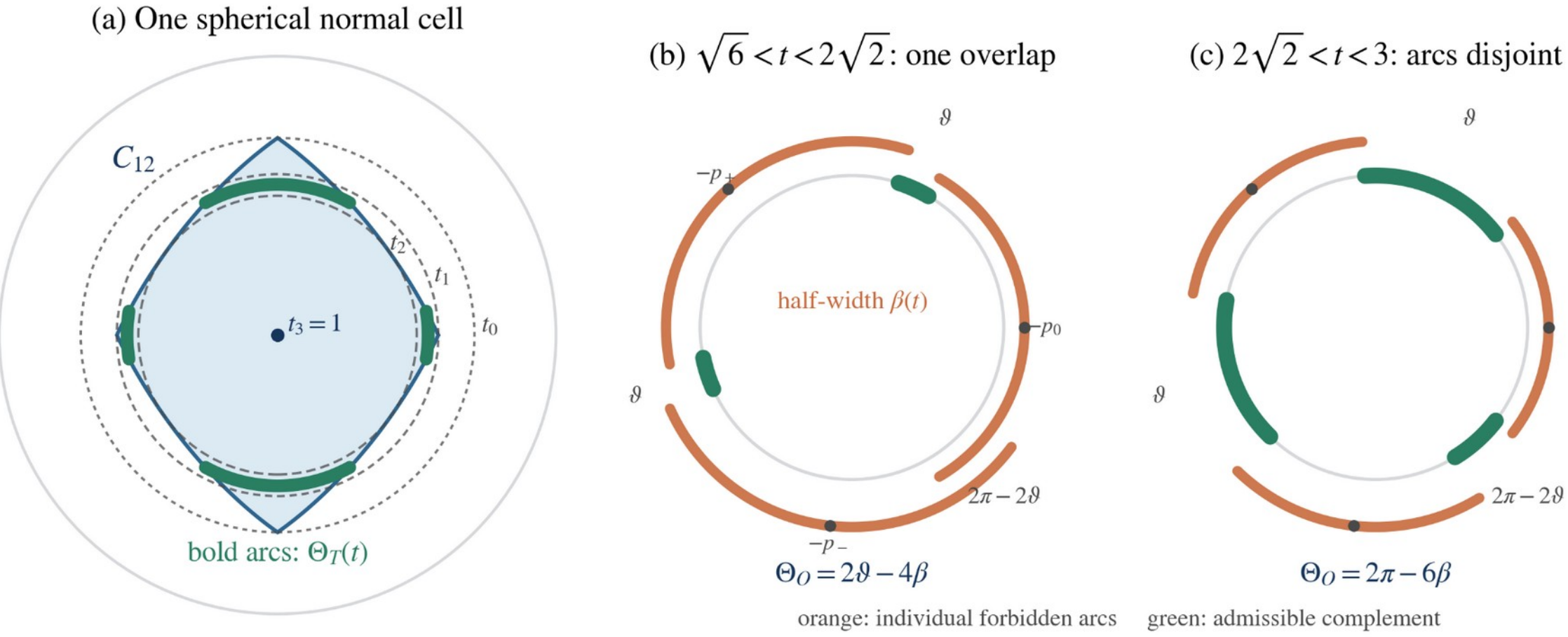

**Figure 1. Exact aperture geometry for tetrahedral and truncated-octahedral normal cells.** Panel (a) is the exact orthographic projection along $e_0$ of the tetrahedral cell $C_{12}$. The concentric curves are the projections of the levels $e_0 \cdot u = t_j$, $j = 0, 1, 2$, while $t_3 = 1$ projects to the axial point. The bold portions of the intermediate level circle are the admissible arcs whose total angular length is $\Theta_T(t)$. Panels (b) and (c) show representative levels of the truncated-octahedral cell in the regimes $\sqrt{6} < t < 2\sqrt{2}$ and $2\sqrt{2} < t < 3$. The three forbidden intervals are centered at $-p_0, -p_+, -p_-$, have common half-width $\beta(t)$, and have cyclic center gaps $\vartheta, \vartheta, 2\pi - 2\vartheta$. Their admissible complements have respective apertures $2\vartheta - 4\beta(t)$ and $2\pi - 6\beta(t)$. The orange intervals are displaced slightly in the radial direction for legibility; angular positions, incidences, and apertures are determined exactly by (6), (7), and (48).

### 3.3. Critical values and analytic branches

For zonotopal sign-cell arrangements, exact incidence tests and the associated branch changes were developed in [1, Corollary 7.10] and [2, Section 8]. The following finite candidate list is the form needed for arbitrary spherical polyhedral cells.

**Proposition 3.3 (finite critical-value structure).** Let $C$ be a spherical polyhedral cell and suppose $F(u) = a \cdot u$ on $C$, with $a \neq 0$. A change in the combinatorial incidence of the level circle $a \cdot u = t$ with $C$ can occur only at one of the following values:

1. the axial values $t = \pm \| a \|$;
2. a boundary-tangency value

$$t = \pm \left\| a - \left(a \cdot \hat{b}\right) \hat{b} \right\|,$$

where $b \cdot u = 0$ supports a boundary arc and $\hat{b} = b / \| b \|$;

3. a vertex value $t = a \cdot v$, where $v$ is a spherical vertex of $\overline{C}$.

After candidates that do not meet the relevant boundary stratum are discarded, the remaining finite set partitions the range into intervals on which the aperture is represented by a fixed real-analytic expression.

*Proof.* The level circle collapses exactly when $|t| = \| a \|$. On a supporting great circle $b \cdot u = 0$, the extremal values of $a \cdot u$ are the extrema of the projection of $a$ onto $b^{\perp}$, which are the stated tangency candidates.

Away from axial collapse and tangency, each active interval endpoint from Proposition 3.2 varies real-analytically. Its cyclic order can change only if two active endpoints coincide. Coincident supporting great circles may first be represented by a single nonredundant constraint. Opposite orientations cannot occur for a two-dimensional cell: the simultaneous inequalities $b_j \cdot u \geq 0$ and $(-b_j) \cdot u \geq 0$ force $b_j \cdot u = 0$, so the cell has empty spherical interior. Thus, for distinct supporting circles, the corresponding normals are linearly independent and the circles meet at the antipodal points parallel to $b_j \times b_k$. If such an intersection is excluded by another cell constraint, neither coincident angle is an active endpoint of the actual intersection near that level, so the aperture formula does not change. Hence an active crossing can occur only at an intersection lying in $\overline{C}$, that is, at a spherical vertex, and its value is $a \cdot v$.

The list is therefore necessary. Whether a tangency or vertex candidate is active is an incidence question on the relevant boundary stratum; in the applications below this is checked directly from the defining inequalities. Once the inactive candidates are removed, no endpoint is created, destroyed, tangent, or reordered between consecutive critical values, so the aperture is a fixed finite sum of differences of analytic endpoint functions there. ▫

The theorem and propositions above provide a finite exact procedure. They do not assert that every candidate is active, nor do they require the cells to arise from a zonotope generator arrangement.

### 3.4. Width density of an arbitrary three-dimensional polytope

The preceding cellwise theorem becomes canonical for polytope widths when the cells are identified with the spherical normal fan of the difference body. Let $P \subset \mathbb{R}^3$ be a full-dimensional polytope and put

$$Q = DP = P - P.$$

For a vertex $q \in \mathit{vert}\, Q$, let

$$N_Q(q) = \{ y \in \mathbb{R}^3 : y \cdot q \geq y \cdot x \text{ for every } x \in Q \}$$

be its normal cone and define the spherical normal cell

$$C_q = N_Q(q) \cap S^2.$$

On the spherical interior of $C_q$, the support function of $Q$ is the single linear form $u \mapsto q \cdot u$.

**Theorem 3.4 (width density of an arbitrary three-dimensional polytope).** Let $P \subset \mathbb{R}^3$ be a full-dimensional polytope, let $Q = P - P$, and let $U$ be uniform on $S^2$. For each $q \in \mathit{vert}\, Q$, put $D_q = \| q \| > 0$ and let $\Theta_q(t)$ be the angular aperture of

$$\int_{S^2} C_q \cap \{ u \in S^2 : q \cdot u = t \}$$

about the axis $q / D_q$. Then the random width $W_P = w_P(U)$ is absolutely continuous and, for almost every $t$,

$$f_{W_P}(t) = \frac{1}{4\pi} \sum_{\substack{q \in \operatorname{vert} Q \\ |t| < D_q}} \frac{\Theta_q(t)}{D_q}.$$

For each normal cell, all branch changes belong to the finite critical-value list of Proposition 3.3. Consequently every three-dimensional polytope width law reduces to finitely many spherical normal cells, finitely many critical values, and finitely many real-analytic aperture branches.

*Proof.* Proposition 2.1 gives $w_P = h_Q$. Since $P$ is full-dimensional, 0 lies in the interior of $Q = P - P$; in particular, no vertex $q$ of $Q$ is the origin, so $D_q > 0$.

For a vertex $q$, one has

$$\int_{S^2} C_q = \int N_Q(q) \cap S^2,$$

and $u$ belongs to this set exactly when $q$ is the unique maximizer of $x \mapsto u \cdot x$ over $Q$. Distinct vertices therefore give pairwise-disjoint spherical interiors. The complement of their union is the set of directions for which the maximizing face has positive dimension; it is contained in the finite union of great circles

$$\{ u \in S^2 : (q - q') \cdot u = 0 \}, q \neq q',$$

and hence has $\sigma$-measure zero. Thus the spherical normal cells satisfy the hypotheses of Theorem 3.1. On $\int_{S^2} C_q$ ,

$$w_P(u) = h_Q(u) = q \cdot u,$$

so Theorem 3.1 gives the displayed density formula and absolute continuity. Proposition 3.3 supplies the finite candidate set for every cell.

For completeness, each vertex $q$ of $Q = P + (-P)$ is a difference of vertices of $P$. Choose $u \in \int N_Q(q)$. The exposed face $F_Q(u) = \{q\}$ is the Minkowski sum $F_P(u) + F_{-P}(u)$. A Minkowski sum is a singleton only when both summands are singletons, so $F_P(u) = \{v_i\}$ and $F_{-P}(u) = \{-v_j\}$ for vertices $v_i, v_j$ of $P$. Thus $q = v_i - v_j$. Hence the normal-fan formulation is precisely the canonical reduction of the max-min width decomposition, with duplicate and empty vertex-pair cells removed. ▫

Theorem 3.4 is the general polytope statement behind the tetrahedral construction in Section 4. The geometry of $P - P$ controls both the linear pieces and their incidence: its vertices index the two-dimensional spherical cells, its normal-fan skeleton supplies the switching set, and Proposition 3.3 reads the possible changes of density branch from that finite fan geometry.

### *3.5. A question on generic singularity types*

The exact laws in Sections 4 and 5 exhibit the same three elementary local mechanisms on opposite sides of Proposition 2.6: finite derivative corners, square-root folds, and finite jumps. This suggests a question at the level of the general cellwise theory rather than a conjecture tied to either example.

**Question 3.5 (generic singularity taxonomy).** Let $F : S^2 \to \mathbb{R}$ be continuous and linear on the two-cells of a finite spherical polyhedral decomposition. Suppose that an active critical value $t_*$ is generated by exactly one local event and that all other cell contributions are real-analytic through $t_*$. Under the corresponding nondegeneracy assumptions, do the following local events exhaust the elementary singularity types of the density?

1. a transverse axial endpoint, producing a finite jump;
2. a nondegenerate quadratic tangency of a level circle to one active boundary arc, producing a square-root fold;
3. a transverse crossing of two active boundary endpoints, producing a finite derivative corner.

More generally, must every other local singularity in dimension three arise from simultaneous active events, higher-order tangency or crossing, or superposition and cancellation of the elementary contributions above?

The question concerns active incidences only. It does not assert that every candidate in Proposition 3.3 is active, and it separates the local classification problem from the global incidence problem of deciding which candidates occur for a given polytope or cell decomposition.

## 4. The width law of the regular tetrahedron

For a unit-edge regular tetrahedron, Finch [8] determined the mean and mean-square width and stated that the probability density was not known. Kabluchko, Litvak and Zaporozhets [10] later obtained higher projection-length moments for regular polytopes. We now determine the full one-dimensional law.

Theorem 3.4 applies to every full-dimensional three-dimensional polytope; here it is carried to a complete closed form for the regular tetrahedron. Proposition 4.1 identifies its difference body as a non-zonotopal cuboctahedron, and the twelve max-min cells below are its spherical normal cells.

### *4.1. Cuboctahedral normal form*

Let

$$v_1=\frac{(1,1,1)}{\sqrt{8}},\, v_2=\frac{(1,-1,-1)}{\sqrt{8}},$$

$$v_3=\frac{(-1,1,-1)}{\sqrt{8}},\, v_4=\frac{(-1,-1,1)}{\sqrt{8}},$$

and set

$$T=conv\{v_1,v_2,v_3,v_4\}.$$

Every nonzero difference $v_i-v_j$ has two coordinates of absolute value $1/\sqrt{2}$ and one zero coordinate, so

$$\| v_i-v_j \| =1 (i\neq j).$$

Thus $T$ is a regular tetrahedron of edge length one.

**Proposition 4.1 (cuboctahedral normal form).** For $u=(x,y,z)\in S^2$,

$$w_T(u)=\frac{1}{\sqrt{2}}max\{|x|+|y|,|x|+|z|,|y|+|z|\}. (5)$$

Its range is

$$range(w_T)=\left[\frac{1}{\sqrt{2}},1\right].$$

The minimum is attained precisely at the six coordinate-axis directions. The maximum is attained precisely when one coordinate vanishes and the other two have equal absolute value.

*Proof.* Proposition 2.1 gives

$$w_T=h_{T-T}.$$

For polytopes,

$$T-T=conv\{v_i-v_j:1\leq i,j\leq 4\};$$

indeed, if $p=\sum_i \alpha_i v_i$ and $q=\sum_j \beta_j v_j$, then $p-q=\sum_{i,j}\alpha_i\beta_j(v_i-v_j)$, while the reverse inclusion is immediate. The twelve nonzero differences are

$$\frac{1}{\sqrt{2}}\{(\pm 1,\pm 1,0),(\pm 1,0,\pm 1),(0,\pm 1,\pm 1)\}.$$

Their convex hull is the regular cuboctahedron associated with the tetrahedral difference body [4, discussion following Corollary 3.4]. The three vertices

$$\frac{(1,1,0)}{\sqrt{2}},\frac{(1,0,1)}{\sqrt{2}},\frac{(0,1,1)}{\sqrt{2}}$$

form a triangular face in the supporting plane $x+y+z=\sqrt{2}$. Hence $T-T$ is not a zonotope, because every two-dimensional face of a zonotope is centrally symmetric. By Proposition 2.6, $w_T$ therefore admits no finite representation as a sum of absolute linear forms. Maximizing the scalar product with $u$ over the twelve difference vertices gives (5).

For the range, let

$$a\geq b\geq c\geq 0$$

be the decreasing rearrangement of $|x|,|y|,|z|$. Then

$$w_T(u)=\frac{a+b}{\sqrt{2}},a^2+b^2+c^2=1.$$

Since $a,b\geq c$,

$$(a+b)^2=1-c^2+2ab\geq 1+c^2\geq 1.$$

Thus $w_T(u)\geq 1/\sqrt{2}$. Equality requires $c=0$ and $ab=0$, hence $(a,b,c)=(1,0,0)$.

Conversely,

$$a+b \le \sqrt{2(a^2+b^2)} \le \sqrt{2},$$

so $w_T(u) \le 1$. Equality holds precisely when $a=b$ and $c=0$, which gives $a=b=1/\sqrt{2}$. ▫

Because $w_T$ is continuous on $S^2$ and normalized spherical measure has full support, the support of the law of $W_T$ is exactly this range.

### *4.2. The twelve max-min cells*

For $i \neq j$, define

$$C_{ij} = \{u \in S^2 : u \cdot v_i \ge u \cdot v_k \text{ for every } k, u \cdot v_j \le u \cdot v_k \text{ for every } k\}.$$

Let $N$ be the union of the great circles on which two tetrahedral vertex values coincide. If $u \notin N$, the maximizing and minimizing vertices are unique. Hence the interiors of the twelve cells $C_{ij}$ are pairwise disjoint and their union is $S^2 \setminus N$. Since $\sigma(N)=0$, this is a cell decomposition for the purposes of Theorem 3.1.

**Lemma 4.2 (spherical-cell reduction).** The twelve cells $C_{ij}$ are congruent, and on $C_{ij}$,

$$w_T(u) = (v_i - v_j) \cdot u.$$

For $C_{12}$, set

$$e_0 = v_1 - v_2 = \frac{(0,1,1)}{\sqrt{2}}, e_1 = (1,0,0), e_2 = \frac{(0,1,-1)}{\sqrt{2}}.$$

On the level $e_0 \cdot u = t$, write

$$u(t,\theta) = t\, e_0 + \sqrt{1-t^2}(\cos\theta\, e_1 + \sin\theta\, e_2).$$

Then $u(t,\theta) \in C_{12}$ precisely when

$$\sqrt{2}|\cos\theta| + |\sin\theta| \le \frac{t}{\sqrt{1-t^2}}. \quad (6)$$

Put

$$\alpha = \arctan\frac{1}{\sqrt{2}}, \gamma(t) = \arccos\left(\frac{t}{\sqrt{3(1-t^2)}}\right),$$

and define

$$t_0 = \frac{1}{\sqrt{2}}, t_1 = \sqrt{\frac{2}{3}}, t_2 = \frac{\sqrt{3}}{2}, t_3 = 1.$$

The angular aperture of the level circle inside $C_{12}$ is

$$\Theta_T(t)=\begin{cases} 0, & t<t_0, \\ 2\pi-4\alpha-4\gamma(t), & t_0\le t<t_1, \\ 2\pi-8\gamma(t), & t_1\le t<t_2, \\ 2\pi, & t_2\le t<t_3, \\ 0, & t\ge t_3. \end{cases} \quad (7)$$

*Proof.* On $C_{ij}$, $v_i$ is a maximizing vertex and $v_j$ a minimizing vertex; hence

$$w_T(u)=u\cdot v_i-u\cdot v_j=(v_i-v_j)\cdot u.$$

The tetrahedral symmetry group acts transitively on ordered pairs of distinct vertices. The same symmetry carries $v_i-v_j$ to the coefficient vector of the image cell, so it transports both the cell and its level-circle aperture. Thus all twelve cells have the same aperture function.

It remains to determine one cell. Put

$$\kappa(t)=\sqrt{1-t^2}.$$

The chosen parametrization gives

$$x=\kappa(t)\cos\theta,\ y=\frac{t+\kappa(t)\sin\theta}{\sqrt{2}},\ z=\frac{t-\kappa(t)\sin\theta}{\sqrt{2}}.$$

The inequalities defining $C_{12}$ reduce to

$$(v_1-v_3)\cdot u\ge 0,(v_1-v_4)\cdot u\ge 0,$$

$$(v_3-v_2)\cdot u\ge 0,(v_4-v_2)\cdot u\ge 0.$$

The remaining comparison

$$(v_1-v_2)\cdot u\ge 0$$

follows from the displayed comparisons: for example, adding $(v_1-v_3)\cdot u\ge 0$ and $(v_3-v_2)\cdot u\ge 0$ gives $(v_1-v_2)\cdot u\ge 0$.

Substitution of the coordinates reduces the four displayed inequalities to

$$t+\kappa(t)\left(\varepsilon_1\sqrt{2}\cos\theta+\varepsilon_2\sin\theta\right)\ge 0,\varepsilon_1,\varepsilon_2\in\{-1,1\}.$$

They hold simultaneously if and only if

$$t\ge\kappa(t)\max_{\varepsilon_1,\varepsilon_2=\pm 1}\left(\varepsilon_1\sqrt{2}\cos\theta+\varepsilon_2\sin\theta\right),$$

and the maximum equals

$$\sqrt{2}|\cos\theta|+|\sin\theta|.$$

This proves (6).

The inequality is invariant under the four reflections of the $\theta$-circle, so it suffices to work on

$$0 \le \theta \le \frac{\pi}{2}.$$

There

$$\sqrt{2}\cos\theta + \sin\theta = \sqrt{3}\cos(\theta - \alpha),$$

because

$$\sin\alpha = \frac{1}{\sqrt{3}}, \cos\alpha = \sqrt{\frac{2}{3}}.$$

Writing

$$\chi(t) = \frac{t}{\sqrt{1-t^2}},$$

condition (6) becomes

$$\cos(\theta - \alpha) \le \frac{\chi(t)}{\sqrt{3}}.$$

The function $\chi$ is strictly increasing, and the values $\chi = 1, \sqrt{2}, \sqrt{3}$ occur at $t_0, t_1, t_2$, respectively.

For $\chi < 1$ there is no admissible angle, since the minimum of $\sqrt{2}\cos\theta + \sin\theta$ on the first quadrant is 1. For $1 \le \chi \le \sqrt{2}$, one has

$$\alpha \le \gamma(t) \le \frac{\pi}{2} - \alpha,$$

and the admissible interval has length

$$\frac{\pi}{2} - \alpha - \gamma(t).$$

For $\sqrt{2} \le \chi \le \sqrt{3}$, one has $0 \le \gamma(t) \le \alpha$, and the two admissible intervals have total length

$$\frac{\pi}{2} - 2\gamma(t).$$

For $\chi \ge \sqrt{3}$, the entire first quadrant is admissible. Multiplying by four gives (7). ▫

The cell, its four critical levels, and a representative admissible aperture are shown in Figure 1(a).

### *4.3. Exact density*

Each cellwise coefficient is a vector $v_i - v_j$ of norm one. By Theorem 3.1 and the twelvefold congruence,

$$f_T(t) = \frac{12}{4\pi}\Theta_T(t) = \frac{3}{\pi}\Theta_T(t).$$

**Theorem 4.3 (exact tetrahedral width density).** Let $U$ be uniform on $S^2$, and set

$$W_T = w_T(U).$$

Then $W_T$ is absolutely continuous, with support

$$\left[\frac{1}{\sqrt{2}},1\right],$$

and density

$$f_T(t)=\begin{cases} 0, & t<t_0, \\ 6-\dfrac{12}{\pi}\left(\alpha+\gamma(t)\right), & t_0\le t<t_1, \\ 6-\dfrac{24}{\pi}\gamma(t), & t_1\le t<t_2, \\ 6, & t_2\le t<t_3, \\ 0, & t\ge t_3. \end{cases} \quad (8)$$

where

$$\alpha=\arctan\frac{1}{\sqrt{2}},\gamma(t)=\arccos\left(\frac{t}{\sqrt{3\left(1-t^2\right)}}\right).$$

*Proof.* Proposition 4.1 gives the support. Lemma 4.2 gives the aperture in one of the twelve pairwise interior-disjoint cells. Since every cellwise coefficient has norm one, formula (4) gives

$$f_T(t)=\frac{3}{\pi}\Theta_T(t).$$

Substitution of (7) yields (8). Absolute continuity is part of Theorem 3.1. ▫

Finch [8] determined the first two moments but did not determine the density. Theorem 4.3 supplies the missing probability law in the same unit-edge normalization. Within the targeted searches under tetrahedral width, breadth, caliper-diameter, Feret-diameter, directional-diameter, and random-width terminology, no later closed analytic density was located. This statement is confined to that documented search scope and does not assert unrestricted historical priority.

### *4.4. Transition structure*

The branch points in (8) arise simultaneously in all twelve congruent cells. Their global contributions therefore add with the common factor 12; no cancellation occurs in the tetrahedral case.

**Lemma 4.4 (critical and endpoint behavior).** The density $f_T$ is real-analytic on

$$\left(t_0,t_1\right),\left(t_1,t_2\right),\left(t_2,t_3\right).$$

At the lower endpoint,

$$f_T(t)=\frac{24}{\pi}\left(t-t_0\right)+O\left(\left(t-t_0\right)^2\right)\left(t\downarrow t_0\right).$$

At $t_1$, the density is continuous and

$$f_T{}'\left(t_1^-\right)=\frac{36\sqrt{3}}{\pi},f_T{}'\left(t_1^+\right)=\frac{72\sqrt{3}}{\pi}.$$

At $t_2$,

$$f_T(t)=6-\frac{96}{\pi 3^{1/4}}\sqrt{t_2-t}+O\left(\left(t_2-t\right)^{3/2}\right)\left(t\uparrow t_2\right).$$

Finally,

$$\lim_{t\uparrow t_3} f_T(t)=6,$$

while the chosen density representative is zero for $t\geq t_3$. This terminal discontinuity carries no atom.

*Proof.* Direct evaluation gives

$$\gamma\left(t_0\right)=\frac{\pi}{2}-\alpha,\gamma\left(t_1\right)=\alpha,\gamma\left(t_2\right)=0,$$

so the adjacent branches agree at the interior transition values.

Differentiation gives

$$\gamma'(t)=-\frac{1}{\left(1-t^2\right)\sqrt{3-4t^2}}.\ (9)$$

At $t_0$, this equals $-2$, giving the linear onset. At $t_1$,

$$\gamma'\left(t_1\right)=-3\sqrt{3},$$

and the coefficients $-12/\pi$ and $-24/\pi$ on the two adjacent branches give the stated one-sided derivatives.

For $t_2=\sqrt{3}/2$, write

$$\rho_T(t)=\frac{t}{\sqrt{3\left(1-t^2\right)}}.$$

Then

$$\rho_T\left(t_2\right)=1,\rho_T{}'\left(t_2\right)=\frac{8}{\sqrt{3}},$$

and hence

$$1-\rho_T(t)=\frac{8}{\sqrt{3}}\left(t_2-t\right)+O\left(\left(t_2-t\right)^2\right).$$

Using

$$\arccos(1-\varepsilon)=\sqrt{2\varepsilon}+O\left(\varepsilon^{3/2}\right),$$

one obtains

$$\gamma(t)=\frac{4}{3^{1/4}}\sqrt{t_2-t}+O\left((t_2-t)^{3/2}\right).$$

Substitution into (8) gives the expansion at $t_2$. The upper-endpoint assertion follows from the constant terminal branch, while the absence of an atom follows from the absolute continuity proved in Theorem 3.1. ▫

Thus the lower endpoint has linear onset, $t_1$ is a finite derivative corner, $t_2$ has square-root fold behavior, and $t_3$ is a terminal step of the density representative. These are the same local mechanism types that appear in the box and zonotopal classifications of [1, Theorem 3.9] and [2, Section 8], but here they arise from non-zonotopal max-min cells. Because all twelve tetrahedral cells undergo the same incidence change with the same sign, their singular contributions reinforce rather than cancel.

The resulting density and its four transition types are displayed in the left panel of Figure 3.

### 4.5. Distribution function and exact normalization

The cell decomposition already implies normalization geometrically. A direct algebraic verification is nevertheless useful because it checks the branch constants independently of the area partition.

Define

$$G(t)=t\gamma(t)-\arctan\sqrt{3-4t^2},\quad t_0\le t\le t_2.$$

Equation (9) gives $G'(t)=\gamma(t)$, and

$$G(t_0)=t_0\left(\frac{\pi}{2}-\alpha\right)-\frac{\pi}{4},\quad G(t_1)=t_1\alpha-\frac{\pi}{6},\quad G(t_2)=0.$$

**Proposition 4.5 (distribution function).** The density in Theorem 4.3 is normalized,

$$\int_{\mathbb{R}} f_T(t)\,dt=1,$$

and its cumulative distribution function is

$$F_T(t)=\begin{cases} 0, & t<t_0,\\ 6(t-t_0)-\dfrac{12}{\pi}\left[\alpha(t-t_0)+G(t)-G(t_0)\right], & t_0\le t<t_1,\\ 6(t-t_0)-\dfrac{12\alpha}{\pi}(t_1-t_0)+\dfrac{12}{\pi}\left[G(t_0)+G(t_1)-2G(t)\right], & t_1\le t<t_2,\\ 6t-5, & t_2\le t<t_3,\\ 1, & t\ge t_3.\end{cases} \tag{10}$$

*Proof.* Integrating the three nonzero branches of (8), using $G'=\gamma$, gives

$$\int_{t_0}^{t_3} f_T(t)\,dt=6(1-t_0)-\frac{12\alpha}{\pi}(t_1-t_0)+\frac{12}{\pi}\left(G(t_0)+G(t_1)\right).$$

The displayed endpoint values reduce the right-hand side exactly to 1. The elementary arithmetic is recorded in Appendix A.3.

Integrating from the lower endpoint to a variable upper limit on each branch gives (10). The adjacent expressions have equal values at $t_1$ and $t_2$, and the final branch satisfies $F_T(t_3)=1$. Differentiation on the three open branch intervals recovers (8). ▫

The terminal branch is particularly simple:

$$F_T(t)=6t-5, t_2\le t\le 1.$$

For example,

$$P(W_T\ge t_2)=6-3\sqrt{3}.$$

The affine terminal branch makes such upper-tail probabilities elementary.

### 4.6. All real power moments

The same cell coordinates yield the complete real-power moment family without integrating the piecewise density. We record the resulting one-dimensional integral because it will provide the natural point of comparison for the cell-coordinate moment formula of Section 6.

On $C_{12}$,

$$u(t,\theta)=t\,e_0+\sqrt{1-t^2}(\cos\theta\, e_1+\sin\theta\, e_2).$$

The coordinate vectors satisfy

$$\left\|\frac{\partial u}{\partial t}\right\|=\frac{1}{\sqrt{1-t^2}}, \left\|\frac{\partial u}{\partial \theta}\right\|=\sqrt{1-t^2},$$

and are orthogonal. Hence

$$dA=dt\,d\theta.$$

On $0\le\theta\le\pi/2$, put

$$g(\theta)=\sqrt{2}\cos\theta+\sin\theta$$

and

$$t_{min}(\theta)=\frac{g(\theta)}{\sqrt{1+g(\theta)^2}}.$$

Condition (6) is equivalent to

$$t_{min}(\theta)\le t\le 1.$$

By Proposition 4.1, $W_T\in[t_0,1]$; hence $W_T$ is bounded and bounded away from zero, so every real power is integrable.

**Proposition 4.6 (all real power moments).** For every $r\in\mathbb{R}\setminus\{-1\}$,

$$E W_T^r = \frac{12}{\pi(r+1)} \int_0^{\pi/2} \left[ 1 - \left( \frac{g(\theta)}{\sqrt{1+g(\theta)^2}} \right)^{r+1} \right] d\theta . \quad (11)$$

At $r = -1$, the continuous limiting form is

$$E W_T^{-1} = \frac{12}{\pi} \int_0^{\pi/2} \log \left( \frac{\sqrt{1+g(\theta)^2}}{g(\theta)} \right) d\theta .$$

*Proof.* The twelve max-min cells are congruent, and the four quadrants of the $\theta$-circle give equal contributions. Therefore, for any real $r$,

$$E W_T^r = \frac{12}{4\pi} 4 \int_0^{\pi/2} \int_{t_{min}(\theta)}^{1} t^r \, dt \, d\theta .$$

For $r \neq -1$, evaluation of the inner integral gives (11). For $r = -1$,

$$\int_{t_{min}(\theta)}^{1} \frac{dt}{t} = -\log t_{min}(\theta) ,$$

which gives the logarithmic expression. ▫

Taking $r = 1$ and $r = 2$ in (11) gives

$$E W_T = \frac{3}{2\pi} \arccos\left( -\frac{1}{3} \right)$$

and

$$E W_T^2 = \frac{1}{3} \left( 1 + \frac{3+\sqrt{3}}{\pi} \right) ,$$

respectively, in agreement with Finch [8]. The complete hand evaluations are given in Appendix A.1-A.2. They use the cell integral (11), not the branchwise density (8), and therefore provide a separate exact check of the spherical-cell calculation.

## 4.7. Cuboctahedral consequence and scaling

The difference body itself is centrally symmetric, so its width law follows directly.

**Corollary 4.7 (unit-edge cuboctahedral width law).** Let

$$Q = T - T .$$

Then $Q$ is the unit-edge cuboctahedron with vertex set

$$\frac{1}{\sqrt{2}} \{ (\pm 1, \pm 1, 0), (\pm 1, 0, \pm 1), (0, \pm 1, \pm 1) \} ,$$

and, for every $u \in S^2$,

$$w_Q(u) = 2 w_T(u) .$$

Consequently,

$$W_Q = 2W_T,$$

$$supp(W_Q) = [\sqrt{2}, 2],$$

and

$$f_Q(s) = \frac{1}{2} f_T\left(\frac{s}{2}\right).$$

*Proof.* Proposition 4.1 gives the displayed vertex set. Adjacent vertices, for example $(1,1,0)/\sqrt{2}$ and $(1,0,1)/\sqrt{2}$, differ by $(0,1,-1)/\sqrt{2}$, of norm one; hence the cuboctahedron has unit edge length. Moreover,

$$h_Q(u) = h_{T-T}(u) = w_T(u).$$

Since $Q = -Q$,

$$w_Q(u) = 2h_Q(u) = 2w_T(u).$$

The probabilistic statements follow by scaling the random variable. ▫

To the best of the author’s knowledge, an explicit probability distribution for the width of the unit-edge cuboctahedron has not previously been recorded. This statement is confined to the explicit width law and makes no claim concerning the classical geometry of the cuboctahedron.

More generally, if $\lambda > 0$, then

$$W_{\lambda T} = \lambda W_T, f_{\lambda T}(t) = \frac{1}{\lambda} f_T\left(\frac{t}{\lambda}\right).$$

Thus the unit-edge formula determines the width law of every regular tetrahedron.

Section 5 turns from tetrahedral width to brightness.

## 5. Tetrahedral brightness and the rhombic-dodecahedral width law

Tetrahedral brightness is the cosine transform of four facet normals and, by Proposition 2.5, pointwise the width of a four-generator zonotope. We identify that zonotope, after explicit scaling, with the unit-edge rhombic dodecahedron. The same spherical co-area framework now has two linear cell families, producing a singularity pattern different from the width law. The broader projected-area lineage includes Walters [17], Vickers [15], Umhauer and Gutsch [14], and Vickers and Brown [16]; the tetrahedral case is located precisely below.

### *5.1. Surface-area measure and brightness normal form*

Let $T$ be the unit-edge regular tetrahedron of Section 4. Each facet is an equilateral triangle of area

$$A_{\triangle} = \frac{\sqrt{3}}{4}.$$

Set

$$\xi_1=(1,1,1),\xi_2=(1,-1,-1),\xi_3=(-1,1,-1),\xi_4=(-1,-1,1),$$

and write $\hat{\xi}_i=\xi_i/\sqrt{3}$ for the corresponding unit vectors. With the vertex convention of Section 4, the outward normal to the facet opposite $v_i$ is $-\hat{\xi}_i$. Thus

$$S_T=A_\triangle\sum_{i=1}^{4}\delta_{-\hat{\xi}_i},\hat{\xi}_i\cdot\hat{\xi}_j=-\frac{1}{3}(i\neq j).$$

Cauchy’s projection formula in Proposition 2.4 gives, for $u=(x,y,z)\in S^2$,

$$b_T(u)\quad=\frac{A_\triangle}{2}\sum_{i=1}^{4}\left|u\cdot\hat{\xi}_i\right|$$

**Proposition 5.1 (brightness normal form and support).** For $u=(x,y,z)\in S^2$,

$$b_T(u)=\frac{1}{4}max\{2\parallel u\parallel_\infty,\parallel u\parallel_1\}.(13)$$

Its range is

$$range(b_T)=\left[\frac{1}{2\sqrt{2}},\frac{1}{2}\right].$$

The minimum is attained precisely at the signed coordinate permutations of $(1,1,0)/\sqrt{2}$; the maximum is attained precisely at the six coordinate-axis directions. At each of the eight tetrahedral face-normal directions,

$$b_T(u)=\frac{\sqrt{3}}{4}.$$

*Proof.* A signed coordinate permutation permutes the unoriented set $\{\pm\xi_1,\ldots,\pm\xi_4\}$, so both sides of (13) are invariant under the signed-permutation group. It is therefore enough to assume

$$x\geq y\geq z\geq 0.$$

On this chamber, the sign of only one linear form in (12), namely $x-y-z$, can change. If $x\geq y+z$, the four absolute values sum to $4x$; if $x\leq y+z$, they sum to $2(x+y+z)$. Since $\parallel u\parallel_\infty=x$ and $\parallel u\parallel_1=x+y+z$, and $2x\gtreqless x+y+z$ exactly when $x\gtreqless y+z$, formula (13) follows.

Put $M=\parallel u\parallel_\infty$ and $L=\parallel u\parallel_1$. Since $\parallel u\parallel_2=1$,

$$1=\sum_{j=1}^{3}\left|u_j\right|^2\leq ML.$$

If $R=max\{2M,L\}$, then $M\leq R/2$ and $L\leq R$, so $1\leq R^2/2$ and therefore $R\geq\sqrt{2}$. Equality forces two nonzero coordinates, both of modulus $1/\sqrt{2}$, giving the stated minimizing directions. Conversely, $2M\leq 2$ and $L\leq\sqrt{3}<2$, so $R\leq 2$, with equality only when $M=1$, namely at a coordinate-axis direction.

Finally, at a tetrahedral face normal all three coordinate moduli equal $1/\sqrt{3}$, so the $l^1$ branch of (13) gives $\sqrt{3}/4$. ▫

Because $b_T$ is continuous on $S^2$ and $\sigma$ has full support, the support of the random brightness $B_T = b_T(U)$ is exactly this range.

The extremal directions of the tetrahedral width and brightness are exactly interchanged. Proposition 4.1 shows that the six coordinate-axis directions minimize $w_T$ and that the signed coordinate permutations of $(1,1,0)/\sqrt{2}$ maximize it; Proposition 5.1 gives the reverse assignment for $b_T$. Thus the two directional observables are extremally anti-aligned on these two symmetry orbits, a geometric counterpart of their placement on opposite sides of Proposition 2.6.

### *5.2. Fundamental chamber and critical values*

Let

$$F = \{(x, y, z) \in S^2 : x \geq y \geq z \geq 0\}.$$

The signed-permutation group has order 48, leaves (13) invariant, and acts freely on the interior $x > y > z > 0$. Hence the interiors of the forty-eight images of $F$ are disjoint and their complement lies in a finite union of great circles.

Inside $F$, define

$$F_{max} = F \cap \{x \geq y + z\},\ F_{\Sigma} = F \cap \{x \leq y + z\}.$$

Introduce the spherical points

$$p_0 = (1,0,0),\ p_1 = \frac{(1,1,0)}{\sqrt{2}},\ p_2 = \frac{(2,1,1)}{\sqrt{6}},\ p_3 = \frac{(1,1,1)}{\sqrt{3}}.$$

**Lemma 5.2 (fundamental chamber and critical set).** The interiors of $F_{max}$ and $F_{\Sigma}$ are disjoint and cover the interior of $F$ apart from the switching arc $x = y + z$. On the two cells,

$$\begin{aligned} b_T(u) &= c_{max} \cdot u = \frac{x}{2}, & c_{max} &= \left(\frac{1}{2}, 0, 0\right), & \| c_{max} \| &= \frac{1}{2}, \\ b_T(u) &= c_{\Sigma} \cdot u = \frac{x+y+z}{4}, & c_{\Sigma} &= \frac{1}{4}(1,1,1), & \| c_{\Sigma} \| &= \frac{\sqrt{3}}{4}. \end{aligned} \tag{14}$$

The irredundant faces of $F_{max}$ are $y = z$, $z = 0$, and $x = y + z$, with vertices $p_0, p_1, p_2$. The irredundant faces of $F_{\Sigma}$ are $x = y$, $y = z$, and $x = y + z$, with vertices $p_1, p_2, p_3$.

If

$$\tau_0 = \frac{1}{2\sqrt{2}},\ \tau_1 = \frac{1}{\sqrt{6}},\ \tau_2 = \frac{\sqrt{3}}{4},\ \tau_3 = \frac{1}{2},$$

then

$$\tau_0 < \tau_1 < \tau_2 < \tau_3,$$

and $\{\tau_0, \tau_1, \tau_2, \tau_3\}$ is the complete active critical set for the two cell families. Each of the forty-eight signed-permutation images of $F$ contains exactly one cell of each family.

*Proof.* The two identities in (14) are the two cases of Proposition 5.1. On $F_{max}$, the inequalities $x \geq y+z$ and $z \geq 0$ imply $x \geq y$, so $x \geq y$ is redundant; the remaining three equalities meet at $p_0, p_1, p_2$. On $F_{\Sigma}$, if $z < 0$ then $y+z < y \leq x$, contradicting $x \leq y+z$; hence $z \geq 0$ is redundant, and the three essential boundaries meet at $p_1, p_2, p_3$.

We now apply Proposition 3.3. For $F_{max}$, the axial value is $\| c_{max} \| = \tau_3$. The switching face $x = y+z$ has tangency value $\tau_1$, while the other two active face tangencies occur at $\tau_3$. The three vertex values at $p_0, p_1, p_2$ are respectively $\tau_3, \tau_0, \tau_1$.

For $F_{\Sigma}$, the axial value is $\| c_{\Sigma} \| = \tau_2$. The switching face again has tangency value $\tau_1$, while the tangencies on $x = y$ and $y = z$ occur at $\tau_2$. The vertex values at $p_1, p_2, p_3$ are $\tau_0, \tau_1, \tau_2$. Each listed tangency lies on the corresponding boundary arc, so all four distinct values are active. Their order follows from the squares $1/8 < 1/6 < 3/16 < 1/4$. No further candidate is available by Proposition 3.3. ▫

The two brightness cells are represented in Figure 2 in the gnomonic coordinates

$$Y = \frac{y}{x}, Z = \frac{z}{x}.$$

Since $x > 0$ throughout the fundamental chamber, this projection is well defined and carries its great-circle boundary arcs to straight line segments. In these coordinates the switching arc $x = y+z$ becomes $Y + Z = 1$, making the two coefficient regions and their common boundary explicit.

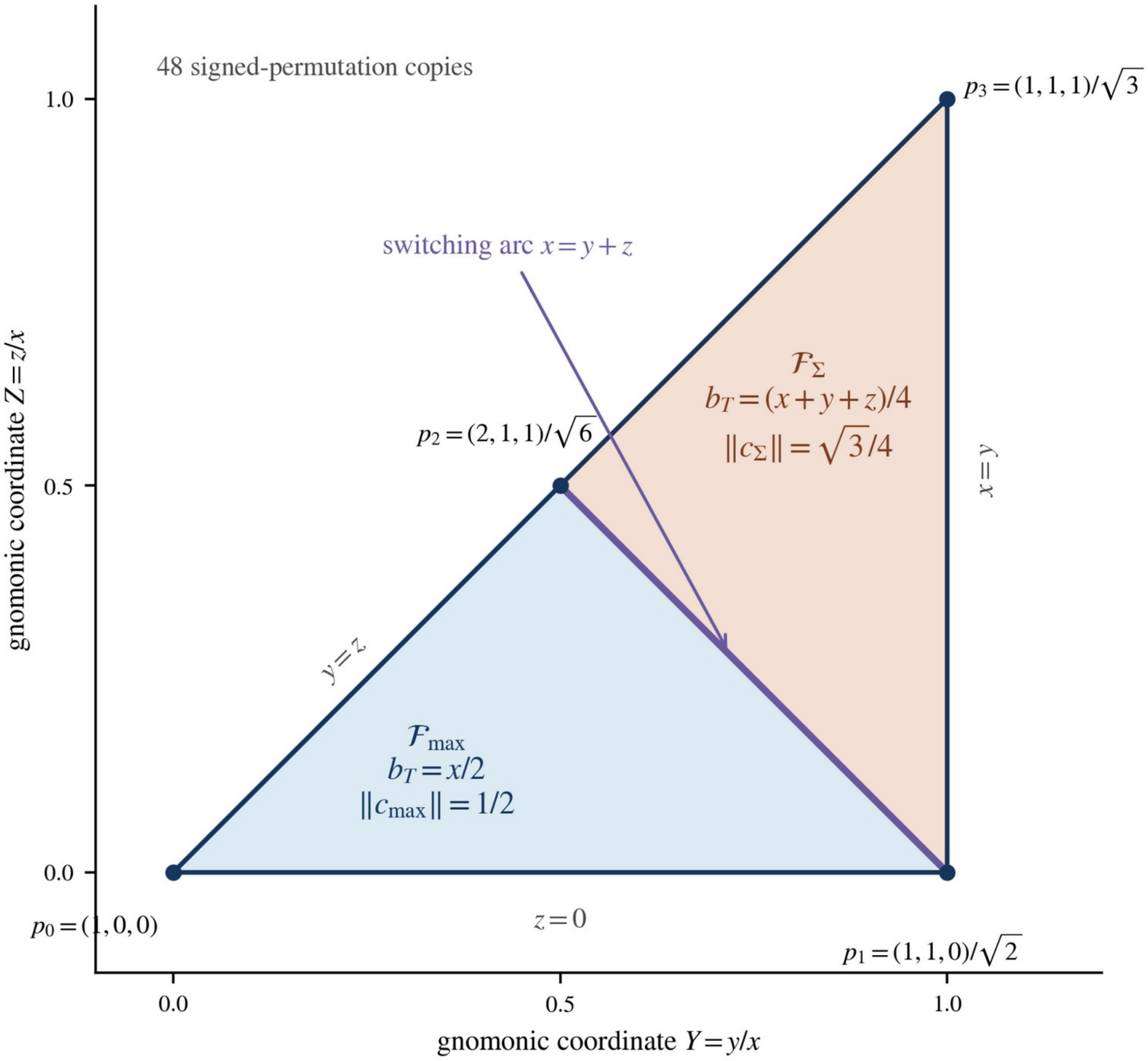

**Figure 2. Fundamental chamber for tetrahedral brightness.** Exact gnomonic image of the spherical chamber $F = \{ x \geq y \geq z \geq 0 \}$. The chamber is divided by $x = y + z$, equivalently $Y + Z = 1$, into $F_{max}$, where $b_T(u) = x/2$ and $\| c_{max} \| = 1/2$, and $F_\Sigma$, where $b_T(u) = (x + y + z)/4$ and $\| c_\Sigma \| = \sqrt{3}/4$. The marked vertices are $p_0, p_1, p_2, p_3$ from Lemma 5.2. The full spherical decomposition consists of forty-eight signed-permutation copies of this chamber.

### *5.3. Exact aperture functions*

The two cell families require different level-circle coordinates.

**Lemma 5.3 (exact aperture functions).** Define

$$\eta(t) = \arcsin\left(\frac{\sqrt{2}t}{\sqrt{1-4t^2}}\right), \delta(t) = \arccos\left(\frac{\sqrt{2}t}{\sqrt{3-16t^2}}\right),$$

for $\tau_0 \leq t \leq \tau_1$. Their endpoint values are

$$\eta(\tau_0) = \frac{\pi}{4}, \eta(\tau_1) = \frac{\pi}{2}, \delta(\tau_0) = \frac{\pi}{3}, \delta(\tau_1) = 0.$$

The angular apertures of the level circles in $F_{max}$ and $F_\Sigma$ are

$$\Theta_{max}(t)=\begin{cases}0, & t<\tau_0,\\ \eta(t)-\dfrac{\pi}{4}, & \tau_0\le t<\tau_1,\\ \dfrac{\pi}{4}, & \tau_1\le t<\tau_3,\\ 0, & t\ge\tau_3,\end{cases} \quad (15)$$

and

$$\Theta_{\Sigma}(t)=\begin{cases}0, & t<\tau_0,\\ \dfrac{\pi}{3}-\delta(t), & \tau_0\le t<\tau_1,\\ \dfrac{\pi}{3}, & \tau_1\le t<\tau_2,\\ 0, & t\ge\tau_2.\end{cases} \quad (16)$$

At the axial values $\tau_2$ and $\tau_3$ the corresponding level circle collapses to a point; the value assigned to an aperture or density representative there is immaterial.

*Proof.* On $F_{max}$, the level equation $x/2=t$ gives

$$x=2t,\ y=l_{max}(t)\cos\theta,\ z=l_{max}(t)\sin\theta,\ l_{max}(t)=\sqrt{1-4t^2}.$$

The inequalities $y\ge z\ge 0$ are equivalent to $0\le\theta\le\pi/4$. The switching condition $x\ge y+z$ becomes

$$\cos\theta+\sin\theta\le\frac{2t}{\sqrt{1-4t^2}}.$$

Since

$$\cos\theta+\sin\theta=\sqrt{2}\sin\left(\theta+\frac{\pi}{4}\right),\theta+\frac{\pi}{4}\in\left[\frac{\pi}{4},\frac{\pi}{2}\right],$$

the intermediate admissible interval is $0\le\theta\le\eta(t)-\pi/4$. Direct substitution gives $\eta(\tau_0)=\pi/4$ and $\eta(\tau_1)=\pi/2$, so the interval is empty below $\tau_0$, reaches the full chamber width $\pi/4$ at $\tau_1$, remains full until the axial collapse at $\tau_3$, and is empty thereafter. This proves (15).

For $F_{\Sigma}$, put

$$\omega_0=\frac{(1,1,1)}{\sqrt{3}},\omega_1=\frac{(1,-1,0)}{\sqrt{2}},\omega_2=\frac{(1,1,-2)}{\sqrt{6}}.$$

The level equation $(x+y+z)/4=t$ is $\omega_0\cdot u=4t/\sqrt{3}$. With

$$l_{\Sigma}(t)=\sqrt{1-\frac{16t^2}{3}},$$

write

$$u = \frac{4t}{\sqrt{3}}\omega_0 + l_\Sigma(t)(\cos\theta\,\omega_1 + \sin\theta\,\omega_2).$$

Direct expansion gives

$$x - y = \sqrt{2}\,l_\Sigma(t)\cos\theta,$$

and

$$y - z = \frac{l_\Sigma(t)}{\sqrt{2}}(-\cos\theta + \sqrt{3}\sin\theta).$$

The inequality $x \geq y$ first gives $\cos\theta \geq 0$, so we may take $-\pi/2 \leq \theta \leq \pi/2$. On this arc, $y \geq z$ is equivalent to $\tan\theta \geq 1/\sqrt{3}$. Hence $x \geq y \geq z$ is exactly

$$\frac{\pi}{6} \leq \theta \leq \frac{\pi}{2}.$$

The switching condition $x \leq y + z$ reduces to

$$\sqrt{3}\cos\theta + \sin\theta \leq \frac{2\sqrt{2}t}{\sqrt{3 - 16t^2}}.$$

Because

$$\sqrt{3}\cos\theta + \sin\theta = 2\cos\left(\theta - \frac{\pi}{6}\right),$$

and the left side decreases from $2$ to $1$ on $[\pi/6, \pi/2]$, the admissible interval has length $\pi/3 - \delta(t)$ for $\tau_0 \leq t < \tau_1$. Direct substitution gives $\delta(\tau_0) = \pi/3$ and $\delta(\tau_1) = 0$; from $\tau_1$ onward the whole chamber interval is admitted until the axial collapse at $\tau_2$. This proves (16). ▫

### 5.4. Tetrahedral brightness density

Each of the forty-eight signed-permutation chambers contains one cell from each family. By Theorem 3.1 and (14), for almost every $t$,

$$f_{B_T}(t) = \frac{48}{4\pi}\left(\frac{\Theta_{max}(t)}{1/2} + \frac{\Theta_\Sigma(t)}{\sqrt{3}/4}\right) = \frac{24}{\pi}\Theta_{max}(t) + \frac{48}{\pi\sqrt{3}}\Theta_\Sigma(t). \quad (17)$$

**Theorem 5.4 (tetrahedral brightness density).** Let $U$ be uniform on $S^2$, and put $B_T = b_T(U)$. Then $B_T$ is absolutely continuous with support $[\tau_0, \tau_3]$. A density representative is

$$f_{B_T}(t)=\begin{cases}0, & t<\tau_0,\\ \dfrac{24}{\pi}\left(\eta(t)-\dfrac{\pi}{4}\right)+\dfrac{48}{\pi\sqrt{3}}\left(\dfrac{\pi}{3}-\delta(t)\right), & \tau_0\le t<\tau_1,\\ 6+\dfrac{16}{\sqrt{3}}, & \tau_1\le t<\tau_2,\\ 6, & \tau_2\le t<\tau_3,\\ 0, & t\ge\tau_3.\end{cases} \quad (18)$$

*Proof.* Proposition 5.1 gives the support, Lemma 5.2 gives the two linear cell families and their multiplicities, and Lemma 5.3 gives their apertures. Substitution of (15)–(16) into (17) yields (18). Absolute continuity follows from Theorem 3.1. ▫

Vickers and Brown [16, p. 292; Appendix B, pp. 304–306], who take each edge of the regular tetrahedron to have unit length, derived the projected-area density used here. Theorem 5.4 gives an independent spherical-cell derivation of the same random variable. To compare the formulae directly, the identities

$$2\eta(t)-\frac{\pi}{2}=\arcsin\left(\frac{8t^2-1}{1-4t^2}\right),\pi-2\delta(t)=\arccos\left(\frac{3-20t^2}{3-16t^2}\right)$$

transform the first branch of (18) into

$$\frac{12}{\pi}\arcsin\left(\frac{8t^2-1}{1-4t^2}\right)+\frac{8\sqrt{3}}{\pi}\arccos\left(\frac{3-20t^2}{3-16t^2}\right)-\frac{8}{\sqrt{3}}. \quad (19)$$

Equation (19) is the unit-edge tetrahedral form displayed in [16, p. 292].

### 5.5. Transition structure

Let

$$c_*=6+\frac{16}{\sqrt{3}}$$

be the plateau value in the second branch of (18).

**Lemma 5.5 (global critical behavior).** The density is real-analytic on

$$(\tau_0,\tau_1),(\tau_1,\tau_2),(\tau_2,\tau_3).$$

At the lower endpoint,

$$f_{B_T}(t)=\frac{192\sqrt{2}}{\pi}(t-\tau_0)+O\left((t-\tau_0)^2\right)(t\downarrow\tau_0). \quad (20)$$

At $\tau_1$, the density is continuous and approaches the plateau through

$$f_{B_T}(t)=c_*-\frac{72\,6^{3/4}}{\pi}\sqrt{\tau_1-t}+O\left((\tau_1-t)^{3/2}\right)(t\uparrow\tau_1). \quad (21)$$

In particular,

$$f_{B_T}'(\tau_1^-)=+\infty.$$

At $\tau_2$, the density has a downward step of size $16/\sqrt{3}$; at $\tau_3$, it has a terminal downward step of size $6$. Neither step represents an atom.

*Proof.* On the first branch,

$$\eta'(t)=\frac{\sqrt{2}}{(1-4t^2)\sqrt{1-6t^2}},\delta'(t)=-\frac{\sqrt{6}}{(3-16t^2)\sqrt{1-6t^2}}.(22)$$

At $t=\tau_0$, (22) gives $\eta'(\tau_0)=4\sqrt{2}$ and $\delta'(\tau_0)=-2\sqrt{6}$. Differentiating the first branch of (18) therefore yields the coefficient in (20).

For the transition at $\tau_1$, define

$$\rho_{max}(t)=\frac{\sqrt{2}t}{\sqrt{1-4t^2}},\rho_{\Sigma}(t)=\frac{\sqrt{2}t}{\sqrt{3-16t^2}}.$$

Both equal $1$ at $\tau_1$, while

$$\rho_{max}'(\tau_1)=3\sqrt{6},\rho_{\Sigma}'(\tau_1)=9\sqrt{6}.$$

Hence

$$1-\rho_{max}(t)=3\sqrt{6}(\tau_1-t)+O\left((\tau_1-t)^2\right),$$

and

$$1-\rho_{\Sigma}(t)=9\sqrt{6}(\tau_1-t)+O\left((\tau_1-t)^2\right).$$

Using $\arccos(1-\varepsilon)=\sqrt{2\varepsilon}+O(\varepsilon^{3/2})$, we obtain

$$\frac{\pi}{2}-\eta(t)=6^{3/4}\sqrt{\tau_1-t}+O\left((\tau_1-t)^{3/2}\right),$$

and

$$\delta(t)=\sqrt{3}6^{3/4}\sqrt{\tau_1-t}+O\left((\tau_1-t)^{3/2}\right).$$

Their weighted sum gives (21). For the one-sided derivative, differentiate the first branch of (18) directly and use (22):

$$f_{B_T}'(t)=\frac{24\sqrt{2}}{\pi(1-4t^2)\sqrt{1-6t^2}}+\frac{48\sqrt{6}}{\pi\sqrt{3}(3-16t^2)\sqrt{1-6t^2}}.$$

Both summands are positive on $(\tau_0,\tau_1)$, and the common factor $(1-6t^2)^{-1/2}$ diverges as $t\uparrow\tau_1$, because $6\tau_1^2=1$. Hence $f'_{B_T}(\tau_1^-)=+\infty$ without differentiating the remainder in (21).

On $(\tau_1, \tau_2)$ both apertures are constant, so the density is $c_*$. At $\tau_2$ the $\Sigma$-family reaches its axial endpoint and its contribution $16/\sqrt{3}$ disappears; at $\tau_3$ the $max$-family reaches its axial endpoint and the remaining contribution $6$ disappears. Absolute continuity from Theorem 3.1 excludes atoms. ▫

The contrast with Section 4 is structural. There the twelve congruent max-min cells have a single coefficient norm and a single axial endpoint, producing a finite corner, a fold, and then a terminal step. Here two cell families have different coefficient norms and different axial values; the collapse of the $\Sigma$-family therefore creates the interior step at $\tau_2$. The four generator directions are in general position. Indeed, $det(\xi_1, \xi_2, \xi_3) = 4$, and since $\xi_1 + \xi_2 + \xi_3 + \xi_4 = 0$, replacing any one vector by $\xi_4$ changes the determinant only by sign; all four triple determinants are therefore $\pm 4$. The observed onset, fold, and step mechanisms are consistent with the stratified zonotopal classification in [2, Section 8]. The non-zonotopal width law of Section 4 exhibits the same local mechanism types without admitting a generator sign-cell representation.

The different ordering of the brightness transitions, including the interior jump caused by termination of the $\Sigma$-family, is displayed in the middle panel of Figure 3.

### 5.6. Distribution function and normalization

For $\tau_0 \le t \le \tau_1$, define

$$J_{max}(t) = t\eta(t) + \frac{1}{2}\arctan\sqrt{2(1-6t^2)},$$

and

$$J_{\Sigma}(t) = t\delta(t) - \frac{\sqrt{3}}{4}\arctan\left(2\sqrt{2(1-6t^2)}\right).$$

Both functions extend continuously to $t = \tau_1$, and direct differentiation on $(\tau_0, \tau_1)$ gives $J_{max}' = \eta$ and $J_{\Sigma}' = \delta$. Put

$$\Phi(t) \quad = \frac{24}{\pi}\left[J_{max}(t) - J_{max}(\tau_0) - \frac{\pi}{4}(t - \tau_0)\right]$$

**Proposition 5.6 (brightness distribution function).** The density in Theorem 5.4 is normalized, and its cumulative distribution function is

$$F_{B_T}(t) = \begin{cases} 0, & t < \tau_0, \\ \Phi(t), & \tau_0 \le t < \tau_1, \\ c_* t - 6, & \tau_1 \le t < \tau_2, \\ 6t - 2, & \tau_2 \le t < \tau_3, \\ 1, & t \ge \tau_3. \end{cases} \quad (23)$$

In particular, $\int_{\mathbb{R}} f_{B_T}(t)\,dt = 1$.

*Proof.* The derivative identities above give $\Phi' = f_{B_T}$ on the first branch and $\Phi(\tau_0) = 0$. Appendix B.1 establishes the endpoint identity

$$\Phi(\tau_1) = c_* \tau_1 - 6.$$

The next branch therefore integrates to $c_* t - 6$. At $t = \tau_2$, this equals $6\tau_2 - 2$, and the last affine branch satisfies $6\tau_3 - 2 = 1$. Hence (23) is continuous and normalized, and differentiating away from the critical values reproduces (18). ▫

### *5.7. Exact first and second brightness moments*

Finch [7, Section 2.1] studied the same random projected area, termed the chorowidth, and obtained exact closed forms for its first two moments. The following proposition records these moments in the present normalization and derives the second by a self-contained two-direction calculation in Appendix B.2.

**Proposition 5.7 (first two brightness moments).** For the unit-edge regular tetrahedron,

$$E\,B_T = \frac{\sqrt{3}}{4},$$

and

$$E\,B_T^2 = \frac{1}{16} + \frac{1}{8\pi}\left(2\sqrt{2} + \arcsin\frac{1}{3}\right). \quad (24)$$

The second expression is equivalent to Finch's form

$$\frac{1}{8} + \frac{\sqrt{2}}{4\pi} - \frac{arcsec(3)}{8\pi}.$$

*Proof.* Rotational invariance gives $E|U \cdot v| = 1/2$ for every fixed $v \in S^2$. Using (12),

$$E\,B_T = \frac{A_\triangle}{2} \sum_{i=1}^{4} E\left|U \cdot \hat{\xi}_i\right| = \frac{\sqrt{3}}{4}.$$

Appendix B.2 proves directly that for unit vectors $\omega, \zeta$ with $|\omega \cdot \zeta| = 1/3$,

$$E\left(|U \cdot \omega||U \cdot \zeta|\right) = \frac{2}{9\pi}\left(2\sqrt{2} + \arcsin\frac{1}{3}\right).$$

Since $\hat{\xi}_i \cdot \hat{\xi}_j = -1/3$ for $i \neq j$, expansion of the square in (12) gives (24). Finally,

$$\arcsin\frac{1}{3} = \frac{\pi}{2} - arcsec(3),$$

which converts (24) to Finch's expression. ▫

The equality $E\,B_T = \sqrt{3}/4 = \tau_2$ is specific to the regular tetrahedron: its mean projected area equals the brightness in a face-normal direction, which is also the interior step location of the density.

### 5.8. Rhombic-dodecahedral realization

The pointwise bridge of Proposition 2.5 now identifies the zonotope underlying the brightness calculation. With the vectors $\xi_i$ of Section 5.1,

$$A_{\triangle}\hat{\xi}_i=\frac{\xi_i}{4},$$

so the associated centered zonotope is

$$Z_T=\sum_{i=1}^{4}\left[-\frac{\xi_i}{16},\frac{\xi_i}{16}\right].(25)$$

The four unoriented generator directions are the long-diagonal directions of a cube, the classical four-generator realization of the rhombic dodecahedron [5].

**Proposition 5.8 (rhombic-dodecahedral realization and pointwise duality).** The zonotope $Z_T$ in (25) has the fourteen vertices

$$\left(\pm\frac{1}{4},0,0\right),\left(0,\pm\frac{1}{4},0\right),\left(0,0,\pm\frac{1}{4}\right),$$

and

$$\left(\pm\frac{1}{8},\pm\frac{1}{8},\pm\frac{1}{8}\right),$$

where the eight sign choices in the second line are independent. Its edge length is $\sqrt{3}/8$. Hence

$$R=\frac{8}{\sqrt{3}}Z_T(26)$$

is a unit-edge rhombic dodecahedron. Moreover,

$$w_R(u)=\frac{8}{\sqrt{3}}b_T(u)\left(u\in S^2\right).$$

*Proof.* Every exposed vertex of (25) has the form

$$\frac{1}{16}\sum_{i=1}^{4}\varepsilon_i\xi_i,\varepsilon_i\in\{-1,1\}.$$

Because $\xi_1+\xi_2+\xi_3+\xi_4=0$, the two all-equal sign choices give the origin. A pattern with one sign different from the other three gives the eight points $(\pm1/8,\pm1/8,\pm1/8)$; a pattern with two positive and two negative signs gives the six axial points displayed above. Each is exposed by a suitable direction, so these are exactly the vertices. An edge joins, for example, $(1/4,0,0)$ to $(1/8,1/8,1/8)$, and its length is $\sqrt{3}/8$, proving (26).

By Proposition 2.5, $w_{Z_T}(u)=b_T(u)$ pointwise. Width is homogeneous under dilation, so (26) gives the stated identity for $R$. ▫

### *5.9. The unit-edge rhombic-dodecahedral width law*

**Theorem 5.9 (unit-edge rhombic-dodecahedral width law).** Let $W_R = w_R(U)$ for the unit-edge rhombic dodecahedron $R$ of Proposition 5.8. Then

$$W_R = \frac{8}{\sqrt{3}} B_T$$

pointwise for the same random direction $U$. In particular, $W_R$ is absolutely continuous. Put

$$\kappa_0 = \frac{2\sqrt{6}}{3}, \kappa_1 = \frac{4\sqrt{2}}{3}, \kappa_2 = 2, \kappa_3 = \frac{4\sqrt{3}}{3},$$

and, for $\kappa_0 \le s \le \kappa_1$, define

$$\eta_R(s) = \arcsin\left(\frac{\sqrt{6}s}{2\sqrt{16-3s^2}}\right), \delta_R(s) = \arccos\left(\frac{\sqrt{2}s}{4\sqrt{4-s^2}}\right).$$

The support of $W_R$ is $[\kappa_0, \kappa_3]$, and a density representative is

$$f_R(s) = \begin{cases} 0, & s < \kappa_0, \\ \frac{3\sqrt{3}}{\pi}\left(\eta_R(s) - \frac{\pi}{4}\right) + \frac{6}{\pi}\left(\frac{\pi}{3} - \delta_R(s)\right), & \kappa_0 \le s < \kappa_1, \\ 2 + \frac{3\sqrt{3}}{4}, & \kappa_1 \le s < \kappa_2, \\ \frac{3\sqrt{3}}{4}, & \kappa_2 \le s < \kappa_3, \\ 0, & s \ge \kappa_3. \end{cases} \quad (27)$$

For the first CDF branch, write

$$\begin{aligned} J_{max,R}(s) &= \frac{\sqrt{3}s}{8}\eta_R(s) + \frac{1}{2}\arctan\sqrt{2\left(1 - \frac{9s^2}{32}\right)}, \\ J_{\Sigma,R}(s) &= \frac{\sqrt{3}s}{8}\delta_R(s) - \frac{\sqrt{3}}{4}\arctan\left(2\sqrt{2\left(1 - \frac{9s^2}{32}\right)}\right), \end{aligned}$$

and

$$\Phi_R(s) = \frac{24}{\pi}\left[J_{max,R}(s) - J_{max,R}(\kappa_0) - \frac{\pi\sqrt{3}}{32}(s - \kappa_0)\right]$$

Then

$$F_R(s)=\begin{cases} 0, & s<\kappa_0, \\ \Phi_R(s), & \kappa_0\le s<\kappa_1, \\ \left(2+\frac{3\sqrt{3}}{4}\right)s-6, & \kappa_1\le s<\kappa_2, \\ \frac{3\sqrt{3}}{4}s-2, & \kappa_2\le s<\kappa_3, \\ 1, & s\ge\kappa_3. \end{cases} \quad (28)$$

*Proof.* Proposition 5.8 gives $W_R=(8/\sqrt{3})B_T$. Thus $f_R(s)=(\sqrt{3}/8)f_{B_T}(\sqrt{3}s/8)$ and $F_R(s)=F_{B_T}(\sqrt{3}s/8)$. Substitution into (18) gives the nonzero branches in (27), while substitution into (23) gives (28). The four critical values are the images of $\tau_0,\tau_1,\tau_2,\tau_3$ under the same dilation. ▫

To the best of the author’s knowledge, the corresponding explicit width law for the unit-edge rhombic dodecahedron has not previously been stated. The analytic density in (27) is the dilation of the tetrahedral projected-area density of Vickers and Brown [16, p. 292] under the pointwise identity of Proposition 5.8; the assertion here therefore concerns the geometric identification and the explicit unit-edge width-law formulation, not the underlying tetrahedral projected-area density. As a direct consequence of Proposition 5.7,

$$E W_R=\frac{8}{\sqrt{3}}E B_T=2.$$

The moment calculus of Section 6 uses both the non-zonotopal width law and the zonotopal brightness representation.

## 6. Functional and moment calculus across the zonotopal divide

The spherical normal-fan calculus of Theorems 3.1 and 3.4 applies to every polytope width. On the zonotopal side, the cosine-transform representation adds a multilinear calculus in the generating measure. We formulate the common functional calculus first and then the stronger second- and third-order information supplied by absolute-product kernels; the overlap for $B_T$ and $W_R$ provides exact consistency checks.

The ordinary width moments also belong to two established $L_p$ frameworks. For every full-dimensional convex body $K$ and every $r\neq 0$, the one-dimensional $r$-mean projection functional of Giannopoulos, Markessinis and Tsolomitis [27, §2(i)] satisfies

$$W_{[1,r]}(K)=\left(\int_{G_{3,1}} vol_1(P_F K)^r d\nu_{3,1}(F)\right)^{1/r}=\left(E W_K^r\right)^{1/r}.$$

Here $\nu_{3,1}$ is Haar probability measure on $G_{3,1}$. For $r\ge 1$, Yuan, Leng and Cheung [28, Introduction] use the $r$-homogeneous normalization

$$\omega_r(L)=2\int_{S^2} h_L(u)^r d\sigma(u).$$

Since $h_{K-K}=w_K$, their convention gives

$$E\,W_K^r=\frac{1}{2}\omega_r(K-K).$$

The exact moment formulae below therefore evaluate these established geometric functionals for the bodies considered here, rather than defining a separate family of invariants.

### 6.1. Cell-coordinate functional calculus

Let $F:S^2\to\mathbb{R}$ satisfy the hypotheses of Theorem 3.1. For $C\in C_2$, write

$$F(u)=a_C\cdot u,\ D_C=\|a_C\|,\ \hat{a}_C=\frac{a_C}{D_C}$$

on $\int C$, and choose an orthonormal basis $(\hat{a}_C, e_{C,1}, e_{C,2})$ of $\mathbb{R}^3$. For $|t|<D_C$, define

$$u_C(t,\theta)=\frac{t}{D_C}\hat{a}_C+\sqrt{1-\frac{t^2}{D_C^2}}(\cos\theta\, e_{C,1}+\sin\theta\, e_{C,2})$$

and

$$\Omega_C=\{(t,\theta):|t|<D_C,u_C(t,\theta)\in\int C\}.$$

**Theorem 6.1 (cell-coordinate functional calculus).** If $g:\mathbb{R}\to[0,\infty]$ is Borel measurable, then

$$E\,g(F(U))=\frac{1}{4\pi}\sum_{C\in C_2}\frac{1}{D_C}\iint_{\Omega_C}g(t)\,dt\,d\theta,(29)$$

with equality in $[0,\infty]$. The same identity holds for real- or complex-valued Borel $g$ whenever $g(F(U))$ is integrable, equivalently whenever the right-hand side of (29) with $|g|$ in place of $g$ is finite.

In particular, for $r\geq 0$ the choice $g(t)=|t|^r$ requires no lower support bound. If $F\geq m>0$, then $g(t)=t^r$ is admissible for every real $r$, including negative powers and the logarithmic limits obtained at $r=-1$ below.

*Proof.* On $\int C$, one has $a_C\cdot u_C(t,\theta)=t$. Direct differentiation gives

$$\left\|\frac{\partial u_C}{\partial t}\times\frac{\partial u_C}{\partial\theta}\right\|=\frac{1}{D_C},$$

so spherical area on the cell is $dA=D_C^{-1}dt\,d\theta$. Since $d\sigma=dA/(4\pi)$, Tonelli's theorem gives (29) for $g\geq 0$ after summing over the pairwise-disjoint cell interiors. Their complement has spherical measure zero by the hypothesis of Theorem 3.1. The integrable signed or complex case follows by applying the same identity to $|g|$ and then to real and imaginary parts. ▫

Equation (29) is simultaneously a distributional and transform calculus. Indicator functions recover distribution functions; power functions give moments. For example, because $F$ is bounded on the finite cell decomposition, its Laplace transform is represented for every real $\lambda$ by

$$E\, e^{-\lambda F(U)} = \frac{1}{4\pi} \sum_{C \in C_2} \frac{1}{D_C} \iint_{\Omega_C} e^{-\lambda t}\, dt\, d\theta .$$

Oscillatory choices give the characteristic function in the same coordinates. Thus Theorem 3.1 and the moment formulas below are specializations of one spherical change of variables. For a polytope width, Theorem 3.4 makes every term in (29) canonical by indexing the cells with vertices of the difference body.

### *6.2. All real tetrahedral brightness moments*

The two chamber families of Lemma 5.3 give an explicit application of Theorem 6.1 to $B_T$. Define

$$t_{max}(\theta) = \frac{\cos\theta + \sin\theta}{2\sqrt{1+(\cos\theta+\sin\theta)^2}}, 0 \le \theta \le \frac{\pi}{4},$$

and

$$t_\Sigma(\theta) = \frac{\sqrt{3}\left(\sqrt{3}\cos\theta + \sin\theta\right)}{2\sqrt{2+4\left(\sqrt{3}\cos\theta+\sin\theta\right)^2}}, \frac{\pi}{6} \le \theta \le \frac{\pi}{2}.$$

The two coefficient norms are $D_{max} = 1/2$ and $D_\Sigma = \sqrt{3}/4$.

**Proposition 6.2 (all real tetrahedral brightness moments).** For every $r \in \mathbb{R} \setminus \{-1\}$,

$$E\, B_T^r = \frac{12}{\pi(r+1)}[ \quad 2\int_0^{\pi/4} \left\{ \left(\frac{1}{2}\right)^{r+1} - t_{max}(\theta)^{r+1} \right\} d\theta$$

At $r = -1$,

$$E\, B_T^{-1} = \frac{12}{\pi}[ \quad 2\int_0^{\pi/4} \log \frac{1/2}{t_{max}(\theta)}\, d\theta$$

Moreover, for every real $r$,

$$E\, W_R^r = \left(\frac{8}{\sqrt{3}}\right)^r E\, B_T^r .$$

*Proof.* In the $F_{max}$ family, the switching inequality from Lemma 5.3 is

$$2t \ge \sqrt{1-4t^2}\, c_{max}(\theta), c_{max}(\theta) = \cos\theta + \sin\theta .$$

Both sides are nonnegative on the chamber. Squaring gives

$$4t^2\left(1 + c_{max}(\theta)^2\right) \ge c_{max}(\theta)^2,$$

hence $t \ge t_{max}(\theta)$. Similarly, in the $F_\Sigma$ family, put

$$c_\Sigma(\theta) = \sqrt{3}\cos\theta + \sin\theta .$$

The switching inequality is equivalent after squaring to

$$3c_{\Sigma}(\theta)^2 \le 8t^2 + 16c_{\Sigma}(\theta)^2 t^2,$$

and therefore to $t \ge t_{\Sigma}(\theta)$. Direct substitution gives

$$t_{max}(0) = t_{\Sigma}(\pi/2) = \tau_0, t_{max}(\pi/4) = t_{\Sigma}(\pi/6) = \tau_1.$$

Each of the forty-eight signed-permutation images of the fundamental chamber contains one cell of each family. Applying Theorem 6.1 and integrating $t^r$ first gives (30), with the logarithmic expression when $r = -1$. The final identity follows from $W_R = (8/\sqrt{3}) B_T$. ▫

The functions $s \mapsto E W_T^{s-1}$ and $s \mapsto E B_T^{s-1}$ are the Mellin transforms of the two positive laws. Because both supports are compact and bounded away from zero, these transforms extend to entire functions of $s$. Equations (11) and (30) therefore provide explicit Mellin representations on opposite sides of Proposition 2.6. In particular, negative moments belong to the same transform family rather than forming an auxiliary extension of the integer moment sequence.

### 6.3. Mean chords and the negative first brightness moment

The case $r = -1$ has an intrinsic integral-geometric interpretation. Let $K \subset \mathbb{R}^3$ be a convex body with positive volume. For $u \in S^2$ and $y \in K \mid u^{\perp}$, let $l_{K,u}(y)$ denote the length of the chord $K \cap (y + \mathbb{R}u)$.

**Proposition 6.3 (direction-averaged mean chord and a width-brightness bound).** The mean length of chords parallel to $u$, with base point uniformly distributed over $K \mid u^{\perp}$, is

$$\bar{l}_K(u) = \frac{V(K)}{b_K(u)}.$$

Consequently,

$$E\bar{l}_K(U) = V(K) E B_K^{-1} \ge \frac{4V(K)}{S(K)}, (31)$$

with equality if and only if $K$ has constant brightness. In every direction,

$$b_K(u) w_K(u) \ge V(K).$$

*Proof.* Fubini's theorem in the direction $u$ gives

$$V(K) = \int_{K \mid u^{\perp}} l_{K,u}(y)\, dy.$$

Division by the projected area $b_K(u)$ gives the first identity. Averaging over $U$ yields the first equality in (31). Cauchy's projection formula implies $E B_K = S(K)/4$, and Jensen's inequality for $x \mapsto x^{-1}$ gives

$$E B_K^{-1} \ge \frac{1}{E B_K} = \frac{4}{S(K)}.$$

Equality in Jensen holds exactly when $B_K$ is almost surely constant; continuity of $b_K$ then makes it constant on $S^2$. Finally, every chord parallel to $u$ has length at most $w_K(u)$, so $\bar{l}_K(u) \le w_K(u)$. Substituting $\bar{l}_K(u) = V(K)/b_K(u)$ gives the pointwise product inequality. ▫

The equality condition in (31) is deliberately stated as constant brightness rather than spherical rigidity. Constant-brightness bodies form a classical rigidity problem; Howard [18] proves that the additional assumption of constant width forces a Euclidean ball. For the tetrahedron, Proposition 6.2 makes the left-hand side of (31) explicitly computable from the same chamber geometry that determines the full brightness law.

### *6.4. Exact direction-averaged mean-chord law of the tetrahedron*

For the unit-edge regular tetrahedron put

$$V_T = \frac{1}{6\sqrt{2}}, L_T = \bar{l}_T(U) = \frac{V_T}{B_T}.$$

**Corollary 6.4 (exact tetrahedral mean-chord law).** The random variable $L_T$ is absolutely continuous with support

$$\left[\frac{1}{3\sqrt{2}}, \frac{1}{3}\right].$$

Its four mapped critical values, in increasing order, are

$$\lambda_0 = \frac{1}{3\sqrt{2}}, \lambda_1 = \frac{\sqrt{6}}{9}, \lambda_2 = \frac{\sqrt{3}}{6}, \lambda_3 = \frac{1}{3}.$$

For $\lambda_0 \le l \le \lambda_3$,

$$F_{L_T}(l) = 1 - F_{B_T}\left(\frac{V_T}{l}\right), f_{L_T}(l) = \frac{V_T}{l^2} f_{B_T}\left(\frac{V_T}{l}\right),$$

with $F_{L_T} = 0$ below the support and $F_{L_T} = 1$ above it. For every real $r$,

$$E\, L_T^r = V_T^r\, E\, B_T^{-r}.$$

Thus Theorem 5.4 and Proposition 5.6 give a complete closed-form density and CDF for $L_T$, while Proposition 6.2 gives all of its real power moments.

*Proof.* The map $x \mapsto V_T / x$ is smooth and strictly decreasing on the positive support of $B_T$. The support and critical values are therefore the images of $\tau_3, \tau_2, \tau_1, \tau_0$, respectively. The CDF and density follow from the one-dimensional monotone change of variables, and the moment identity is immediate from $L_T = V_T / B_T$. ▫

The square-root and jump mechanisms of the brightness law transfer under this smooth change of variable. The corollary therefore supplies an exact law for a direction-dependent chord statistic, not merely its averaged value.

### 6.5. Polarized moments and Minkowski composition

For $k \geq 1$ and $n_1, \ldots, n_k \in S^2$, define

$$I_k(n_1, \ldots, n_k) = E \prod_{j=1}^{k} |U \cdot n_j|.$$

For finite positive Borel measures $\nu_1, \ldots, \nu_k$ on $S^2$, define

$$P_k(\nu_1, \ldots, \nu_k) = E \prod_{j=1}^{k} A_{\nu_j}(U).$$

**Lemma 6.5 (polarized expansion and Minkowski composition).** For finite positive Borel measures $\nu_1, \ldots, \nu_k$,

$$P_k(\nu_1, \ldots, \nu_k) = \int_{(S^2)^k} I_k(n_1, \ldots, n_k)\, d\nu_1(n_1) \cdots d\nu_k(n_k). \quad (32)$$

The functional $P_k$ is symmetric and separately linear. If $Z_1, \ldots, Z_m$ are zonoids with unique even generating measures $\mu_1, \ldots, \mu_m$, and $\lambda_1, \ldots, \lambda_m \geq 0$, then for $Z = \lambda_1 Z_1 + \cdots + \lambda_m Z_m$,

$$E W_Z^k = \sum_{\alpha_1 + \cdots + \alpha_m = k} \binom{k}{\alpha_1, \ldots, \alpha_m} \prod_{j=1}^{m} \lambda_j^{\alpha_j} P_k\left(\mu_1^{[\alpha_1]}, \ldots, \mu_m^{[\alpha_m]}\right), \quad (33)$$

where $\mu_j^{[\alpha_j]}$ denotes $\alpha_j$ repeated arguments.

*Proof.* Tonelli's theorem applied to the product of the nonnegative cosine-transform integrals gives (32). Symmetry and separate linearity follow immediately. Generating measures are additive under Minkowski addition and homogeneous under positive dilation, so $\mu_Z = \sum_j \lambda_j \mu_j$. Expanding $P_k(\mu_Z, \ldots, \mu_Z)$ by multilinearity gives (33). ▫

The expansion (33) is the composition rule that makes $P_k$ a polarized moment calculus: moments of Minkowski sums decompose into mixed directional moments of their summands.

### 6.6. Zonoidal and brightness moment identities

**Corollary 6.6 (zonoidal and brightness moment identities).** Let $Z \subset \mathbb{R}^3$ be a zonoid with unique even generating measure $\mu_Z$, and let $K \subset \mathbb{R}^3$ be a convex body with surface-area measure $S_K$. For every positive integer $k$,

$$E W_Z^k = \int_{(S^2)^k} I_k \, d\mu_Z^{\otimes k}, \quad E B_K^k = \frac{1}{2^k} \int_{(S^2)^k} I_k \, d S_K^{\otimes k}.$$

For a polyhedron $P$,

$$E B_P^k = \frac{1}{2^k} \sum_{F_1, \ldots, F_k} A_{F_1} \cdots A_{F_k} I_k(n_{F_1}, \ldots, n_{F_k}).$$

At $k=1$ the brightness identity is Cauchy's classical mean-projection formula

$$E\,B_K=\frac{S(K)}{4}.$$

*Proof.* Apply Proposition 2.4 and Lemma 6.5 with $\nu=\mu_Z$ and $\nu=S_K/2$, respectively. The polyhedral formula follows from the atomic surface-area measure. ▫

### *6.7. Second-order kernel, variance, and inversion*

For $n_1,n_2\in S^2$, the quantity $I_2(n_1,n_2)$ depends only on

$$\rho=|n_1\cdot n_2|\in[0,1],$$

and we write $I_2(\rho)$ accordingly. The corresponding bivariate Gaussian absolute moment is classical; see Nabeya [11, p. 3]. Under Lemma 6.9 with $k=2$, the Gaussian formula is multiplied by the spherical radial factor $1/3$. We retain a direct spherical derivation because strict monotonicity and the resulting inversion are central here.

**Theorem 6.7 (second-order kernel, variance, and inversion).** The second-order spherical kernel is

$$I_2(\rho)=\frac{2}{3\pi}\left(\sqrt{1-\rho^2}+\rho\arcsin\rho\right),(34)$$

and

$$I_2'(\rho)=\frac{2}{3\pi}\arcsin\rho.$$

Hence $I_2$ is strictly increasing and maps $[0,1]$ bijectively onto $[2/(3\pi),1/3]$.

Define

$$k_2(\rho)=I_2(\rho)-\frac{1}{4}.$$

If $Z\subset\mathbb{R}^3$ is a zonoid with unique even generating measure $\mu_Z$, and $K\subset\mathbb{R}^3$ is a convex body, then

$$\begin{aligned}Var\,W_Z&=\int_{S^2}❑\int_{S^2}k_2(|n\cdot n'|)\,d\mu_Z(n)\,d\mu_Z(n'),\\ Var\,B_K&=\frac{1}{4}\int_{S^2}❑\int_{S^2}k_2(|n\cdot n'|)\,dS_K(n)\,dS_K(n').\end{aligned}\quad(35)$$

For a polyhedron $P$ the second identity becomes

$$Var\,B_P=\frac{1}{4}\sum_{F,F'}A_F A_{F'}I_2(|n_F\cdot n_{F'}|)-\frac{S(P)^2}{16}.$$

The kernel $(n,n')\mapsto k_2(|n\cdot n'|)$ is positive semidefinite on unoriented directions, although it takes negative values. Moreover, among full-dimensional zonoids, $Var\,W_Z=0$ if and only if $Z$ is a Euclidean ball up to translation, while $Var\,B_K=0$ if and only if $K$ has constant brightness.

*Proof.* Rotation invariance permits

$$n_1=(1,0,0),\, n_2=(\cos\vartheta,\sin\vartheta,0),\, 0\le\vartheta\le\frac{\pi}{2},$$

with $\rho=\cos\vartheta$. Writing

$$U=\left(\sqrt{1-z^2}\cos\psi,\sqrt{1-z^2}\sin\psi,z\right)$$

and $d\sigma=(4\pi)^{-1}d\psi\, dz$ separates the integral. Since

$$\int_{-1}^{1}(1-z^2)\,dz=\frac{4}{3}$$

and

$$\int_0^{2\pi}\left|\cos\psi\cos(\psi-\vartheta)\right|d\psi=2\left[\sin\vartheta+\left(\frac{\pi}{2}-\vartheta\right)\cos\vartheta\right],$$

formula (34) follows. Differentiation gives the displayed derivative because the two algebraic terms cancel. The endpoint values are $I_2(0)=2/(3\pi)$ and $I_2(1)=1/3$.

Corollary 6.6 at $k=1,2$ gives

$$E\,W_Z=\frac{1}{2}\mu_Z(S^2),\, E\,B_K=\frac{1}{4}S(K).$$

Subtracting the squares of these means from the corresponding second-moment integrals gives (35) and its polyhedral specialization. For arbitrary real coefficients $c_i$,

$$\sum_{i,j}c_i c_j k_2(|n_i\cdot n_j|)=Var\left(\sum_i c_i|U\cdot n_i|\right)\ge 0,$$

so the kernel is positive semidefinite. It takes negative values because $k_2(0)=2/(3\pi)-1/4<0$, whereas $k_2(1)=1/12>0$.

If $Var\,W_Z=0$, continuity makes $w_Z$ constant on $S^2$. After translating $Z$ so that its center is the origin, central symmetry gives $h_Z=w_Z/2$, so the centered translate is a Euclidean ball; the converse is immediate. The brightness statement follows similarly from continuity: zero variance is equivalent to constant brightness. ▫

The brightness rigidity is intentionally weaker than the width rigidity. Howard [18] proves that constant width together with constant brightness forces a Euclidean ball; Theorem 6.7 requires only constant brightness in its second zero-variance clause.

For the tetrahedral brightness and the associated unit-edge rhombic dodecahedron, Proposition 5.7 and the pointwise scaling give

$$Var\,B_T=\frac{2\sqrt{2}+\arcsin(1/3)}{8\pi}-\frac{1}{8},\, Var\,W_R=\frac{8}{3\pi}\left(2\sqrt{2}+\arcsin\frac{1}{3}\right)-\frac{8}{3}.$$

These are specializations of the general quadratic form (35), not separate tetrahedral calculations.

### *6.8. Second-order pair geometry*

For a finite positive Borel measure $\nu$ on $S^2$, define the pair-angle pushforward

$$\Pi_\nu = (|n \cdot n'|)_*(\nu \otimes \nu)$$

on $[0,1]$, and the second-order kernel-value pushforward

$$\Gamma_\nu = (I_2(|n \cdot n'|))_*(\nu \otimes \nu)$$

on $[2/(3\pi), 1/3]$.

**Corollary 6.8 (second-order pair geometry and absolute Gram recovery).** The measure $\Gamma_\nu$ determines $\Pi_\nu$ uniquely through

$$\Pi_\nu = (I_2^{-1})_* \Gamma_\nu .$$

In particular, for labelled directions $n_1, \dots, n_m \in S^2$, the collection

$$\{ I_2(n_i, n_j) : 1 \le i, j \le m \}$$

determines the absolute Gram matrix

$$(|n_i \cdot n_j|)_{i,j=1}^m .$$

*Proof.* Theorem 6.7 shows that $I_2 : [0,1] \to [2/(3\pi), 1/3]$ is a homeomorphism, so the first assertion is the pushforward identity under its inverse. Applying the same inverse entrywise gives the finite labelled statement. ▫

The absolute Gram matrix is invariant under independent sign changes $n_i \mapsto \pm n_i$. Recovering a compatible signed Gram matrix, and determining when it is unique up to this sign gauge and a global orthogonal transformation, is the reconstruction problem treated in Section 7.

### *6.9. Gaussian reduction*

Let $G \sim N(0, I_3)$ and write

$$G = R_G U, R_G = \| G \| .$$

Then $U$ is uniform on $S^2$, $R_G$ and $U$ are independent, and

$$E R_G^k = 2^{k/2} \frac{\Gamma((k+3)/2)}{\Gamma(3/2)} .$$

**Lemma 6.9 (Gaussian reduction).** For every positive integer $k$ and all $n_1, \dots, n_k \in S^2$,

$$I_k(n_1, \dots, n_k) = \frac{\sqrt{\pi}}{2^{k/2+1} \Gamma((k+3)/2)} E \prod_{j=1}^{k} |G \cdot n_j| . \quad (36)$$

The formula remains valid when the Gram matrix $(n_i \cdot n_j)_{i,j}$ is singular.

*Proof.* Since $G = R_G U$,

$$\prod_{j=1}^{k} |G \cdot n_j| = R_G^k \prod_{j=1}^{k} |U \cdot n_j|.$$

Independence of $R_G$ and $U$ gives the product of expectations, and substitution of the displayed chi moment together with $\Gamma(3/2) = \sqrt{\pi}/2$ gives (36). No inverse covariance matrix is used, so singular Gram matrices require no separate hypothesis. ▫

Equation (36) is a radial reduction, not a novelty claim for Gaussian absolute moments. Its role is to place classical Gaussian formulas and spherical directional geometry in the same normalization.

### *6.10. The three-direction kernel*

Let $\rho_{ij} = n_i \cdot n_j$ and

$$\Delta = 1 + 2\rho_{12}\rho_{13}\rho_{23} - \rho_{12}^2 - \rho_{13}^2 - \rho_{23}^2.$$

When no pair is collinear, define

$$\rho_{ij\,|\,k} = \frac{\rho_{ij} - \rho_{ik}\rho_{jk}}{\sqrt{(1-\rho_{ik}^2)(1-\rho_{jk}^2)}}.$$

**Proposition 6.10 (three-direction integral with singular extension).** If $\Delta > 0$, then

$$I_3(n_1, n_2, n_3) = \frac{1}{4\pi}\Big[\sqrt{\Delta} + (\rho_{12} + \rho_{13}\rho_{23}) \arcsin \rho_{12\,|\,3} + (\rho_{23} + \rho_{12}\rho_{13}) \arcsin \rho_{23\,|\,1}\Big]. \tag{37}$$

The expression extends continuously to $\Delta = 0$ whenever $|\rho_{12}|, |\rho_{13}|, |\rho_{23}| < 1$, with $\sqrt{\Delta} = 0$ and the limiting partial correlations interpreted as $\pm 1$. If one pair is collinear, say $|\rho_{12}| = 1$, and $r = |\rho_{13}| = |\rho_{23}|$, then

$$I_3(n_1, n_2, n_3) = \frac{1 + r^2}{8}.$$

In particular, three mutually orthogonal directions give $I_3 = 1/(4\pi)$, while three distinct tetrahedral directions give

$$I_3 = \frac{1}{36} + \frac{1}{3\pi\sqrt{3}}.$$

*Proof.* For $\Delta > 0$, the vector $(G \cdot n_1, G \cdot n_2, G \cdot n_3)$ is a nonsingular standardized trivariate normal vector with the displayed Gram matrix. Nabeya's exact $(1, 1, 1)$ absolute-moment formula [12, p. 19], followed by Lemma 6.9 with $k = 3$, gives (37).

For the rank-two boundary,

$$1-\rho_{ij\mid k}^{2}=\frac{\Delta}{\left(1-\rho_{ik}^{2}\right)\left(1-\rho_{jk}^{2}\right)}.$$

Thus each partial correlation tends to $\pm 1$ as a noncollinear Gram matrix approaches a rank-two boundary point. Nonsingular Gram configurations are dense in a neighbourhood of every such rank-two configuration, and $I_3$ is continuous in $\left(n_1, n_2, n_3\right)$ by dominated convergence. Hence the limiting expression is the unique continuous extension of the nonsingular formula.

For a collinear pair, after a sign change take $n_2=n_1$. Put $X=G\cdot n_1$ and $Y=G\cdot n_3$, with correlation $r\geq 0$. Conditional normality gives

$$E\left(X^2 \mid Y\right)=1-r^2+r^2 Y^2,$$

and therefore

$$E\left(X^2|Y|\right)=\left(1-r^2\right)E|Y|+r^2 E|Y|^3=\sqrt{\frac{2}{\pi}}\left(1+r^2\right).$$

Lemma 6.9 at $k=3$ yields $\left(1+r^2\right)/8$. The orthogonal value follows from (37) with all correlations zero; the tetrahedral value follows from $\rho_{ij}=-1/3$, $\Delta=16/27$, and $\rho_{ij\mid k}=-1/2$. ▫

Nabeya's formula supplies only the nonsingular Gaussian input. The singular spherical cases, including the entire collinear branch, are part of the present extension.

### *6.11. Tetrahedral third moments*

For the unit-edge rhombic dodecahedron of Proposition 5.8, the pointwise zonotopal representation is

$$W_R=\sum_{i=1}^{4}\left|U\cdot\hat{\xi}_i\right|, B_T=\frac{\sqrt{3}}{8}W_R.$$

For distinct tetrahedral directions, Section 5.5 gives $\hat{\xi}_i\cdot\hat{\xi}_j=-1/3$, and Proposition 6.10 gives

$$I_3\left(\hat{\xi}_i,\hat{\xi}_j,\hat{\xi}_k\right)=\frac{1}{36}+\frac{1}{3\pi\sqrt{3}}.$$

For repeated indices,

$$E\left|U\cdot\hat{\xi}_i\right|^3=\frac{1}{4}, E\left|U\cdot\hat{\xi}_i\right|^2\left|U\cdot\hat{\xi}_j\right|=\frac{5}{36}\left(i\neq j\right).$$

**Corollary 6.11 (tetrahedral and rhombic-dodecahedral third moments).** The exact third moments are

$$E B_T^3=\frac{5\sqrt{3}}{128}+\frac{3}{64\pi}, E W_R^3=\frac{20}{3}+\frac{8}{\pi\sqrt{3}}. \quad (38)$$

*Proof.* Expanding $W_R^3$ gives four pure cubes, twelve ordered unequal pairs with multinomial weight $3$, and four distinct triples with multinomial weight $6$. Hence

$$E W_R^3 \quad = 4\left(\frac{1}{4}\right) + 3 \cdot 12\left(\frac{5}{36}\right) + 6 \cdot 4\left(\frac{1}{36} + \frac{1}{3\pi\sqrt{3}}\right)$$

Multiplication by $\left(\sqrt{3}/8\right)^3$ gives the first identity in (38). ▫

The value in (38) agrees with Proposition 6.2 at $r = 3$, so the cell and polarized calculi meet exactly.

### 6.12. Spectral decomposition of the second-order variance

The quadratic form in Theorem 6.7 has a canonical harmonic diagonalization. We use spherical harmonics orthonormal in $L^2\left(S^2, \sigma\right)$. For a finite signed Borel measure $v$ define

$$\hat{v}_{lm} = \int_{S^2} \overline{Y_{lm}(n)}\, dv(n).$$

The spectral action of the spherical cosine transform is classical; see, for example, Ournycheva and Rubin [19]. The normalization used here is derived directly below.

**Theorem 6.12 (spectral variance decomposition).** Put

$$a_l = \frac{1}{2}\int_{-1}^{1} |t| P_l(t)\, dt. (39)$$

Then $a_0 = 1/2$, $a_l = 0$ for odd $l$, and, for $j \geq 1$,

$$a_{2j} = (-1)^{j-1} \frac{(2j-3)!!}{(2j+2)!!}.$$

In particular,

$$a_2 = \frac{1}{8}, a_4 = -\frac{1}{48}, a_6 = \frac{1}{128}, a_8 = -\frac{1}{256}.$$

If $Z$ is a zonoid with unique even generating measure $\mu_Z$, then

$$Var\, W_Z = \sum_{j \geq 1} a_{2j}^2 \sum_{m=-2j}^{2j} \left|\widehat{\mu}_{Z2j,m}\right|^2. (40)$$

If $K$ is a convex body, then

$$Var\, B_K = \frac{1}{4}\sum_{j \geq 1} a_{2j}^2 \sum_{m=-2j}^{2j} \left|\widehat{(S_K)}_{e2j,m}\right|^2.$$

Equivalently, the covariance kernel of Theorem 6.7 has the $L^2$ zonal expansion

$$k_2(|n \cdot n'|) = \sum_{j \geq 1} (4j+1) a_{2j}^2 P_{2j}(n \cdot n'). (41)$$

For every fixed even $L \geq 2$, (40) yields the finite-band stability estimate

$$\sum_{\substack{2\le l\le L\\ l\,even}} \sum_{m=-l}^{l} \left|\widehat{\mu}_{Zlm}\right|^2 \le \frac{Var\,W_Z}{\min\limits_{\substack{2\le l\le L\\ l\,even}} \left|a_l\right|^2}.$$

There is no corresponding uniform all-degree estimate of this form, because $a_{2j} \to 0$.

*Proof.* Funk–Hecke with normalized spherical measure gives, for every spherical harmonic $Y_{lm}$,

$$\int_{S^2} |u\cdot n| Y_{lm}(n)\, d\sigma(n) = a_l Y_{lm}(u),$$

with $a_l$ as in (39). Odd multipliers vanish by parity. For even degree, the Legendre recurrence and integration over $[0,1]$ give the displayed double-factorial formula. Hence, in $L^2(S^2,\sigma)$,

$$A_v = \sum_{l,m} a_l \hat{v}_{lm} Y_{lm}.$$

For the even measure $\mu_Z$, the degree-zero term is the mean width. Parseval's identity therefore gives (40) after subtracting the mean. Since $B_K = A_{S_K/2} = A_{(S_K)_e/2}$, the brightness formula follows in the same way.

The addition theorem in the present normalization is

$$\sum_{m=-l}^{l} Y_{lm}(n)\overline{Y_{lm}(n')} = (2l+1) P_l(n\cdot n'),$$

which converts the quadratic form into (41). The finite-band estimate follows from (40) by bounding the finitely many retained multipliers below by their minimum absolute value. ▫

Theorem 6.12 makes the zero-variance statements of Theorem 6.7 structural. For a zonoid, vanishing variance removes every even harmonic component of $\mu_Z$ above degree zero, so the generating measure is uniform and the centered body is a ball. For brightness, only the even part of $S_K$ is visible, and zero variance removes its nonconstant even modes; this is exactly the constant-brightness condition. The spectral form also resolves the apparent paradox that $k_2$ takes negative values while defining a positive-semidefinite quadratic form: its harmonic coefficients in (41) are nonnegative squares.

Appendix C specializes the spectral identity (40) to finite weighted direction systems. It identifies the degree-two variance component with the squared Frobenius norm of the traceless part of the weighted frame operator. The resulting tight-frame criterion is applied to the tetrahedral normal system and the six generators of the regular truncated octahedron.

## 7. Finite normal tomography

The second-order theory of Section 6 recovers the absolute Gram data of a labelled family of directions but cannot distinguish the independent sign changes $n_i \mapsto \pm n_i$. The third-order kernel contains the additional information required to resolve this ambiguity on an explicitly characterized class of finite configurations. The reconstruction problem is therefore naturally divided into two stages: pair data recover the magnitudes of the Gram entries, while triple data recover the gauge-invariant cycle signs needed to determine their signs.

Let

$$n_1, \ldots, n_m \in S^2$$

be labelled directions and write

$$\rho_{ij} = n_i \cdot n_j .$$

The associated Gram matrix is denoted by

$$G = \left(\rho_{ij}\right)_{i,j=1}^{m} .$$

Since the observables considered in Section 6 are unchanged under independent reversals of the $n_i$, the natural reconstruction object is the equivalence class of $G$ under diagonal sign conjugacy.

### *7.1. Sign gauge and cycle data*

For

$$\varepsilon = \left(\varepsilon_1, \ldots, \varepsilon_m\right) \in \{\pm 1\}^m ,$$

put

$$D_\varepsilon = diag\left(\varepsilon_1, \ldots, \varepsilon_m\right).$$

Replacing $n_i$ by $\varepsilon_i n_i$ transforms the Gram matrix according to

$$G \mapsto D_\varepsilon G D_\varepsilon .$$

This is the switching operation in the signed-graph framework of Zaslavsky [21].

Define the **correlation graph** $\Gamma = \Gamma\left(n_1, \ldots, n_m\right)$ on the vertex set $\{1, \ldots, m\}$ by

$$ij \in E\left(\Gamma\right) \Leftrightarrow \rho_{ij} \neq 0 .$$

Each edge carries the sign

$$s_{ij} = sgn\, \rho_{ij} .$$

For a cycle

$$C = i_1 i_2 \cdots i_r i_1$$

define

$$\chi\left(C\right) = \prod_{k=1}^{r} s_{i_k i_{k+1}}, i_{r+1} = i_1 .$$

**Proposition 7.1 (sign gauge and cycle data).** Two sign assignments on the same correlation graph are related by diagonal sign conjugacy if and only if they have the same sign on every cycle. Consequently, once the absolute entries $\left|\rho_{ij}\right|$ are fixed, the cycle signs are exactly the gauge-invariant discrete information contained in the signed Gram matrix.

*Proof.* Diagonal sign conjugacy multiplies the sign of an edge $ij$ by $\varepsilon_i \varepsilon_j$. Around a cycle every vertex factor occurs twice, so every cycle sign is invariant.

Conversely, fix a spanning forest $T$ of $\Gamma$. By switching successively within each connected component, any signing can be transformed so that every edge of $T$ is positive. If $e \notin T$, then $e$ together with the unique path in $T$ joining its endpoints forms a fundamental cycle $C_e$. After the tree edges have been made positive, the sign of $e$ is precisely $\chi(C_e)$. Two signings with the same cycle signs therefore agree on every non-tree edge after the same normalization and hence are switching equivalent. ▫

Corollary 6.8 determines the absolute Gram matrix and therefore determines $\Gamma$, the magnitudes $|\rho_{ij}|$, and the continuous part of the reconstruction problem. It remains to recover the cycle signs.

### 7.2. Strict separation of triple switching classes

For three nonzero absolute correlations

$$a, b, c \in (0,1),$$

there are at most two switching classes. After independent reversals of the three directions, representatives may be written as

$$G_{\pm}(a,b,c) = \begin{pmatrix} 1 & a & b \\ a & 1 & \pm c \\ b & \pm c & 1 \end{pmatrix}, \Delta_{\pm} = 1 - a^2 - b^2 - c^2 \pm 2abc. \quad (42)$$

The signs + and - are the two possible values of the gauge-invariant triangle sign

$$\chi_{123} = sgn(\rho_{12}\rho_{13}\rho_{23}).$$

Whenever the corresponding Gram matrix is positive semidefinite, write

$$J_{\pm}(a,b,c) = I_3(G_{\pm}(a,b,c)).$$

**Theorem 7.2 (strict triple-sign separation).** Let $a, b, c \in (0,1)$. The negative switching class is feasible precisely when $\Delta_{-} \geq 0$, and its feasibility implies feasibility of the positive class. If $\Delta_{-} < 0$, every realizable triple with absolute correlations $a, b, c$ has positive triangle sign. If $\Delta_{-} \geq 0$, then

$$J_{+}(a,b,c) - J_{-}(a,b,c) = \frac{1}{4\pi} \int_0^a \int_0^b \int_{-c}^{c} \frac{dz\,dy\,dx}{\sqrt{1 - x^2 - y^2 - z^2 + 2xyz}}. \quad (43)$$

In particular,

$$J_{+}(a,b,c) > J_{-}(a,b,c)$$

and

$$J_{+}(a,b,c) - J_{-}(a,b,c) \geq \frac{abc}{2\pi}.$$

Thus the labelled data

$$\left(|\rho_{12}|, |\rho_{13}|, |\rho_{23}|, I_3(n_1, n_2, n_3)\right)$$

determine $sgn(\rho_{12}\rho_{13}\rho_{23})$ uniquely whenever all three correlations are nonzero.

*Proof.* Since $a, b, c < 1$, every $2 \times 2$ principal minor of the matrices in (42) is positive. Feasibility is therefore equivalent to nonnegativity of the determinant $\Delta_{\pm}$. Since

$$\Delta_{+} - \Delta_{-} = 4\,abc > 0\,,$$

feasibility of the negative class implies feasibility of the positive class. If $\Delta_{-} < 0$, the negative class is impossible.

It remains to compare the two feasible values. Let

$$\Sigma(x, y, z) = \begin{pmatrix} 1 & x & y \\ x & 1 & z \\ y & z & 1 \end{pmatrix}$$

and, on its positive-definite domain, let

$$\Phi(x, y, z) = E|X_1 X_2 X_3|, (X_1, X_2, X_3) \sim N(0, \Sigma(x, y, z)).$$

For $\delta > 0$, set

$$f_\delta(t) = \sqrt{t^2 + \delta^2}$$

and define

$$\Phi_\delta(x, y, z) = E\prod_{j=1}^{3} f_\delta(X_j), (X_1, X_2, X_3) \sim N(0, \Sigma(x, y, z)).$$

Price’s covariance-differentiation identity [20], applied successively to the three correlation parameters, gives

$$\frac{\partial^3 \Phi_\delta}{\partial x \partial y \partial z} = E\left[f_\delta''(X_1) f_\delta''(X_2) f_\delta''(X_3)\right],$$

where

$$f_\delta''(t) = \frac{\delta^2}{(t^2 + \delta^2)^{3/2}} \geq 0\,.$$

The measures $f_\delta''(t)\,dt$ converge weakly to $2\delta_0$. On every compact subset of the positive-definite correlation domain, the corresponding Gaussian densities and their first derivatives are uniformly bounded in a neighbourhood of the origin. Hence

$$E\left[f_\delta''(X_1) f_\delta''(X_2) f_\delta''(X_3)\right] \to 8\, p_{\Sigma(x, y, z)}(0)$$

uniformly on such compact subsets. Moreover $f_\delta(t) \leq |t| + \delta$, so dominated convergence gives $\Phi_\delta \to \Phi$. Integrating the preceding Price identity over a compact rectangular box in the positive-definite domain and passing to the limit therefore gives

$$\frac{\partial^3 \Phi}{\partial x \partial y \partial z} = 8\, p_{\Sigma(x,y,z)}(0) = \frac{8}{(2\pi)^{3/2}\sqrt{1-x^2-y^2-z^2+2xyz}}$$

in the integrated sense used below. For completeness, the factor $8$ can also be read directly from the distributional identity

$$\partial_{12}\partial_{13}\partial_{23}\left(|x_1||x_2||x_3|\right) = 8\,\delta_0(x_1)\,\delta_0(x_2)\,\delta_0(x_3).$$

Assume first that $\Delta_- > 0$. For

$$0 \le x \le a,\, 0 \le y \le b,\, |z| \le c,$$

one has

$$x^2 + y^2 + z^2 - 2xyz \le a^2 + b^2 + c^2 + 2abc < 1,$$

so the entire rectangular parameter region lies inside the positive-definite correlation domain. Moreover,

$$\Phi(0, y, c) = \Phi(0, y, -c),\, \Phi(x, 0, c) = \Phi(x, 0, -c),$$

because when one of the other correlations is zero the change $c \mapsto -c$ is realized by reversing one Gaussian coordinate, which leaves the absolute product unchanged. Integrating the mixed derivative over the parameter box therefore yields

$$\Phi(a,b,c) - \Phi(a,b,-c) = \int_0^a \int_0^b \int_{-c}^c \frac{8\, dz\, dy\, dx}{(2\pi)^{3/2}\sqrt{1-x^2-y^2-z^2+2xyz}}.$$

Lemma 6.9 with $k = 3$ gives

$$I_3 = \frac{\sqrt{\pi}}{8\sqrt{2}}\Phi,$$

and the constants reduce to $1/(4\pi)$, proving (43).

If $\Delta_- = 0$, then

$$\frac{\partial \Delta_-}{\partial c} = -2(c + ab) < 0.$$

Consequently, replacing $c$ by $(1-\varepsilon)c$ strictly increases $\Delta_-$ for every sufficiently small $\varepsilon > 0$, placing the perturbed triple in the positive-definite case. Applying the preceding identity and letting $\varepsilon \downarrow 0$, Proposition 6.10 gives continuity of $I_3$ at the rank-two boundary, while the right-hand side of (43) converges to the corresponding improper integral.

Finally,

$$x^2 + y^2 + z^2 - 2xyz = (x - yz)^2 + y^2(1 - z^2) + z^2 \ge 0.$$

Hence the determinant in the denominator of (43) is at most $1$, so its reciprocal square root is at least $1$. The integration box has volume $2abc$, and therefore

$$J_+ - J_- \geq \frac{2abc}{4\pi} = \frac{abc}{2\pi}.$$

The separation is strict because $a, b, c > 0$. ▫

Theorem 7.2 removes the genericity qualification from the nonzero triple-sign problem. There is no exceptional interior locus on which the two switching classes have the same third-order datum.

### 7.3. Reconstruction of finite normal configurations

Let $C(\Gamma)$ denote the cycle space of the correlation graph over $F_2$, and let

$$T(\Gamma) \subseteq C(\Gamma)$$

be the subspace generated by its triangles.

**Theorem 7.3 (finite normal tomography).** Let $n_1, \ldots, n_m \in S^2$ be labelled directions, no two of which are parallel or antiparallel, and let $\Gamma$ be their correlation graph. Assume that

$$T(\Gamma) = C(\Gamma).$$

Then the labelled pair data

$$\{I_2(n_i, n_j) : 1 \leq i < j \leq m\}$$

together with the labelled triple data

$$\{I_3(n_i, n_j, n_k) : 1 \leq i < j < k \leq m\}$$

determine the direction family up to transformations

$$n_i \mapsto \varepsilon_i Q n_i, \varepsilon_i \in \{\pm 1\}, Q \in O(3).$$

Equivalently, the pair and triple data determine the Gram matrix uniquely up to diagonal sign conjugacy.

If the correlation graph is complete, it is enough to use the triples containing one fixed reference direction. Thus, after choosing the reference index $1$, the pair data together with

$$\{I_3(n_1, n_i, n_j) : 2 \leq i < j \leq m\}$$

already determine the configuration. A complete nonorthogonal configuration therefore requires only

$$\binom{m}{2} + \binom{m-1}{2} = (m-1)^2$$

nontrivial labelled pair and triple values.

*Proof.* Corollary 6.8 recovers every magnitude

$$a_{ij} = |\rho_{ij}|$$

and hence the correlation graph $\Gamma$. For every triangle $ijk$ of $\Gamma$, all three correlations are nonzero, so Theorem 7.2 determines

$$\chi_{ijk} = sgn(\rho_{ij} \rho_{ik} \rho_{jk}).$$

Thus the pair and triple data recover the sign of every triangle cycle. By hypothesis the triangle cycles generate the full cycle space, so every cycle sign is determined. Proposition 7.1 then determines the signing of the Gram matrix up to diagonal sign conjugacy.

Choose any representative signed Gram matrix $G$ in that switching class. One representative is the Gram matrix of the original direction family and is therefore positive semidefinite of rank at most three. Diagonal sign conjugacy preserves both positive semidefiniteness and rank. Hence every representative in the recovered switching class satisfies

$$G \succcurlyeq 0 , rank\, G \leq 3 .$$

A Gram factorization produces vectors $\tilde{n}_1, \ldots, \tilde{n}_m \in \mathbb{R}^3$ with Gram matrix $G$. If two labelled vector families have the same Gram matrix, the correspondence between their spans sending one family to the other preserves all inner products. It is therefore an isometry of the spans and extends to an element of $O(3)$. Restoring the switching freedom gives precisely

$$\tilde{n}_i = \varepsilon_i Q n_i .$$

If $\Gamma$ is complete, fix the sign gauge by requiring

$$\rho_{1j} > 0 , j = 2 , \ldots , m .$$

Then, for $2 \leq i < j \leq m$,

$$sgn\, \rho_{ij} = sgn(\rho_{1i} \rho_{1j} \rho_{ij}) = \chi_{1ij} .$$

Hence the anchored triples determine every remaining edge sign. ▫

Theorem 7.3 is a labelled reconstruction theorem. It does not assert reconstruction from the scalar distribution of a width or brightness random variable. Its data retain their pair and triple labels, and that distinction is essential.

### *7.4. Orthogonality strata and the third-order information boundary*

The hypothesis $T(\Gamma) = C(\Gamma)$ is automatic for complete correlation graphs and, more generally, whenever every cycle sign is generated by triangle signs. It also identifies precisely what third-order data can fail to see when orthogonalities remove edges from the correlation graph.

**Proposition 7.4 (triangle-cycle completeness).** For an arbitrary labelled family of distinct unoriented directions, the complete pair and triple data determine the correlation graph, the magnitude $|\rho_{ij}|$ on every edge, and the cycle sign of every triangle of $\Gamma$. No additional switching invariant is contained in a triple whose induced correlation graph is a forest. Consequently, before the positive-semidefinite and rank constraints of Euclidean Gram matrices are imposed, the sign information not resolved by second- and third-order data is precisely the cycle information outside

$$T(\Gamma) \subseteq C(\Gamma) .$$

If $T(\Gamma) = C(\Gamma)$, the switching class is completely determined by pair and triple data. If the inclusion is strict, any remaining uniqueness must arise from additional Gram-feasibility constraints or from information of order higher than three.

*Proof.* The correlation graph and the absolute edge values follow from Corollary 6.8. The sign of every triangle follows from Theorem 7.2.

Consider a triple whose induced correlation graph is not a triangle. It is a forest. Every signing of a forest is switching equivalent to the all-positive signing, since signs may be removed successively along each tree component. Because $I_3$ is invariant under independent reversals of the directions, such a triple carries no sign invariant beyond its absolute correlations.

Proposition 7.1 identifies switching classes with cycle-sign data. Hence third-order measurements recover exactly the portion generated by triangle cycles. ▫

The triangle-generation hypothesis in Theorem 7.3 is sufficient for reconstruction. Its necessity after the positive-semidefinite rank-three constraints on Euclidean Gram matrices are imposed is not asserted here.

Repeated or antipodal directions, detected by $|\rho_{ij}|=1$, represent the same unoriented line and may be collapsed before applying the reconstruction theorem. Rank-two configurations require no separate uniqueness argument once the signed Gram matrix has been recovered: Gram factorization remains unique up to a global orthogonal transformation of the ambient space.

### *7.5. Stability on fixed correlation strata*

The exact reconstruction theorem also identifies the two mechanisms by which stability can deteriorate. Pair inversion becomes ill-conditioned near orthogonality, while sign recovery becomes ill-conditioned only when the nonzero triangle correlations approach zero.

**Proposition 7.5 (stability on fixed correlation strata).** Fix a correlation graph $\Gamma$ whose cycle space is generated by triangles, and let $\eta>0$. On the class of configurations satisfying

$$|\rho_{ij}|\geq\eta\left(ij\in E(\Gamma)\right),$$

pair reconstruction obeys

$$|a-\tilde{a}|\leq\frac{3\pi}{2\arcsin\eta}\left|I_2(a)-I_2(\tilde{a})\right|,(44)$$

for $a,\tilde{a}\in[\eta,1]$, while every feasible pair of triple switching classes satisfies

$$J_+(a,b,c)-J_-(a,b,c)\geq\frac{\eta^3}{2\pi}.(45)$$

*Proof.* Theorem 6.7 gives

$$I_2'(r)=\frac{2}{3\pi}\arcsin r,0\leq r\leq 1.$$

On $[\eta,1]$, this derivative is bounded below by $2\arcsin\eta/(3\pi)$, and the mean-value theorem gives (44). For a triangle with edge magnitudes $a,b,c\geq\eta$, Theorem 7.2 gives

$$J_+(a,b,c)-J_-(a,b,c)\geq\frac{abc}{2\pi}\geq\frac{\eta^3}{2\pi},$$

which proves (45). ▫

These estimates imply that, on a fixed correlation stratum separated from orthogonality, pair inversion is continuous and the two triple switching branches remain separated by a positive margin. The recovered triangle signs, and hence the switching class, are therefore locally constant under sufficiently small perturbations of the labelled data. Once the switching class is fixed, the signed Gram matrix depends continuously on those data.

If one further restricts to a fixed-rank positive-semidefinite stratum whose smallest positive eigenvalue stays bounded away from zero, the corresponding spectral projection and positive square-root factor vary continuously. Gram realizations consequently vary continuously modulo the relevant orthogonal action. Thus the recovered direction configuration is stable, in this restricted sense, modulo independent sign changes and the action of $O(3)$.

The present section concerns recovery of the direction configuration. Recovering unknown positive weights in an observable

$$F(u)=\sum_{i=1}^{m} c_i |u \cdot n_i|, c_i > 0,$$

is a distinct inverse problem and is not included in Theorem 7.3. No fourth-order absolute-product data are required for the reconstruction class covered here.

The resulting information chain is exact:

labelled second-order data → absolute Gram matrix → triangle cycle signs →

switching class → signed Gram matrix → direction configuration modulo sign and $O(3)$.

This is the inverse counterpart of the forward moment calculus of Section 6.

## 8. The regular truncated octahedron: exact width law and tomography

This is the first section in which the normal-fan density theorem, the polarized moment calculus, and finite normal tomography all close on the same object. The spherical-cell theory determines the complete scalar width law, the polarized moments expose its correlation orbits, and the labelled second- and third-order data reconstruct the six generating lines. The regular truncated octahedron is sufficiently non-orthogonal to test all three mechanisms while retaining enough symmetry for a concise exact calculation.

The regular permutohedron is the graphical zonotope of the complete graph; in dimension three this is the type-$A_3$ permutohedral realization. Postnikov gives the general graphical-zonotope description of the regular permutohedron [22, Proposition 2.3]. Independently, Lángi records the regular truncated octahedron as the six-generator type-(5) three-dimensional parallelohedron, with generators represented by the segments joining midpoints of opposite edges of a cube [23, §2.1]. Lángi further proves that, among unit-volume three-dimensional parallelohedra, the minimum mean width is attained precisely by regular truncated octahedra, up to congruence [23, Theorem 1]. Scaling gives the corresponding statement at every fixed positive volume. Thus the present section determines the complete width law of the fixed-volume mean-width minimizer in this class. We use a unit-generator normalization adapted to the calculation.

### 8.1. Zonotopal realization and width normal form

Set

$$r_1 = \frac{(1,1,0)}{\sqrt{2}}, \quad r_2 = \frac{(1,-1,0)}{\sqrt{2}},$$
$$r_3 = \frac{(1,0,1)}{\sqrt{2}}, \quad r_4 = \frac{(1,0,-1)}{\sqrt{2}},$$
$$r_5 = \frac{(0,1,1)}{\sqrt{2}}, \quad r_6 = \frac{(0,1,-1)}{\sqrt{2}},$$

and define

$$O = \sum_{j=1}^{6} \left[ -\frac{r_j}{2}, \frac{r_j}{2} \right]. \quad (46)$$

**Proposition 8.1 (truncated-octahedral normal form).** The zonotope $O$ is a unit-edge regular truncated octahedron, equivalently a unit-generator realization of the three-dimensional regular permutohedron. Its vertices are the twenty-four signed permutations of

$$\left( \sqrt{2}, \frac{1}{\sqrt{2}}, 0 \right).$$

For $u = (x, y, z) \in S^2$,

$$w_O(u) = \sqrt{2}\left[ max(|x|,|y|) + max(|x|,|z|) + max(|y|,|z|) \right].$$

Equivalently, if $a \geq b \geq c \geq 0$ is the decreasing rearrangement of $(|x|,|y|,|z|)$, then

$$w_O(u) = \sqrt{2}(2a + b). \quad (47)$$

The support of $W_O = w_O(U)$ is

$$supp(W_O) = [\sqrt{6}, \sqrt{10}].$$

*Proof.* Proposition 2.2 applied to (46) gives $w_O(u) = \sum_{j=1}^{6} |u \cdot r_j|$. For arbitrary real $p, q$,

$$|p+q| + |p-q| = 2\,max(|p|,|q|).$$

Applying this identity to the three coordinate pairs $(x, y)$, $(x, z)$, and $(y, z)$ gives (47).

On the chamber $x \geq y \geq |z|$, all six scalar products $u \cdot r_j$ are nonnegative. The corresponding support vertex is

$$\frac{1}{2} \sum_{j=1}^{6} r_j = \left( \sqrt{2}, \frac{1}{\sqrt{2}}, 0 \right).$$

The signed-permutation symmetry of the six generating directions produces its twenty-four signed permutations, each with a two-dimensional spherical normal cell. These are therefore precisely the vertices of $O$.

For example, the supporting plane $x=\sqrt{2}$ contains exactly

$$\left(\sqrt{2},\pm 1/\sqrt{2},0\right),\left(\sqrt{2},0,\pm 1/\sqrt{2}\right),$$

in cyclic order; successive differences have norm 1 and adjacent differences are orthogonal. It therefore cuts out a unit square. The supporting plane $x+y+z=3/\sqrt{2}$ contains the six permutations of $\left(\sqrt{2},1/\sqrt{2},0\right)$; consecutive permutations differ by a vector of norm 1, and the six successive edge angles are equal. It therefore cuts out a unit regular hexagon. Signed permutations now show that the six supporting planes $x_i=\pm\sqrt{2}$ cut out square facets, while the eight supporting planes

$$\pm x\pm y\pm z=\frac{3}{\sqrt{2}}$$

cut out regular hexagonal facets. Thus $O$ has six unit squares and eight unit regular hexagons and is the unit-edge regular truncated octahedron.

It remains to determine the support of its width. Since $a\ge b\ge c$,

$$1=a^2+b^2+c^2\le a^2+2b^2\le\frac{\left(2a+b\right)^2}{3},$$

where the second inequality is equivalent to

$$\left(a-b\right)\left(a+5b\right)\ge 0.$$

Hence $2a+b\ge\sqrt{3}$, with equality only when $a=b=c=1/\sqrt{3}$. This gives the lower endpoint $\sqrt{6}$.

At the other end,

$$2a+b\le\sqrt{5}\sqrt{a^2+b^2}\le\sqrt{5},$$

with equality precisely when $c=0$ and $a:b=2:1$. Hence the upper endpoint is $\sqrt{10}$. ▫

The lower endpoint is attained at the eight cube-diagonal directions. The upper endpoint is attained at the twenty-four normal-cell axes, the signed permutations of $\left(2,1,0\right)/\sqrt{5}$.

### 8.2. A single normal cell and its aperture

The full density follows from one normal cell. Let

$$C=\{\left(x,y,z\right)\in S^2 : x\ge y\ge|z|\}.$$

On $C$, equation (47) becomes

$$w_O\left(u\right)=q_O\cdot u, q_O=\left(2\sqrt{2},\sqrt{2},0\right), D_O:=\|q_O\|=\sqrt{10}.$$

The coefficient $q_O$ is the vertex $2\left(\sqrt{2},1/\sqrt{2},0\right)$ of the difference body $O-O=2O$ whose normal cell is $C$. The three walls of $C$ have inward unit normals

$$b_0=\frac{(1,-1,0)}{\sqrt{2}},b_+=\frac{(0,1,-1)}{\sqrt{2}},b_-=\frac{(0,1,1)}{\sqrt{2}}.$$

**Lemma 8.2 (truncated-octahedral aperture).** Put

$$\vartheta=\arccos\left(-\frac{2}{3}\right),\beta(t)=\arccos\frac{t}{3\sqrt{10-t^2}}\left(\sqrt{6}\le t\le 3\right).\ (48)$$

The aperture of the level circle $C\cap\{u:q_O\cdot u=t\}$ is

$$\Theta_O(t)=\begin{cases}0, & t<\sqrt{6},\\ 2\vartheta-4\beta(t), & \sqrt{6}\le t<2\sqrt{2},\\ 2\pi-6\beta(t), & 2\sqrt{2}\le t<3,\\ 2\pi, & 3\le t<\sqrt{10},\\ 0, & t\ge\sqrt{10}.\end{cases}$$

*Proof.* Write $\hat{q}_O=q_O/\sqrt{10}$. A point of the level circle has the form

$$u=\frac{t}{\sqrt{10}}\hat{q}_O+\frac{\sqrt{10-t^2}}{\sqrt{10}}v,v\in q_O^{\perp},\|v\|=1.$$

Each wall normal satisfies $b_j\cdot\hat{q}_O=1/\sqrt{10}$. Projecting $b_j$ to $q_O^{\perp}$, define

$$p_j=\frac{\sqrt{10}}{3}\left(b_j-\frac{1}{\sqrt{10}}\hat{q}_O\right).$$

Then $\|p_j\|=1$, and the wall inequality $b_j\cdot u\ge 0$ becomes

$$p_j\cdot v\ge-\frac{t}{3\sqrt{10-t^2}}.$$

For $t\le 3$, each wall therefore removes from the level circle an interval centered at $-p_j$ with half-width $\beta(t)$.

A direct computation gives

$$p_0\cdot p_+=p_0\cdot p_-=-\frac{2}{3},p_+\cdot p_-=-\frac{1}{9}.$$

Consequently the three cyclic gaps between the forbidden-interval centres are

$$\vartheta,\vartheta,2\pi-2\vartheta,$$

because $\cos(2\vartheta)=-1/9$.

At the lower support, $\beta(\sqrt{6})=\vartheta/2$, so the three forbidden arcs cover the circle and the admissible aperture is zero. As $t$ increases, $\beta(t)$ decreases. In the first active regime the only nonempty pairwise overlap is the overlap across the cyclic gap $2\pi-2\vartheta$; the two overlaps across the gaps $\vartheta$ are empty. This overlap disappears when $2\beta = 2\pi - 2\vartheta$, equivalently when $\beta = \pi - \vartheta$. Since $\cos(\pi - \vartheta) = 2/3$, the transition occurs at $t = 2\sqrt{2}$. The union of the three forbidden intervals then has length

$$6\beta - (2\beta - (2\pi - 2\vartheta)) = 4\beta + 2\pi - 2\vartheta,$$

and the complementary aperture is $2\vartheta - 4\beta$.

For $2\sqrt{2} \le t < 3$, the three forbidden intervals are disjoint, so the complementary aperture is $2\pi - 6\beta$. At $t = 3$, $\beta = 0$; for $t > 3$ the right-hand side of each wall inequality is less than $-1$, and the entire level circle lies in the cell. The axial value $D_O = \sqrt{10}$ terminates the cell. ▫

The disappearance of the single overlap at $2\sqrt{2}$ and the subsequent disjoint-arc regime are shown in Figures 1(b) and 1(c).

The four critical values are therefore

$$\sqrt{6}, 2\sqrt{2}, 3, \sqrt{10}.$$

### *8.3. Exact width density*

The twenty-four normal cells of $O$ are the signed-permutation images of $C$, and every cell has coefficient norm $\sqrt{10}$. The global co-area theorem therefore reduces the complete law to Lemma 8.2.

**Theorem 8.3 (exact width law of the unit-edge regular truncated octahedron).** Let $U$ be uniform on $S^2$ and let $W_O = w_O(U)$. Then $W_O$ is supported on $[\sqrt{6}, \sqrt{10}]$ and has density

$$f_O(t) = \begin{cases} 0, & t < \sqrt{6}, \\ \dfrac{12}{\pi\sqrt{10}}(\vartheta - 2\beta(t)), & \sqrt{6} \le t < 2\sqrt{2}, \\ \dfrac{12}{\pi\sqrt{10}}(\pi - 3\beta(t)), & 2\sqrt{2} \le t < 3, \\ \dfrac{12}{\sqrt{10}}, & 3 \le t < \sqrt{10}, \\ 0, & t \ge \sqrt{10}, \end{cases} \quad (49)$$

where $\vartheta$ and $\beta$ are defined in (48).

*Proof.* The vertex $q_O \in O - O$ has spherical normal cell $C$, and its signed-permutation orbit consists of the twenty-four vertices of $O - O = 2O$. Theorem 3.4 therefore gives

$$f_O(t) = \frac{24}{4\pi} \frac{\Theta_O(t)}{\sqrt{10}} = \frac{6}{\pi\sqrt{10}} \Theta_O(t).$$

Substitution of the four branches from Lemma 8.2 yields (49). The spherical normal cells cover $S^2$ up to their one-dimensional boundaries, so Theorem 3.1 also gives normalization of the resulting probability density. ▫

An elementary antiderivative useful for an independent normalization check is

$$J(t) = t\beta(t) - \sqrt{10}\arctan\sqrt{9 - t^2}, J'(t) = \beta(t).$$

The endpoint values are

$$\begin{aligned} J(\sqrt{6}) &= \frac{\sqrt{6}}{2}\vartheta - \frac{\pi\sqrt{10}}{3}, \\ J(2\sqrt{2}) &= 2\sqrt{2}(\pi - \vartheta) - \frac{\pi\sqrt{10}}{4}, \\ J(3) &= 0. \end{aligned}$$

Direct integration of the three nonzero branches of (49) gives

$$\int_{\sqrt{6}}^{\sqrt{10}} f_O(t)\,dt = 1.$$

The same primitive gives the CDF in elementary form if required. In particular, the terminal branch is

$$F_O(t) = \frac{12t}{\sqrt{10}} - 11, 3 \le t \le \sqrt{10}.$$

To the best of the author's knowledge, no explicit closed-form random-direction width density for the regular truncated octahedron has previously been recorded. This scoped statement is based on the width, caliper, Feret-diameter, and permutohedral searches documented in the prior-art ledger; it is not an unrestricted claim of historical priority.

### *8.4. Critical geometry and transition structure*

The derivative

$$\beta'(t) = -\frac{\sqrt{10}}{(10 - t^2)\sqrt{9 - t^2}} \quad (50)$$

makes all four global transitions explicit.

**Proposition 8.4 (transition and resonance structure).** At the lower support,

$$f_O(t) = \frac{2\sqrt{3}}{\pi}(t - \sqrt{6}) + O\left((t - \sqrt{6})^2\right).$$

At $t = 2\sqrt{2}$, the density is continuous, but its derivative changes from

$$f_O'(2\sqrt{2}^-) = \frac{12}{\pi}$$

to

$$f_O'(2\sqrt{2}^+) = \frac{18}{\pi}.$$

At $t = 3$,

$$f_O(t) = \frac{12}{\sqrt{10}} - \frac{12\sqrt{6}}{\pi}\sqrt{3-t} + O\left((3-t)^{3/2}\right), t \uparrow 3.$$

Finally, $f_O$ has a terminal jump from $12/\sqrt{10}$ to $0$ at $t = \sqrt{10}$.

The spherical normal fan has twenty-four two-cells, thirty-six edges, and fourteen vertices. The critical value $\sqrt{6}$ is supported by the eight cube-diagonal fan vertices, with six normal cells meeting at each one. The value $2\sqrt{2}$ is supported by the six coordinate-axis fan vertices, with four cells meeting at each one. The value $3$ is the common tangency value of the thirty-six spherical fan edges, while $\sqrt{10}$ is the common axial value of the twenty-four two-cells. Thus the four transition types of Question 3.5 occur on four different fan strata or incidence events in a single highly resonant normal fan.

*Proof.* Equation (50) follows by differentiating (48). Substitution into the first two branches of (49) gives the lower linear expansion and the two one-sided derivatives at $2\sqrt{2}$.

As $t \uparrow 3$,

$$\beta(t) = \sqrt{\frac{20}{3}}\sqrt{3-t} + O\left((3-t)^{3/2}\right),$$

and substitution into the third branch of (49) gives the stated square-root expansion.

For the incidence counts, the fourteen fan vertices are the normals to the eight hexagonal and six square facets of $O$. A regular hexagon has six vertices and a square has four, giving respectively $8 \cdot 6 = 48$ and $6 \cdot 4 = 24$ cell-vertex incidences. The thirty-six fan edges correspond to the thirty-six edges of $O$, and the twenty-four fan cells correspond to its twenty-four vertices. The Euler relation

$$14 - 36 + 24 = 2$$

provides the consistency check. ▫

Thus the same exact density contains a linear onset, a finite derivative corner, a square-root fold, and a terminal step, with the resonance multiplicities determined by the octahedral symmetry.

Figure 3 compares the three complete densities obtained in Theorems 4.3, 5.4, and 8.3. The tetrahedral and truncated-octahedral width laws exhibit the same ordered sequence of a linear onset, a finite derivative corner, a square-root fold, and a terminal step, although the first is governed by a non-zonotopal difference body and the second by a zonotope. Tetrahedral brightness has a different sequence: its first interior transition is a square-root fold, followed by the interior jump created by axial termination of the $\Sigma$-cell family. The comparison therefore exhibits within three exact laws all of the local mechanisms appearing in Question 3.5. These conclusions follow from Lemmas 4.4 and 5.5 and Proposition 8.4; the plots serve only to display the proved formulas.

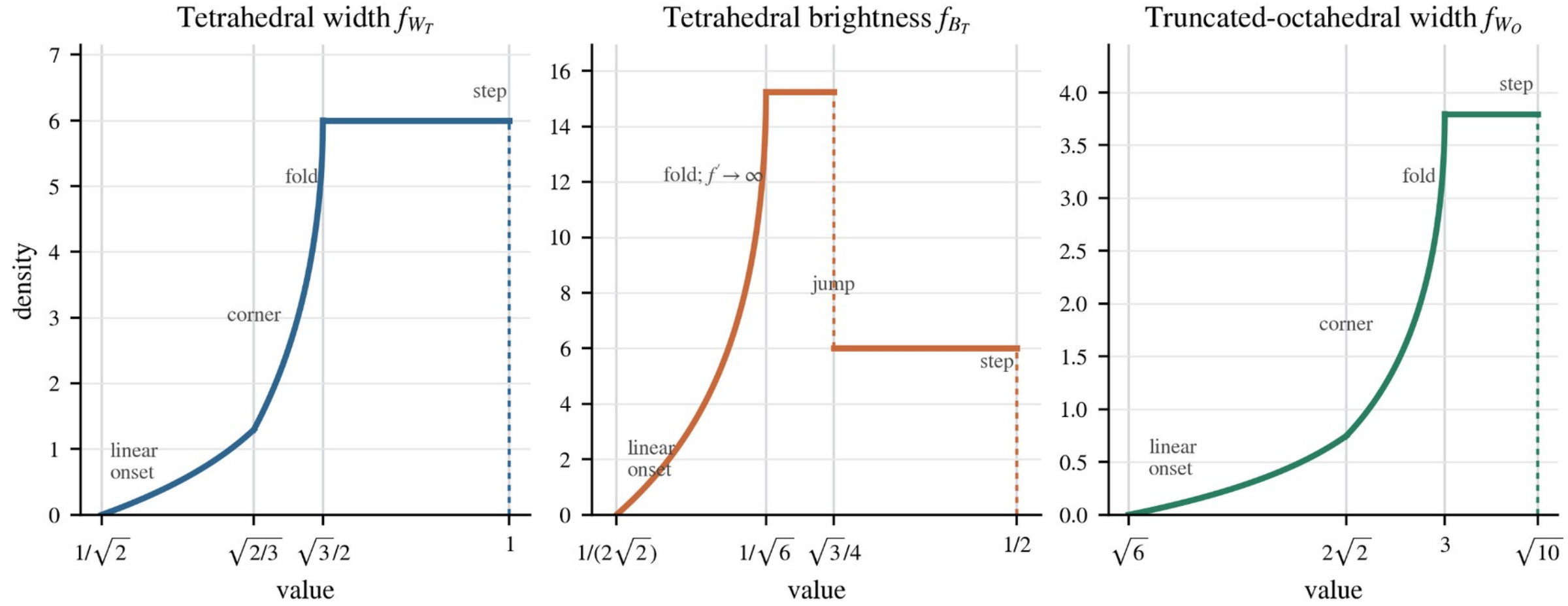

**Figure 3. Exact flagship densities and their transition types.** The densities of the unit-edge regular tetrahedron width $W_T$, the unit-edge regular tetrahedron brightness $B_T$, and the unit-edge regular truncated-octahedron width $W_O$, as given by (8), (18), and (49). Vertical markers indicate the complete critical sets. For $W_T$, the transitions are a linear onset at $1/\sqrt{2}$, a derivative corner at $\sqrt{2/3}$, a square-root fold at $\sqrt{3}/2$, and a terminal step at $1$. For $B_T$, they are a linear onset at $1/(2\sqrt{2})$, a square-root fold with $f'_{B_T}(t) \to +\infty$ at $1/\sqrt{6}$, an interior jump at $\sqrt{3}/4$, and a terminal step at $1/2$. For $W_O$, they are a linear onset at $\sqrt{6}$, a derivative corner at $2\sqrt{2}$, a square-root fold at $3$, and a terminal step at $\sqrt{10}$. Dashed vertical segments mark jumps of the chosen density representatives and do not represent atoms. By Theorem 5.9, the middle panel also determines the unit-edge rhombic-dodecahedral width density through

$$f_R(s) = \frac{\sqrt{3}}{8} f_{B_T}\left(\frac{\sqrt{3}s}{8}\right).$$

## 8.5. Exact moments, cumulants, and harmonic concentration

The moment calculation has a much smaller orbit structure than a direct expansion of the density suggests. The six generator lines form three orthogonal pairs,

$$(r_1, r_2), (r_3, r_4), (r_5, r_6),$$

and every other pair has absolute correlation $1/2$.

**Proposition 8.5 (first three moments and cumulants).** For the unit-edge regular truncated octahedron,

$$\begin{aligned} E W_O &= 3, \\ E W_O^2 &= \frac{10}{3} + \frac{4 + 8\sqrt{3}}{\pi}, \\ E W_O^3 &= \frac{87}{4} + \frac{12\sqrt{2}}{\pi}. \end{aligned} \tag{51}$$

Consequently its first three cumulants are

$$\begin{aligned}\kappa_1(W_O) &= 3,\\ \kappa_2(W_O) &= -\frac{17}{3}+\frac{4+8\sqrt{3}}{\pi},\\ \kappa_3(W_O) &= \frac{183}{4}+\frac{12\sqrt{2}-36-72\sqrt{3}}{\pi}.\end{aligned} \quad (52)$$

*Proof.* Since $W_O = \sum_{j=1}^{6} |U \cdot r_j|$,

$$E\,W_O = 6 \cdot \frac{1}{2} = 3.$$

There are three unordered orthogonal pairs and twelve unordered pairs with absolute correlation $1/2$. Theorem 6.7 gives

$$I_2(0) = \frac{2}{3\pi}$$

and

$$I_2\left(\frac{1}{2}\right) = \frac{\sqrt{3}}{3\pi} + \frac{1}{18}.$$

Expanding $W_O^2$ gives the second identity in (51).

For the third moment, repeated directions contribute $I_3(n,n,n) = 1/4$, and Proposition 6.10 gives

$$I_3(n,n,m) = \frac{1+\rho^2}{8}.$$

Thus the repeated-pair values are $1/8$ for an orthogonal pair and $5/32$ when $|\rho| = 1/2$.

Among the twenty triples of distinct generator lines, twelve contain one of the three orthogonal pairs. Their correlation type is $(1/2, 0, 1/2)$, for which Proposition 6.10 gives

$$I_{path} = \frac{2\sqrt{2} + 4\arcsin(1/\sqrt{3}) - \arcsin(1/3)}{16\pi}.$$

The remaining eight triples contain one generator from each orthogonal pair. Their absolute correlations are all $1/2$. Four have positive triangle sign and four have negative triangle sign. The corresponding values are

$$I_+ = \frac{2\sqrt{2} + 9\arcsin(1/3)}{16\pi}$$

and

$$I_- = \frac{3}{32}.$$

The second value lies on the rank-two boundary $\Delta_- = 0$.

Using

$$\arcsin\frac{1}{3}+2\arcsin\frac{1}{\sqrt{3}}=\frac{\pi}{2},$$

the multinomial expansion of $W_O^3$ simplifies to

$$E W_O^3=\frac{87}{4}+\frac{12\sqrt{2}}{\pi}.$$

The inverse-sine identity follows by writing $\alpha=\arcsin\left(1/\sqrt{3}\right)$, for which $\cos\left(2\alpha\right)=1/3$ and $0<2\alpha<\pi/2$.

Finally,

$$\kappa_2=E W_O^2-\left(E W_O\right)^2$$

and

$$\kappa_3=E W_O^3-3 E W_O E W_O^2+2\left(E W_O\right)^3,$$

which gives (52). ▫

The mean $3$ agrees with the classical six-generator mean-width formula in the present normalization [23, Remark 1]. The higher moments provide independent checks on the complete density (49). The singular extension of Proposition 6.10 is first needed here at the rank-two value $I_{-}=3/32$; the other orbit types have nonsingular Gram matrices.

The third cumulant in (52) is the first translation-invariant odd cumulant of $W_O$. It records third-order asymmetry information that is absent from the mean and variance and places the exact law within the cumulant framework developed in [2].

The same six directions also explain where the variance is spectrally concentrated. Their frame operator satisfies

$$\sum_{j=1}^{6} r_j r_j^T=2I.$$

Appendix C shows from first principles that this tight-frame identity is equivalent to the vanishing of the degree-two term in Theorem 6.12. Hence the variance of $W_O$ begins at spherical-harmonic degree four. The tetrahedral facet-normal system has the same property, with frame operator $4I/3$; Appendix C gives the exact degree-four and degree-six contributions for both systems.

### *8.6. Tomography of the six generating directions*

The generator system provides a noncomplete correlation graph on which Theorem 7.3 nevertheless closes exactly. Let $\Gamma_O$ be the correlation graph of the six unoriented generator lines. Its three missing edges are the three orthogonal pairs, so

$$\Gamma_O \cong K_6 \setminus 3K_2 \cong K_{2,2,2}.$$

**Theorem 8.6 (finite tomography of the truncated-octahedral generators).** The labelled second- and third-order data of the six generating directions of $O$ determine the complete direction family up to independent sign changes and a global orthogonal transformation. More precisely, the pair data identify the three orthogonal pairs and determine that every remaining absolute correlation is $1/2$. The eight triangles of $\Gamma_O$ span its full seven-dimensional cycle space, and their switching classes identify exactly the four coplanar triples that generate the eight hexagonal facets of $O$; the three orthogonal pairs generate the six square facets.

*Proof.* The graph $\Gamma_O$ has six vertices and twelve edges and is connected. Hence

$$dim\, C\left(\Gamma_O\right)=12-6+1=7.$$

Write the three parts of $K_{2,2,2}$ as $A, B, C$, each with two vertices. There are eight triangles, one for each choice of one vertex from each part.

Suppose an $F_2$-linear combination of these eight triangles vanishes. For each $AB$-edge, the coefficients of the two triangles containing that edge must agree. Repeating the same argument with the $AC$- and $BC$-edges forces all eight triangle coefficients to be equal. Thus the only nontrivial relation is the sum of all eight triangles, since every edge occurs in exactly two of them. Their span therefore has dimension seven, and

$$T\left(\Gamma_O\right)=C\left(\Gamma_O\right).$$

Corollary 6.8 distinguishes the two pair values $I_2(0)$ and $I_2(1/2)$, and hence identifies the three missing edges of $\Gamma_O$.

For the triangle signs, write

$$A_\varepsilon=\frac{(1,\varepsilon,0)}{\sqrt{2}}, B_\delta=\frac{(1,0,\delta)}{\sqrt{2}}, C_\gamma=\frac{(0,1,\gamma)}{\sqrt{2}},$$

where $\varepsilon, \delta, \gamma \in \{\pm 1\}$. Then

$$A_\varepsilon \cdot B_\delta=\frac{1}{2}, A_\varepsilon \cdot C_\gamma=\frac{\varepsilon}{2}, B_\delta \cdot C_\gamma=\frac{\delta\gamma}{2},$$

so the triangle sign is $\varepsilon\delta\gamma$. Four of the eight triangles therefore have positive switching class and four have negative switching class. Theorem 7.2 separates these classes strictly. Theorem 7.3 now reconstructs the full generator family modulo sign gauge and $O(3)$.

The facet incidence is read from the same signs. Direct expansion gives

$$det\left(A_\varepsilon, B_\delta, C_\gamma\right)=-\frac{\delta+\varepsilon\gamma}{2\sqrt{2}}.$$

Thus the determinant vanishes exactly when $\varepsilon\delta\gamma=-1$, precisely the negative switching class. These four coplanar triples are the rank-two cases producing $I_-=3/32$; their two choices of supporting sign generate the eight regular hexagonal facets. The three missing edges of $\Gamma_O$ are exactly the orthogonal generator pairs, and their two supporting signs generate the six square facets.

The coincidence between negative switching class and coplanarity uses the common nonzero absolute correlation $1/2$, for which

$$\Delta_{-}\left(\frac{1}{2},\frac{1}{2},\frac{1}{2}\right)=1-\frac{3}{4}-\frac{1}{4}=0.$$

It is therefore a special property of this configuration rather than a general consequence of switching sign. For the truncated-octahedral system, however, the moment-theoretic switching classes, the rank boundary of the three-direction kernel, and the facet structure are the same combinatorial datum. ▫

This example is materially different from the complete-correlation case following Theorem 7.3: three pair correlations vanish, but the remaining triangles still generate the entire cycle space. The orthogonality pattern therefore does not obstruct third-order reconstruction.

### *8.7. Brightness comparison*

The face structure of the same solid gives a second absolute-normal representation. Let $e_1, e_2, e_3$ be the three coordinate-axis directions, corresponding to the three pairs of square facets. Let $m_1, \ldots, m_4$ be one normal from each opposite pair of hexagonal facets, chosen as the tetrahedral directions

$$\frac{(1,1,1)}{\sqrt{3}},\frac{(1,-1,-1)}{\sqrt{3}},\frac{(-1,1,-1)}{\sqrt{3}},\frac{(-1,-1,1)}{\sqrt{3}}.$$

Each square facet has area $1$, while each unit-edge regular hexagonal facet has area $3\sqrt{3}/2$. Proposition 2.4 therefore gives

$$b_O(u)=\sum_{r=1}^{3}\left|u\cdot e_r\right|+\frac{3\sqrt{3}}{2}\sum_{i=1}^{4}\left|u\cdot m_i\right|. \quad (53)$$

The two normal systems satisfy $e_r\cdot e_s=0$ for $r\neq s$, $m_i\cdot m_j=-1/3$ for $i\neq j$, and $\left|e_r\cdot m_i\right|=1/\sqrt{3}$.

**Proposition 8.7 (first two brightness moments).** For the unit-edge regular truncated octahedron,

$$E\,B_O=\frac{3}{2}+3\sqrt{3}=\frac{S(O)}{4}$$

and

$$E\,B_O^2=16+\frac{4+60\sqrt{2}+6\arcsin(1/3)}{\pi}. \quad (54)$$

*Proof.* The first identity follows directly from (53) and $E\left|U\cdot n\right|=1/2$, equivalently from Cauchy's formula $E\,B_O=S(O)/4$.

For the second moment, expand the square of (53). There are three orthogonal $e$-pairs, six tetrahedral $m$-pairs with absolute correlation $1/3$, and twelve mixed pairs with absolute correlation $1/\sqrt{3}$. Theorem 6.7 gives their exact $I_2$-values. Substitution and the identity

$$2\arcsin\frac{1}{\sqrt{3}}+\arcsin\frac{1}{3}=\frac{\pi}{2}$$

reduce the result to (54). ▫

Thus the same regular truncated octahedron supplies two finite normal systems of different types: its width is generated by six equal root directions with correlation graph $K_{2,2,2}$, whereas its brightness is generated by an orthogonal triple together with a weighted tetrahedral quadruple. The former is completely reconstructed by the pair-triple tomography of Section 7; the latter displays directly how the surface-area measure introduces a different weighted normal geometry. No complete brightness density is claimed here: Proposition 8.7 records exactly the first two brightness moments needed for the comparison.

The regular truncated octahedron therefore closes the three forward and inverse strands developed earlier: Theorem 3.4 produces its exact density from the difference-body normal fan, Section 6 reduces its moments and harmonic variance to a small correlation-orbit calculation, and Theorem 8.6 reconstructs the six generating lines and recovers the facet incidence from the same pair-triple data. This is the second non-orthogonal width flagship required by the general theory.

# 9. Intrinsic volumes and the loss of directional information

The intrinsic volumes retain fundamental scalar geometry, but they do not generally retain a complete random-direction law. In the normalization used here, the Steiner formulas in dimensions three and two are

$$Vol_3(K+\rho B^3)=V_3(K)+2\rho V_2(K)+\pi\rho^2 V_1(K)+\frac{4\pi}{3}\rho^3 V_0(K),$$

and

$$Area(P+\rho B^2)=V_2(P)+2\rho V_1(P)+\pi\rho^2 V_0(P). \quad (55)$$

Thus $V_3$, $V_2$, and $V_1$ in three dimensions encode volume, half the surface area, and twice the mean width, respectively. For rectangular boxes these three quantities determine the edge multiset and hence the complete width law [1, Proposition 6.7]. The purpose of this section is to show, by an exact centrally symmetric construction, that this closure is special to the box class: even the width and brightness second moments are not determined by the intrinsic volumes in the full class of centrally symmetric convex bodies.

## 9.1. Intrinsic-volume normalization and planar separation

The underlying mechanism is already visible in the plane. Let $m<n$ be distinct even integers, and choose

$$0<a<\frac{1}{\sqrt{(m^2-1)(n^2-1)}}, b=a\sqrt{\frac{m^2-1}{n^2-1}}.$$

Define the $2\pi$-periodic functions

$$h_m(\theta)=1+a\cos m\theta, h_n(\theta)=1+b\cos n\theta. \quad (56)$$

**Proposition 9.1 (planar harmonic separation at fixed intrinsic volumes).** The functions in (56) are the support functions of smooth strictly convex centrally symmetric plane bodies $P_m$ and $P_n$. These bodies have the same planar intrinsic volumes $V_0, V_1, V_2$, but

$$E\,w_{P_m}(U_\circ)^2 - E\,w_{P_n}(U_\circ)^2 = \frac{2a^2(n^2-m^2)}{n^2-1} > 0\,, (57)$$

where $U_\circ$ is uniform on $S^1$.

*Proof.* Put $q_k(\theta) = h_k(\theta) + h_k''(\theta)$. The support parametrization

$$x_k(\theta) = h_k(\theta)(\cos\theta, \sin\theta) + h_k'(\theta)(-\sin\theta, \cos\theta)$$

satisfies

$$x_k'(\theta) = q_k(\theta)(-\sin\theta, \cos\theta).$$

For $k = m$,

$$q_m(\theta) = 1 - (m^2-1)\,a\cos m\theta > 0\,,$$

because the bound on $a$ is stronger than $a < 1/(m^2-1)$. For $k = n$,

$$q_n(\theta) = 1 - (n^2-1)\,b\cos n\theta = 1 - a\sqrt{(m^2-1)(n^2-1)}\cos n\theta > 0\,.$$

Hence both parametrizations are regular with everywhere positive curvature radius and define smooth strictly convex bodies. Since $m$ and $n$ are even, $h_k(\theta+\pi) = h_k(\theta)$, so both bodies are centrally symmetric.

Their perimeters are

$$L(P_k) = \int_0^{2\pi} q_k(\theta)\,d\theta = \int_0^{2\pi} h_k(\theta)\,d\theta = 2\pi\,.$$

Integration of $det(x_k, x_k') = h_k q_k$ gives

$$Area(P_k) = \frac{1}{2}\int_0^{2\pi} h_k(\theta)\,q_k(\theta)\,d\theta = \pi\left(1 - \frac{k^2-1}{2}c_k^2\right),$$

where $c_m = a$ and $c_n = b$. The definition of $b$ makes these two areas equal. Equation (55) now gives equality of $V_0, V_1, V_2$.

Central symmetry gives $w_{P_k}(\theta) = 2h_k(\theta)$. Orthogonality of the trigonometric modes therefore yields

$$E\,w_{P_k}(U_\circ)^2 = \frac{1}{2\pi}\int_0^{2\pi} 4h_k(\theta)^2\,d\theta = 4 + 2c_k^2.$$

Consequently

$$E\,w_{P_m}(U_\circ)^2 - E\,w_{P_n}(U_\circ)^2 = 2(a^2-b^2) = \frac{2a^2(n^2-m^2)}{n^2-1}\,,$$

which is positive because $m<n$. ▫

The proposition isolates the information loss. Perimeter retains the constant Fourier mode of the support function, while area retains one weighted quadratic sum of its nonconstant modes. These data do not determine the unweighted quadratic sum governing the width second moment. This complements the planar noninjectivity mechanisms of [2] with an intrinsic-volume-preserving family.

A concrete member of the family in Proposition 9.1 is shown in Figure 4. We take

$$m=2, n=4, a=\frac{13}{100}, b=\frac{13}{100\sqrt{5}}.$$

The strict-convexity condition holds because $13/100<1/\sqrt{45}$. This specialization is used only to display the planar construction; Proposition 9.1 remains stated for every pair of distinct even integers $m<n$ satisfying its hypotheses.

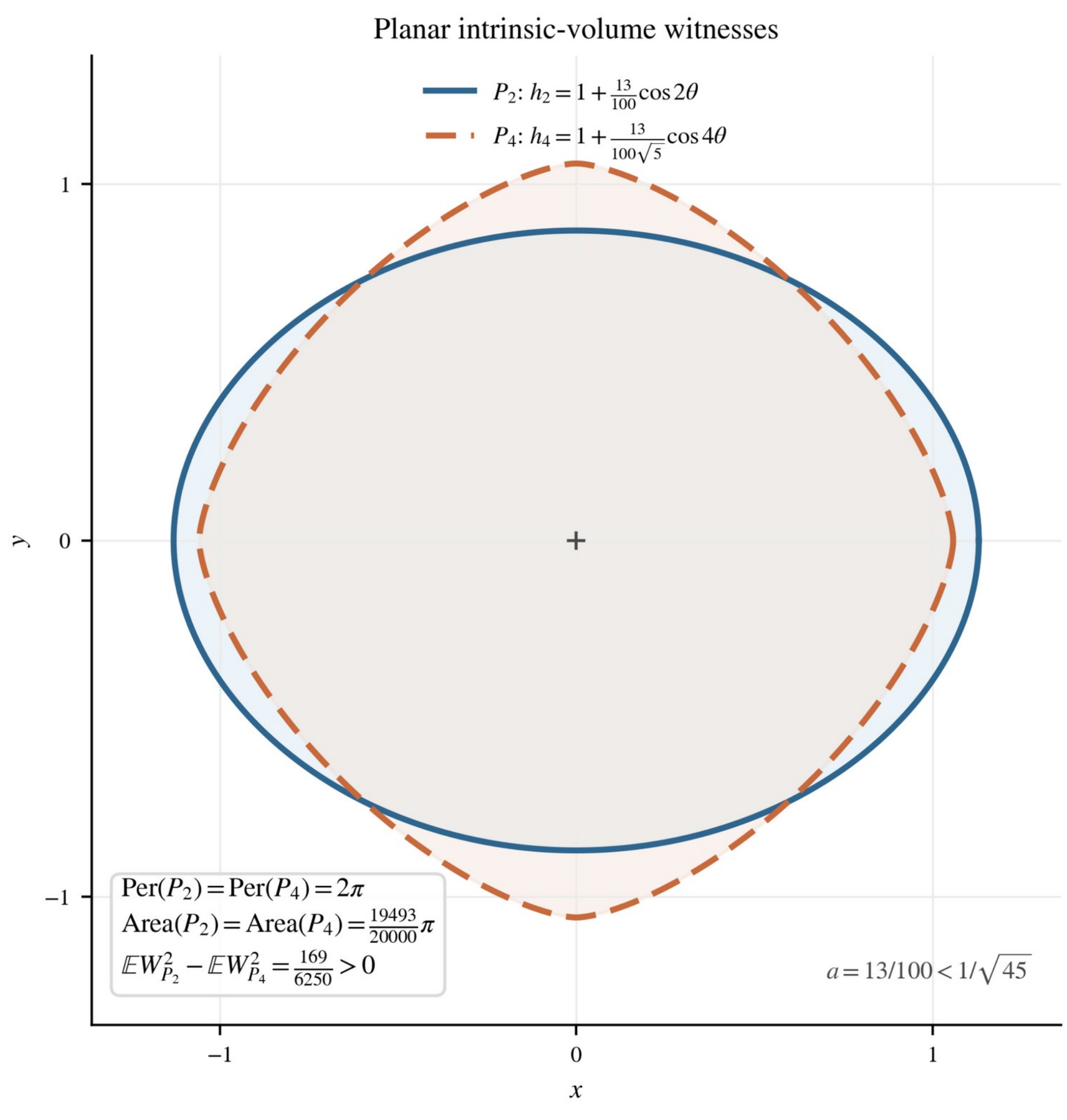

**Figure 4. Planar bodies with equal intrinsic volumes and different width laws.** The exact support-function boundaries $P_2$ and $P_4$ obtained from

$$h_2(\theta)=1+\frac{13}{100}\cos 2\theta, h_4(\theta)=1+\frac{13}{100\sqrt{5}}\cos 4\theta,$$

using the support parametrization

$$x_k(\theta)=h_k(\theta)(\cos\theta,\sin\theta)+h_k'(\theta)(-\sin\theta,\cos\theta).$$

Here the subscripts 2 and 4 denote the Fourier frequencies in Proposition 9.1. Both bodies are smooth, strictly convex, and centrally symmetric, and they satisfy

$$L(P_2)=L(P_4)=2\pi,$$

$$Area(P_2)=Area(P_4)=\frac{19493\pi}{20000}.$$

Their width laws differ, since

$$E\,w_{P_2}(U_\circ)^2-E\,w_{P_4}(U_\circ)^2=\frac{169}{6250}>0.$$

The curves are generated from the exact support parametrization and are not fitted approximations.

## *9.2. Centrally symmetric spatial separation*

We now use the first two even modes in Proposition 9.1. Fix

$$0<t<\frac{1}{3\sqrt{5}}, H>0,$$

and let $P_1, P_2\subset\mathbb{R}^2$ have support functions

$$h_1(\theta)=1+t\cos 2\theta, h_2(\theta)=1+\frac{t}{\sqrt{5}}\cos 4\theta.$$

Set

$$K_j=P_j\times[-H/2,H/2], j=1,2. \quad (58)$$

By Proposition 9.1 the planar bodies are smooth, strictly convex, and centrally symmetric, with common perimeter and area

$$L=2\pi, A=\pi\left(1-\frac{3}{2}t^2\right).$$

The prisms $K_1, K_2$ are centrally symmetric convex bodies. Their flat bases and vertical generators are retained in the construction; no three-dimensional smoothness or strict convexity is asserted.

**Theorem 9.2 (intrinsic volumes do not determine width or brightness laws).** The bodies $K_1, K_2\subset\mathbb{R}^3$ defined in (58) satisfy

$$V_3(K_j)=AH, V_2(K_j)=A+\pi H, V_1(K_j)=\pi+H, \quad (59)$$

and hence have identical intrinsic volumes. Nevertheless,

$$E W_{K_1}^2 - E W_{K_2}^2 = \frac{16}{15} t^2 > 0 , E B_{K_1}^2 - E B_{K_2}^2 = \frac{16}{15} H^2 t^2 > 0 . (60)$$

In particular, the intrinsic-volume nonclosure exhibited in [1, Theorem 7.5] persists within the centrally symmetric class, and neither the random-direction width law nor the random-direction brightness law is determined by $V_1, V_2, V_3$ there.

*Proof.* The prism volume is $AH$, and its surface area is $2A + LH$, so the first two identities in (59) follow from the Steiner normalization. To obtain the third, write a uniform direction $U \in S^2$ as

$$U = \left( \sqrt{1 - Z^2} U_\circ , Z \right),$$

where $Z$ is uniform on $[-1, 1]$, $U_\circ$ is uniform on $S^1$, and the two variables are independent. Minkowski additivity of width gives the pointwise identity

$$W_{K_j} = \sqrt{1 - Z^2} w_{P_j}\left(U_\circ\right) + H |Z| . (61)$$

The planar Cauchy formula gives $E w_{P_j}\left(U_\circ\right) = L / \pi$, while direct integration gives

$$E \sqrt{1 - Z^2} = \frac{\pi}{4} , E|Z| = \frac{1}{2} .$$

Thus $E W_{K_j} = L/4 + H/2$. Since $V_1 = 2 E W$ in the normalization (55), $V_1\left(K_j\right) = L/2 + H = \pi + H$, proving (59).

For the second moment in (61), the required one-variable averages are

$$E\left(1 - Z^2\right) = \frac{2}{3} , E\left(|Z| \sqrt{1 - Z^2}\right) = \frac{1}{3} , E Z^2 = \frac{1}{3} . (62)$$

Proposition 9.1 gives

$$E w_{P_1}\left(U_\circ\right)^2 - E w_{P_2}\left(U_\circ\right)^2 = \frac{8}{5} t^2 .$$

The mixed and vertical terms in the square of (61) are equal for the two bodies because their planar mean widths and heights agree. Hence

$$E W_{K_1}^2 - E W_{K_2}^2 = \frac{2}{3} \cdot \frac{8}{5} t^2 = \frac{16}{15} t^2 .$$

It remains to separate brightness. Cauchy’s surface-area formula applied to a right prism gives, pointwise,

$$B_{K_j} = A |Z| + H \sqrt{1 - Z^2} w_{P_j}\left(U_\circ^\perp\right), (63)$$

where $U_\circ^\perp$ is obtained from $U_\circ$ by a rotation through $\pi / 2$. Indeed, the two horizontal bases contribute $A|Z|$. On the lateral surface, whose outward normals are $(v, 0)$, the same formula reduces the contribution to

$$\frac{H \sqrt{1 - Z^2}}{2} \int_{\partial P_j} \left|U_\circ \cdot v\right| ds .$$

The planar Cauchy formula identifies the last half-integral with the projection length of $P_j$ onto $U_\circ^\perp$, namely $w_{P_j}\left(U_\circ^\perp\right)$. This proves (63) from the surface-area representation.

The area terms and the mixed terms in the square of (63) agree because $A$, $L$, and $H$ agree, and $U_\circ^\perp$ is again uniform on $S^1$. Equations (57) and (62) therefore give

$$E\,B_{K_1}^2 - E\,B_{K_2}^2 = H^2 \frac{2}{3} \cdot \frac{8}{5} t^2 = \frac{16}{15} H^2 t^2.$$

Both gaps are strictly positive, so the corresponding laws are different. ▫

### *9.3. Harmonic and geometric information boundaries*

Theorem 9.2 strengthens the nonclosure phenomenon established in [1, Theorem 7.5] in two directions: both witnesses are centrally symmetric, and the same pair separates the brightness laws as well as the width laws. The construction also exposes the mechanism. Intrinsic-volume equality fixes the constant support mode and one curvature-weighted quadratic combination of the planar harmonic amplitudes, whereas the width and brightness second moments retain a different quadratic combination. The factor $1/\sqrt{5}$ in the fourth harmonic is exactly what preserves area while changing that second-moment energy.

From the perspective of Section 6, the two planar generating patterns have equal low-order scalar geometry but different second-order directional energy after the prism lift. From the perspective of Section 7, the failure is caused by discarding directional labels before reconstruction: three intrinsic volumes cannot replace the pair-triple correlation information used by finite normal tomography. Thus the positive closure within rectangular boxes and the negative conclusion here are compatible. The box class has enough algebraic rigidity for three scalars to recover its finite parameter multiset; centrally symmetric convex bodies in general do not.

The witnesses in Theorem 9.2 are deliberately elementary right prisms, so the proof is exact and entirely by hand. Whether the same simultaneous separation can be obtained by smooth strictly convex three-dimensional bodies while preserving all three intrinsic volumes remains open here.

## 10. Folded-normal identifiability and the hierarchy of directional information

The coordinatewise absolute value of a centered Gaussian vector is a centered multivariate folded-normal vector. The inverse problem considered here is to determine which finite labelled product moments of that vector identify its correlation matrix, and to specify the unavoidable ambiguity. Under an explicit condition on the correlation graph, second- and third-order moments suffice in every dimension. For a complete correlation graph, the reconstruction uses $(m-1)^2$ labelled observations and no moment of order four or higher.

Classical work of Nabeya [11, 12] gives explicit low-order absolute moments of Gaussian vectors. More recently, Kan and Robotti [24] developed recurrences for arbitrary folded and truncated multivariate-normal moments, Liu, Jin, Yang and Pan [25] studied distributional properties and parameter estimation for the multivariate folded normal law, and Benko, Hübnerová and Witkovský [26] derived its characteristic and moment-generating functions. These works concern forward evaluation or statistical estimation. The reconstruction theorem below addresses the inverse problem. Its switching mechanism also clarifies the relation between labelled Gaussian moments, the finite normal tomography of Section 7, direction-labelled cosine transforms, and scalar width and brightness laws.

### *10.1. Pair and triple absolute moments*

Let

$$X=(X_1,\ldots,X_m)$$

be a centered Gaussian vector with unit coordinate variances and correlation matrix

$$\Sigma=(\rho_{ij})_{i,j=1}^{m}.$$

Its correlation graph $\Gamma$ has vertex set $\{1,\ldots,m\}$ and edge $ij$ whenever $\rho_{ij}\neq 0$. As in Section 7, let $C(\Gamma)$ be its cycle space over $F_2$, and let $T(\Gamma)$ be the subspace generated by its triangles.

For $0\le r\le 1$, put

$$g_2(r)=\frac{2}{\pi}\left(\sqrt{1-r^2}+r\arcsin r\right).$$

The classical bivariate formula and its derivative are

$$E|X_iX_j|=g_2(|\rho_{ij}|),\ g_2'(r)=\frac{2}{\pi}\arcsin r.\ (64)$$

Thus the labelled pair moments recover every $|\rho_{ij}|$, including the zero pattern. To see what is added by a triple moment, fix $a,b,c\in(0,1)$ and write

$$\Sigma_\pm(a,b,c)=\begin{pmatrix}1 & a & b\\ a & 1 & \pm c\\ b & \pm c & 1\end{pmatrix},\ \Delta_\pm=1-a^2-b^2-c^2\pm 2abc.$$

Whenever the indicated matrix is positive semidefinite, let $G_\pm(a,b,c)$ denote the corresponding value of $E|X_1X_2X_3|$. The Gaussian form of Theorem 7.2 is

$$G_+(a,b,c)-G_-(a,b,c)=\frac{2\sqrt{2}}{\pi^{3/2}}\int_0^a\int_0^b\int_{-c}^{c}\frac{dz\,dy\,dx}{\sqrt{1-x^2-y^2-z^2+2xyz}},\ (65)$$

when the negative class is feasible, with the rank-two boundary understood as an improper integral. In particular,

$$G_+(a,b,c)-G_-(a,b,c)\ge\frac{4\sqrt{2}}{\pi^{3/2}}abc>0.\ (66)$$

The strict inequality, rather than merely the availability of a closed formula, is the decisive inverse fact: it makes the triangle sign observable from one labelled third absolute moment.

### 10.2. Gaussian absolute-moment tomography

For $\varepsilon = (\varepsilon_1, \ldots, \varepsilon_m) \in \{\pm 1\}^m$, write

$$D_\varepsilon = diag(\varepsilon_1, \ldots, \varepsilon_m).$$

**Theorem 10.1 (Gaussian absolute-moment tomography).** Let $X$ be a centered Gaussian vector with unit coordinate variances, correlation matrix $\Sigma = (\rho_{ij})$, and $|\rho_{ij}| < 1$ for $i \neq j$. If

$$T(\Gamma) = C(\Gamma),$$

then the labelled moments

$$\{E|X_i X_j| : 1 \le i < j \le m\}$$

and

$$\{E|X_i X_j X_k| : 1 \le i < j < k \le m\}$$

determine $\Sigma$ up to diagonal sign conjugacy

$$\Sigma \mapsto D_\varepsilon \Sigma D_\varepsilon.$$

If $\Gamma$ is complete, it is enough to use the pair moments together with the anchored triple moments

$$\{E|X_1 X_i X_j| : 2 \le i < j \le m\}.$$

Thus, for a complete correlation graph, the displayed reconstruction uses

$$\binom{m}{2} + \binom{m-1}{2} = (m-1)^2$$

nontrivial labelled observations. No rank restriction on $\Sigma$ is required, and no absolute product moment of order four or higher is used. Equivalently, under the stated graph condition, the correlation matrix of the standard centered multivariate folded-normal vector is identifiable from its labelled second- and third-order product moments up to diagonal sign conjugacy.

*Proof.* Equation (64) and the strict monotonicity of $g_2$ recover every $|\rho_{ij}|$ and hence $\Gamma$. Consider a triangle $ijk$ in $\Gamma$, and put

$$a = |\rho_{ij}|, b = |\rho_{ik}|, c = |\rho_{jk}|.$$

If $\Delta_- < 0$, the negative switching class is not positive semidefinite, so the triangle sign is positive. If $\Delta_- \ge 0$, both switching classes are feasible and (65)–(66) show that their triple absolute moments are distinct. The observed value of $E|X_i X_j X_k|$ therefore determines

$$sgn(\rho_{ij} \rho_{ik} \rho_{jk}).$$

The triangle moments consequently recover the sign of every triangle cycle. Since the triangle cycles generate $C(\Gamma)$, they recover every cycle sign. Proposition 7.1 then determines the signing of $\Sigma$ up to diagonal sign conjugacy. A centered Gaussian law is determined by its covariance matrix, including when that matrix is singular.

If $\Gamma = K_m$, fix the gauge by taking $\rho_{1j} > 0$ for $j > 1$. The sign of $\rho_{ij}$ then equals the sign of $\rho_{1i}\rho_{1j}\rho_{ij}$, which is recovered from the anchored triple $(1, i, j)$. The anchored triangles

$$1ij, 2 \le i < j \le m,$$

form a basis of $C(K_m)$. Indeed, they are linearly independent because each non-anchor edge $ij$ occurs in exactly one such triangle, and their number is

$$\binom{m-1}{2} = \binom{m}{2} - m + 1 = dim\, C(K_m).$$

Thus the anchored family is minimal among triangle families spanning the abstract cycle space. This is a minimality statement about the displayed cycle-basis construction, not a claim of global optimality among arbitrary real-valued observation schemes. Finally, the remaining ambiguity is unavoidable: if $Y = D_\varepsilon X$, then $Y$ has covariance $D_\varepsilon \Sigma D_\varepsilon$, while $|Y_i| = |X_i|$ for every coordinate. ▫

The theorem is dimension-free because its proof uses only two- and three-coordinate Gaussian marginals and the cycle structure of the correlation graph. It is therefore not a reformulation of the rank-three directional reconstruction in Theorem 7.3, although both results are governed by the same switching invariants.

The quantitative separation is also inherited. If $a, \tilde{a} \in [\eta, 1]$, then (64) and the mean-value theorem give

$$|a - \tilde{a}| \le \frac{\pi}{2 \arcsin \eta} |g_2(a) - g_2(\tilde{a})|.$$

If the three edge magnitudes of a feasible triangle are at least $\eta$, equation (66) separates its two Gaussian branches by at least

$$\frac{4\sqrt{2}}{\pi^{3/2}} \eta^3.$$

Thus the reconstruction is stable on fixed correlation strata bounded away from zero. As in Proposition 7.4, triangle generation is asserted as a sufficient graph condition. Positive-semidefinite constraints may force additional signs when $T(\Gamma)$ is smaller than $C(\Gamma)$, and no necessity claim is made in that regime.

### *10.3. A four-cycle obstruction at third order*

The triangle-span condition cannot be omitted from a general dimension-free reconstruction theorem based only on pair and triple absolute moments. The following positive-definite example already occurs on four coordinates.

**Proposition 10.2 (four-cycle obstruction at third order).** Put $r = 2/5$ and

$$\Sigma_+^{\square} = \begin{pmatrix} 1 & r & 0 & r \\ r & 1 & r & 0 \\ 0 & r & 1 & r \\ r & 0 & r & 1 \end{pmatrix}, \Sigma_-^{\square} = \begin{pmatrix} 1 & r & 0 & -r \\ r & 1 & r & 0 \\ 0 & r & 1 & r \\ -r & 0 & r & 1 \end{pmatrix}.$$

Both matrices are positive definite. They are not diagonally sign-conjugate, but every corresponding principal submatrix of order at most three is diagonally sign-conjugate. Consequently, if $X^{+}$ and $X^{-}$ are centered Gaussian vectors with correlation matrices $\Sigma_{+}^{\square}$ and $\Sigma_{-}^{\square}$, respectively, then all labelled absolute product moments involving at most three distinct coordinates agree, although $\Sigma_{+}^{\square}$ and $\Sigma_{-}^{\square}$ belong to different sign-gauge classes.

*Proof.* Write $\Sigma_{\pm}^{\square}=I+r\,A_{\pm}^{\square}$, where $A_{\pm}^{\square}$ are the signed adjacency matrices of the four-cycle. For the positive signing, the vectors

$$(1,1,1,1),(1,-1,1,-1),(1,0,-1,0),(0,1,0,-1)$$

give the eigenvalues $2,-2,0,0$ of $A_{+}^{\square}$. Hence the eigenvalues of $\Sigma_{+}^{\square}$ are

$$1+2r,1-2r,1,1,$$

which are positive for $r=2/5$. Direct multiplication gives

$$\left(A_{-}^{\square}\right)^{2}=2I.$$

Since $A_{-}^{\square}$ is real symmetric with trace zero, its eigenvalues are $\sqrt{2},\sqrt{2},-\sqrt{2},-\sqrt{2}$. The eigenvalues of $\Sigma_{-}^{\square}$ are therefore $1+r\sqrt{2}$ and $1-r\sqrt{2}$, each with multiplicity two, and are again positive.

The correlation graph is the four-cycle. Its unique cycle has sign $+1$ for $\Sigma_{+}^{\square}$ and $-1$ for $\Sigma_{-}^{\square}$, so Proposition 7.1 shows that the two matrices are not diagonally sign-conjugate. Every induced subgraph on at most three vertices is a forest. Every signing of a forest is switching-equivalent to the all-positive signing, again by Proposition 7.1. Thus the corresponding principal submatrices of $\Sigma_{+}^{\square}$ and $\Sigma_{-}^{\square}$ are diagonally sign-conjugate. The associated Gaussian marginals are therefore related by coordinate sign changes, which leave their coordinatewise absolute values unchanged. Every labelled absolute product moment supported on such a marginal is consequently the same for $X^{+}$ and $X^{-}$. $\square$

Both matrices in Proposition 10.2 have rank four. Consequently, this example does not resolve whether the four-cycle ambiguity survives under the rank-three positive-semidefinite Gram constraint for direction families in $\mathbb{R}^{3}$. The corresponding directional reconstruction problem remains open.

This example leaves a natural higher-order question. For $l\geq 4$, do labelled absolute product moments involving at most $l$ distinct coordinates determine the sign of every cycle of length at most $l$, and hence recover the correlation matrix up to diagonal sign conjugacy whenever those cycles generate $C(\Gamma)$? Theorem 10.1 answers the corresponding question at order three.

### *10.4. Inversion of direction-labelled cosine transforms*

The Gaussian theorem retains coordinate labels. The analogous fully direction-labelled object in convex geometry is the cosine transform itself. For a finite signed Borel measure $\nu$ on $S^{2}$, let $\nu_{e}$ denote its even part, as in Definition 2.3.

The injectivity of the spherical cosine transform on even measures is the measure form of the classical projection theorem of Aleksandrov; see [9, Theorem 3.3.6; 13, Section 3.4]. The classical result is dimension-independent. We state the three-dimensional form required here and include the short harmonic proof because Theorem 6.12 already supplies the normalized multipliers.

**Proposition 10.3 (inversion of the labelled cosine transform).** For finite signed Borel measures $\nu$ and $\eta$ on $S^2$,

$$A_\nu = A_\eta \Leftrightarrow \nu_e = \eta_e . (67)$$

*Proof.* The reverse implication follows from the evenness of the kernel. For the forward implication, put $\theta = \nu_e - \eta_e$. Then $\theta$ is even and $A_\theta = 0$. The harmonic expansion in the proof of Theorem 6.12 gives

$$\widehat{A}_{\theta l m} = a_l \hat{\theta}_{l m} .$$

All odd coefficients of $\theta$ vanish by parity, while every even multiplier $a_{2j}$, including $a_0$, is nonzero by (39). Hence every spherical-harmonic coefficient of $\theta$ vanishes. Since finite linear combinations of spherical harmonics are uniformly dense in $C(S^2)$, the measure $\theta$ is zero. ▫

Consequently, the direction-labelled width function of a zonoid determines its unique even generating measure and the zonoid up to translation. A direction-labelled brightness function determines the even part of the surface-area measure. If $K$ is centrally symmetric and has nonempty interior, then $S_K$ is even, so Minkowski uniqueness determines $K$ up to translation. These conclusions have sharp limitations: for an arbitrary convex body, the width function determines the difference body $K - K$, not $K$, and the brightness function of a general body determines only $(S_K)_e$.

### 10.5. Scalar laws and moment contractions

Let $\nu$ be a finite positive Borel measure on $S^2$, and write

$$\lambda_\nu = (A_\nu)_\# \sigma$$

for the scalar law obtained after the direction label has been discarded. For $k \geq 1$, its moments are

$$m_k(\nu) = \int_{\mathbb{R}} t^k \, d\lambda_\nu(t) = \int_{S^2} A_\nu(u)^k \, d\sigma(u) = P_k(\nu, \dots, \nu) . (68)$$

Because $0 \leq A_\nu \leq \nu(S^2)$, the law $\lambda_\nu$ has compact support and is determined by the complete sequence $(m_k(\nu))_{k \geq 0}$.

For a finite even direction system

$$\nu = \frac{1}{2} \sum_{i=1}^{N} c_i (\delta_{n_i} + \delta_{-n_i}), c_i > 0 ,$$

define

$$T^{(k)}_{i_1, \dots, i_k} = I_k(n_{i_1}, \dots, n_{i_k}) .$$

Then

$$m_k(v) = \sum_{i_1,\dots,i_k=1}^{N} c_{i_1} \cdots c_{i_k} T^{(k)}_{i_1,\dots,i_k}. \quad (69)$$

Thus an ordinary scalar moment is a complete contraction of a labelled tensor. The contraction retains enough information to determine the one-dimensional law when all orders are known, but it generally does not recover the labels or the tensor entries from which it was formed. Theorem 9.2 and Theorem 10.1 therefore occupy different information levels: the former proves loss after scalarization to classical invariants, whereas the latter proves finite recovery when pair and triple labels are retained.

### *10.6. The information hierarchy*

The preceding results give the following hierarchy:

$$\begin{array}{c} \text{direction-labelled } A_v \\ \Updownarrow \\ v_e \\ \Downarrow \\ \text{labelled polarized data} \\ \Downarrow \\ \lambda_v \\ \Updownarrow \\ \left(m_k(v)\right)_{k\ge 0} \\ \Downarrow \\ \text{finitely many ordinary moments}. \end{array} \quad (70)$$

The first equivalence is Proposition 10.3, and the second is compact-support moment determinacy. The downward arrows are unconditional passages obtained by evaluation, contraction, or forgetting labels. Their reverses require additional structure and must be supplied by separate theorems.

For finite direction systems in $\mathbb{R}^3$, Theorem 7.3 supplies one such reverse implication:

$$\left.\begin{array}{c} \text{labelled } I_2 \text{ and } I_3 \text{ data} \\ T(\Gamma)=C(\Gamma) \end{array}\right\} \Rightarrow \text{direction configuration modulo sign and } O(3). \quad (71)$$

For Gaussian vectors of arbitrary dimension, Theorem 10.1 gives the parallel statement

$$\left.\begin{array}{c} \text{labelled pair and triple absolute moments} \\ T(\Gamma)=C(\Gamma) \end{array}\right\} \Rightarrow \text{correlation matrix modulo diagonal sign conjugacy}. \quad (72)$$

In the opposite direction, Theorem 9.2 shows that even the complete intrinsic-volume vector does not determine the width or brightness law among centrally symmetric convex bodies. The four levels must therefore remain distinct. A direction-labelled transform can determine an even measure; low-order labelled data can reconstruct finite configurations under an explicit graph condition; a scalar law retains every ordinary moment but forgets where its values occur; and finitely many classical scalar invariants retain still less. This separation is the inverse counterpart of the exact forward density theory developed in Sections 3–8.

# Appendix A. Hand evaluation of the tetrahedral width checks

The following calculations supply independent hand checks of the normalization and first two moments used in Sections 4.5–4.6.

### A.1. Recovery of the mean width

Setting $r=1$ in (11), with $g(\theta)=\sqrt{2}\cos\theta+\sin\theta=\sqrt{3}\cos(\theta-\alpha)$ and $\alpha=\arctan(1/\sqrt{2})$, gives

$$E\,W_T=\frac{6}{\pi}\int_0^{\pi/2}\frac{d\theta}{1+g(\theta)^2}.$$

Put $x=\theta-\alpha$ and $y=\tan x$. The endpoints become $-1/\sqrt{2}$ and $\sqrt{2}$, while $dx=dy/(1+y^2)$ and $1+3\cos^2 x=(y^2+4)/(1+y^2)$. Hence

$$E\,W_T \quad =\frac{6}{\pi}\int_{-1/\sqrt{2}}^{\sqrt{2}}\frac{dy}{y^2+4}$$

Since $\cos(2\arctan\sqrt{2})=-1/3$ and $2\arctan\sqrt{2}\in(0,\pi)$,

$$E\,W_T=\frac{3}{2\pi}\arccos\left(-\frac{1}{3}\right),$$

which is Finch’s unit-edge mean width [8].

### A.2. Recovery of the mean-square width

For $r=2$, equation (11) gives

$$E\,W_T^2=\frac{4}{\pi}\int_0^{\pi/2}\left(1-t_{min}(\theta)^3\right)d\theta\,,\,t_{min}(\theta)=\frac{g(\theta)}{\sqrt{1+g(\theta)^2}}.$$

With $x=\theta-\alpha$ and $y=\sin x$, the endpoints are $-1/\sqrt{3}$ and $\sqrt{2/3}$, and

$$\int_0^{\pi/2}t_{min}(\theta)^3\,d\theta=3\sqrt{3}\int_{-1/\sqrt{3}}^{\sqrt{2/3}}\frac{1-y^2}{(4-3y^2)^{3/2}}\,dy.$$

An antiderivative of the displayed integrand, including the prefactor $3\sqrt{3}$, is

$$H(y)=-\frac{\sqrt{3}}{4}\frac{y}{\sqrt{4-3y^2}}+\arcsin\left(\frac{\sqrt{3}\,y}{2}\right).$$

Direct evaluation gives

$$H\left(\sqrt{\frac{2}{3}}\right)=-\frac{1}{4}+\frac{\pi}{4}\,,\,H\left(-\frac{1}{\sqrt{3}}\right)=\frac{\sqrt{3}}{12}-\frac{\pi}{6}\,,$$

and therefore

$$E W_T^2 \quad = \frac{4}{\pi}\left(\frac{\pi}{2} - \frac{5\pi}{12} + \frac{1}{4} + \frac{\sqrt{3}}{12}\right)$$

This is Finch's unit-edge mean-square width [8].

### A.3. Direct normalization of the density

For $t_0 \le t \le t_2$, let

$$G(t) = t\gamma(t) - \arctan\sqrt{3 - 4t^2}, G'(t) = \gamma(t).$$

The endpoint values are

$$G(t_0) = t_0\left(\frac{\pi}{2} - \alpha\right) - \frac{\pi}{4}, G(t_1) = t_1\alpha - \frac{\pi}{6}, G(t_2) = 0.$$

Integrating the three branches of (8) and substituting these values gives

$$\begin{aligned} \int_{t_0}^{1} f_T(t)\,dt \quad &= 6(1 - t_0) - \frac{12\alpha}{\pi}(t_1 - t_0) + \frac{12}{\pi}\left(G(t_0) + G(t_1)\right) \\ &= 1. \end{aligned}$$

Thus (8) is normalized directly, independently of the geometric normalization supplied by Theorem 3.1.

## Appendix B. Auxiliary calculations for the tetrahedral brightness law

### B.1. Antiderivatives, normalization, and the CDF

For

$$\eta(t) = \arcsin\left(\frac{\sqrt{2}t}{\sqrt{1 - 4t^2}}\right), \delta(t) = \arccos\left(\frac{\sqrt{2}t}{\sqrt{3 - 16t^2}}\right),$$

differentiation followed by integration by parts and the substitution $v = \sqrt{1 - 6t^2}$ gives the antiderivatives used in (22):

$$\begin{aligned} J_{max}(t) \quad &= \int \eta(t)\,dt = t\eta(t) + \frac{1}{2}\arctan\sqrt{2(1 - 6t^2)}, \\ J_{\Sigma}(t) \quad &= \int \delta(t)\,dt = t\delta(t) - \frac{\sqrt{3}}{4}\arctan\left(2\sqrt{2(1 - 6t^2)}\right). \end{aligned}$$

Indeed,

$$\eta'(t) = \frac{\sqrt{2}}{(1 - 4t^2)\sqrt{1 - 6t^2}}, \delta'(t) = -\frac{\sqrt{6}}{(3 - 16t^2)\sqrt{1 - 6t^2}},$$

while $1 - 4t^2 = (1 + 2v^2)/3$ and $3 - 16t^2 = (1 + 8v^2)/3$. The endpoint values from Lemma 5.3 yield

$$J_{max}(\tau_0) = \frac{\pi\tau_0}{4} + \frac{1}{2}\arctan\frac{1}{\sqrt{2}}, \quad J_{max}(\tau_1) = \frac{\pi\tau_1}{2},$$
$$J_\Sigma(\tau_0) = \frac{\pi\tau_0}{3} - \frac{\sqrt{3}}{4}\arctan\sqrt{2}, \quad J_\Sigma(\tau_1) = 0.$$

Using $\arctan(1/\sqrt{2}) + \arctan\sqrt{2} = \pi/2$, the first branch has mass

$$B_0 = \left(6 + \frac{16}{\sqrt{3}}\right)\tau_1 - 6.$$

With $c_* = 6 + 16/\sqrt{3}$, the other two masses are $B_1 = c_*(\tau_2 - \tau_1)$ and $B_2 = 6(\tau_3 - \tau_2)$; hence

$$B_0 + B_1 + B_2 = c_*\tau_2 - 6 + 6\tau_3 - 6\tau_2 = 1.$$

The same antiderivatives give Proposition 5.6. In particular, $F_{B_T}(\tau_1) = c_*\tau_1 - 6$, and the constant subsequent density branches give

$$F_{B_T}(t) = c_* t - 6\,(\tau_1 \le t \le \tau_2),\, F_{B_T}(t) = 6t - 2\,(\tau_2 \le t \le \tau_3).$$

### *B.2. Direct two-direction integral and the second moment*

Let $\omega, \zeta \in S^2$, and write

$$\rho = |\omega \cdot \zeta|.$$

Because the integrand is unchanged when $\zeta$ is replaced by $-\zeta$, assume $0 \le \rho \le 1$ and put

$$\rho = \cos\vartheta,\, 0 \le \vartheta \le \frac{\pi}{2}.$$

Choose coordinates in which

$$\omega = (1, 0, 0),\, \zeta = (\cos\vartheta, \sin\vartheta, 0).$$

A uniform direction $U \in S^2$ may be parametrized as

$$U = \left(\sqrt{1-z^2}\cos\psi, \sqrt{1-z^2}\sin\psi, z\right),$$

with normalized area element

$$d\sigma = \frac{1}{4\pi}\, d\psi\, dz,\, -1 \le z \le 1,\, 0 \le \psi < 2\pi.$$

Hence

$$E(|U\cdot\omega||U\cdot\zeta|) = \frac{1}{4\pi}\int_{-1}^{1}(1-z^2)\,dz$$

The first integral is $4/3$. For the angular integral, the sign changes occur at the zeros of the two cosine factors. Splitting one period at those four points and integrating

$$\cos\psi\cos(\psi-\vartheta)=\frac{1}{2}[\cos\vartheta+\cos(2\psi-\vartheta)]$$

gives

$$\int_0^{2\pi}|\cos\psi\cos(\psi-\vartheta)|d\psi=2\left[\sin\vartheta+\left(\frac{\pi}{2}-\vartheta\right)\cos\vartheta\right].$$

Consequently,

$$E(|U\cdot\omega||U\cdot\zeta|)\quad=\frac{2}{3\pi}\left[\sin\vartheta+\left(\frac{\pi}{2}-\vartheta\right)\cos\vartheta\right]$$

For the tetrahedral normals, $\rho=1/3$, so

$$E\left(|U\cdot\hat{\xi}_i||U\cdot\hat{\xi}_j|\right)=\frac{2}{9\pi}\left(2\sqrt{2}+\arcsin\frac{1}{3}\right),i\neq j.$$

Now write

$$X_i=|U\cdot\hat{\xi}_i|.$$

From (12),

$$B_T=\frac{\sqrt{3}}{8}\sum_{i=1}^{4}X_i.$$

Since

$$E\,X_i^2=E\left(U\cdot\hat{\xi}_i\right)^2=\frac{1}{3},$$

we have

$$\begin{aligned}E\,B_T^2\quad&=\frac{3}{64}\left[\sum_{i=1}^{4}E\,X_i^2+2\sum_{1\le i<j\le 4}E\left(X_iX_j\right)\right]\\&=\frac{1}{16}+\frac{1}{8\pi}\left(2\sqrt{2}+\arcsin\frac{1}{3}\right),\end{aligned}$$

as asserted in (24).

## Appendix C. Tight frames and the degree-two width variance

This appendix proves the tight-frame statement used in Section 8.5 and evaluates the first nonzero harmonic contributions for the truncated-octahedral and tetrahedral direction systems. Throughout, the spherical harmonics are orthonormal in $L^2(S^2,\sigma)$, where $\sigma(S^2)=1$. In this normalization the addition theorem is

$$\sum_{m=-l}^{l} Y_{lm}(n)\overline{Y_{lm}(n')} = (2l+1) P_l(n\cdot n').$$

This convention accounts for the absence of a $1/(4\pi)$ factor below.

### *C.1. Atomic harmonic energy and irreducible symmetry*

Let $n_1,\dots,n_N \in S^2$, let $c_1,\dots,c_N > 0$, and put

$$\mu = \frac{1}{2}\sum_{p=1}^{N} c_p\left(\delta_{n_p} + \delta_{-n_p}\right).$$

Define the weighted frame operator and total weight by

$$M = \sum_{p=1}^{N} c_p n_p n_p^T, s = \sum_{p=1}^{N} c_p.$$

**Proposition C.1 (degree-two energy and tight frames).** For every even $l$,

$$\sum_{m=-l}^{l} \left|\hat{\mu}_{lm}\right|^2 = (2l+1) \sum_{p,q=1}^{N} c_p c_q P_l(n_p\cdot n_q).$$

At degree two,

$$\sum_{p,q=1}^{N} c_p c_q P_2(n_p\cdot n_q) \quad = \frac{1}{2}\left(3\ \| M \|_F^2 - s^2\right)$$

Consequently, the degree-two contribution to the variance in (40) is

$$\frac{15}{128}\left\|M - \frac{s}{3} I\right\|_F^2.$$

It vanishes if and only if the weighted directions form a unit-norm tight frame:

$$M = \frac{s}{3} I.$$

*Proof.* Since $l$ is even, the antipodal atoms have the same degree-$l$ harmonic values. Hence

$$\hat{\mu}_{lm} = \sum_{p=1}^{N} c_p \overline{Y_{lm}(n_p)}.$$

Expanding the squared norm and applying the addition theorem in the normalization fixed above gives

$$\sum_{m=-l}^{l} \left|\hat{\mu}_{lm}\right|^2 \quad = \sum_{p,q=1}^{N} c_p c_q \sum_{m=-l}^{l} \overline{Y_{lm}(n_p)} Y_{lm}(n_q)$$

Since $P_2(t) = (3t^2 - 1)/2$,

$$\sum_{p,q=1}^{N} c_p c_q P_2(n_p \cdot n_q) = \frac{1}{2}\left(3\sum_{p,q=1}^{N} c_p c_q (n_p \cdot n_q)^2 - s^2\right).$$

Moreover,

$$\sum_{p,q=1}^{N} c_p c_q (n_p \cdot n_q)^2 \quad = \sum_{p,q=1}^{N} c_p c_q tr\left(n_p n_p^T n_q n_q^T\right)$$

Because $tr\, M = s$,

$$\left\| M - \frac{s}{3} I \right\|_F^2 = \| M \|_F^2 - \frac{s^2}{3}.$$

This gives both forms of the degree-two energy. Finally, $a_2 = 1/8$, while the addition theorem contributes the factor $2l+1=5$. The degree-two term in (40) is therefore

$$a_2^2 \cdot 5 \cdot \frac{3}{2} \left\| M - \frac{s}{3} I \right\|_F^2 = \frac{15}{128} \left\| M - \frac{s}{3} I \right\|_F^2.$$

Its vanishing is equivalent to $M = (s/3) I$. ▫

**Corollary C.2 (irreducible symmetry criterion).** Suppose that the even weighted measure $\mu$ is invariant under a subgroup $G \subset O(3)$ whose action on $\mathbb{R}^3$ is irreducible. Then

$$M = \frac{s}{3} I.$$

Consequently, the weighted direction system is a unit-norm tight frame and its degree-two variance contribution vanishes.

*Proof.* Invariance of $\mu$ gives

$$gM g^T = M \,(g \in G),$$

and hence $gM = Mg$. Since $M$ is symmetric, $\mathbb{R}^3$ is the orthogonal sum of its eigenspaces. Each eigenspace is $G$-invariant because $M$ commutes with $G$. Irreducibility therefore permits only one eigenvalue, so $M = \lambda I$. Taking traces gives $3\lambda = s$. Proposition C.1 completes the proof. ▫

The corollary concerns symmetry-invariant weighted direction measures. It does not assert that the width of every symmetric convex body has an atomic generating measure of this form.

### *C.2. The truncated-octahedral and tetrahedral systems*

The full octahedral group preserves the six generator lines in (46), while the tetrahedral rotation group preserves the four tetrahedral normal lines of Section 5. Both actions on $\mathbb{R}^3$ are irreducible. For the octahedral action this follows from the coordinate permutations and sign changes. For the tetrahedral action, rotations about two nonparallel vertex axes admit no common invariant line; an invariant plane would have an invariant orthogonal line. Corollary C.2 therefore gives

$$\sum_{j=1}^{6} r_j r_j^T = 2I, \sum_{i=1}^{4} \hat{\xi}_i \hat{\xi}_i^T = \frac{4}{3} I$$

from symmetry and the respective traces $6$ and $4$.

For the truncated-octahedral system there are three orthogonal unordered pairs and twelve unordered pairs of absolute correlation $1/2$. Since its degree-two term vanishes, equation (40) becomes

$$Var\, W_O = \sum_{j\geq 2} (4j+1) a_{2j}^2 \left[ 6 + 6 P_{2j}(0) + 24 P_{2j}\left(\frac{1}{2}\right) \right].$$

For the tetrahedral normal system, every distinct pair has absolute correlation $1/3$, and hence

$$Var\, W_R = \sum_{j\geq 2} (4j+1) a_{2j}^2 \left[ 4 + 12 P_{2j}\left(\frac{1}{3}\right) \right].$$

The required Legendre values are

$$\begin{aligned} P_4(0) &= \frac{3}{8}, & P_4\left(\frac{1}{2}\right) &= -\frac{37}{128}, & P_4\left(\frac{1}{3}\right) &= \frac{1}{81}, \\ P_6(0) &= -\frac{5}{16}, & P_6\left(\frac{1}{2}\right) &= \frac{331}{1024}, & P_6\left(\frac{1}{3}\right) &= \frac{47}{243}. \end{aligned}$$

Together with

$$a_4 = -\frac{1}{48}, a_6 = \frac{1}{128},$$

these values give

| observable | degree four | degree six | exact total variance |
|---|---|---|---|
| $W_O$ | $\frac{21}{4096}$ | $\frac{19773}{2097152}$ | $-\frac{17}{3} + \frac{4+8\sqrt{3}}{\pi}$ |
| $W_R$ | $\frac{7}{432}$ | $\frac{13}{2592}$ | $\frac{8}{3\pi}\left(2\sqrt{2} + \arcsin\frac{1}{3}\right) - \frac{8}{3}$ |

The degree-four and degree-six terms account for approximately $84.6\,\%$ of $Var\, W_O$ and $93.7\,\%$ of $Var\, W_R$. Thus symmetry removes the leading degree-two contribution, while the next two even degrees contain most of the remaining dispersion in both examples.

Because

$$B_T = \frac{\sqrt{3}}{8} W_R,$$

each harmonic contribution to $Var\, B_T$ is the corresponding contribution for $W_R$ multiplied by $3/64$. This statement concerns the tetrahedral facet-normal generating measure and hence the zonotopal observable $B_T$, equivalently $W_R$. The non-zonotopal tetrahedral width $W_T$ is not represented by this measure and is not covered by the tight-frame generating-measure criterion.

**Acknowledgements**. The material in this paper has been in development for many years and appears in print only now. I thank my family for their patience over that time. I am grateful to the early readers for their careful reading and detailed suggestions, which materially improved the mathematical exposition, historical attribution, and organization of the paper. In particular, their comments led to a sharper treatment of the singular cases, the tight-frame interpretation of the degree-two harmonic term, and the positioning of the Gaussian reconstruction theorem; I owe a lasting debt to my late mathematics teacher, Amos Matalon, who showed me and seeded the love for advanced mathematics at an early age, insisted that it be taken seriously, and set me on the path to research. Any remaining errors are my own.

**Use of AI tools**. All derivations in this paper were carried out by hand. A large language model was used as an auxiliary tool in two respects: to run independent numerical checks confirming the closed forms stated here, and to assist with prior-art and bibliographic searching. Figures 1–4 are deterministic figures generated with Matplotlib directly from the formulas and exact geometric constructions established in this paper. The plotting script was drafted with the assistance of a large language model and independently checked by the author. No computational output is used as part of any proof. No figure in this paper is AI-generated imagery. The model produced no text and no mathematical argument in the paper, and the author takes full responsibility for its content.

**Declarations**. Funding: This research received no external funding. Competing interests: The author declares none. Data availability: The script generating Figures 1–4 is archived at 10.5281/zenodo.21870870. It is used for plotting and numerical confirmation only; no formula in this paper is derived from it.